\documentclass[a4paper,12pt]{article}
\usepackage[all]{xy}
\usepackage[T1]{fontenc}
\usepackage[dvips]{graphicx}
\usepackage{amsmath,amssymb,amsfonts,amsthm,latexsym,mathrsfs,textcomp,verbatim}
\xyoption{web}

\title{Lattices of theories}
\author{{Olivia Caramello} \vspace{3 mm}\\ {\small DPMMS, University of Cambridge,}\\{\small Wilberforce Road, Cambridge CB3 0WB, UK}\\{\small O.Caramello@dpmms.cam.ac.uk}}
\date{May 4, 2009}
\begin{document}
\bgroup           
\let\footnoterule\relax  
\maketitle
\flushleft  
\begin{abstract}
We show that there is a bijection between the subtoposes of the classifying topos of a geometric theory $\mathbb{T}$ over a signature $\Sigma$ and the closed geometric theories over $\Sigma$ which are `quotients' of the theory $\mathbb{T}$; next, we analyze how classical topos-theoretic constructions on the lattice of subtoposes of a given topos can be transferred, via the bijection above, to logical constructions in the corresponding lattice of theories.  
\end{abstract} 
\egroup 
\flushleft
\vspace{5 mm}

\begin{center}
Dedicated to Peter Johnstone and Martin Hyland\\ on the occasion of their sixtieth birthdays. 
\end{center}
\newpage
\tableofcontents
\newpage

\def\Monthnameof#1{\ifcase#1\or
   January\or February\or March\or April\or May\or June\or
   July\or August\or September\or October\or November\or December\fi}
\def\today{\number\day~\Monthnameof\month~\number\year}

%
%
%
\def\pushright#1{{
   \parfillskip=0pt            
   \widowpenalty=10000         
   \displaywidowpenalty=10000  
   \finalhyphendemerits=0      
  %
   \leavevmode                 
   \unskip                     
   \nobreak                    
   \hfil                       
   \penalty50                  
   \hskip.2em                  
   \null                       
   \hfill                      
   {#1}                        
  %
   \par}}                      

\def\qed{\pushright{$\square$}\penalty-700 \smallskip}

\newtheorem{theorem}{Theorem}[section]

\newtheorem{proposition}[theorem]{Proposition}

\newtheorem{scholium}[theorem]{Scholium}

\newtheorem{lemma}[theorem]{Lemma}

\newtheorem{corollary}[theorem]{Corollary}

\newtheorem{conjecture}[theorem]{Conjecture}

\newenvironment{proofs}%
 {\begin{trivlist}\item[]{\bf Proof }}%
 {\qed\end{trivlist}}

  \newtheorem{rmk}[theorem]{Remark}
\newenvironment{remark}{\begin{rmk}\em}{\end{rmk}}

  \newtheorem{rmks}[theorem]{Remarks}
\newenvironment{remarks}{\begin{rmks}\em}{\end{rmks}}

  \newtheorem{defn}[theorem]{Definition}
\newenvironment{definition}{\begin{defn}\em}{\end{defn}}

  \newtheorem{eg}[theorem]{Example}
\newenvironment{example}{\begin{eg}\em}{\end{eg}}

  \newtheorem{egs}[theorem]{Examples}
\newenvironment{examples}{\begin{egs}\em}{\end{egs}}


\mathcode`\<="4268  
\mathcode`\>="5269  
\mathcode`\.="313A  
\mathchardef\semicolon="603B 
\mathchardef\gt="313E
\mathchardef\lt="313C

\newcommand{\app}
 {{\sf app}}

\newcommand{\Ass}
 {{\bf Ass}}

\newcommand{\ASS}
 {{\mathbb A}{\sf ss}}

\newcommand{\Bb}
{\mathbb}

\newcommand{\biimp}
 {\!\Leftrightarrow\!}

\newcommand{\bim}
 {\rightarrowtail\kern-1em\twoheadrightarrow}

\newcommand{\bjg}
 {\mathrel{{\dashv}\,{\vdash}}}

\newcommand{\bstp}[3]
 {\mbox{$#1\! : #2 \bim #3$}}

\newcommand{\cat}
 {\!\mbox{\t{\ }}}

\newcommand{\cinf}
 {C^{\infty}}

\newcommand{\cinfrg}
 {\cinf\hy{\bf Rng}}

\newcommand{\cocomma}[2]
 {\mbox{$(#1\!\uparrow\!#2)$}}

\newcommand{\cod}
 {{\rm cod}}

\newcommand{\comma}[2]
 {\mbox{$(#1\!\downarrow\!#2)$}}

\newcommand{\comp}
 {\circ}

\newcommand{\cons}
 {{\sf cons}}

\newcommand{\Cont}
 {{\bf Cont}}

\newcommand{\ContE}
 {{\bf Cont}_{\cal E}}

\newcommand{\ContS}
 {{\bf Cont}_{\cal S}}

\newcommand{\cover}
 {-\!\!\triangleright\,}

\newcommand{\cstp}[3]
 {\mbox{$#1\! : #2 \cover #3$}}

\newcommand{\Dec}
 {{\rm Dec}}

\newcommand{\DEC}
 {{\mathbb D}{\sf ec}}

\newcommand{\den}[1]
 {[\![#1]\!]}

\newcommand{\Desc}
 {{\bf Desc}}

\newcommand{\dom}
 {{\rm dom}}

\newcommand{\Eff}
 {{\bf Eff}}

\newcommand{\EFF}
 {{\mathbb E}{\sf ff}}

\newcommand{\empstg}
 {[\,]}

\newcommand{\epi}
 {\twoheadrightarrow}

\newcommand{\estp}[3]
 {\mbox{$#1 \! : #2 \epi #3$}}

\newcommand{\ev}
 {{\rm ev}}

\newcommand{\Ext}
 {{\rm Ext}}

\newcommand{\fr}
 {\sf}

\newcommand{\fst}
 {{\sf fst}}

\newcommand{\fun}[2]
 {\mbox{$[#1\!\to\!#2]$}}

\newcommand{\funs}[2]
 {[#1\!\to\!#2]}

\newcommand{\Gl}
 {{\bf Gl}}

\newcommand{\hash}
 {\,\#\,}

\newcommand{\hy}
 {\mbox{-}}

\newcommand{\im}
 {{\rm im}}

\newcommand{\imp}
 {\!\Rightarrow\!}

\newcommand{\Ind}[1]
 {{\rm Ind}\hy #1}

\newcommand{\iten}[1]
{\item[{\rm (#1)}]}

\newcommand{\iter}
 {{\sf iter}}

\newcommand{\Kalg}
 {K\hy{\bf Alg}}

\newcommand{\llim}
 {{\mbox{$\lower.95ex\hbox{{\rm lim}}$}\atop{\scriptstyle
{\leftarrow}}}{}}

\newcommand{\llimd}
 {\lower0.37ex\hbox{$\pile{\lim \\ {\scriptstyle
\leftarrow}}$}{}}

\newcommand{\Mf}
 {{\bf Mf}}

\newcommand{\Mod}
 {{\bf Mod}}

\newcommand{\MOD}
{{\mathbb M}{\sf od}}

\newcommand{\mono}
 {\rightarrowtail}

\newcommand{\mor}
 {{\rm mor}}

\newcommand{\mstp}[3]
 {\mbox{$#1\! : #2 \mono #3$}}

\newcommand{\Mu}
 {{\rm M}}

\newcommand{\name}[1]
 {\mbox{$\ulcorner #1 \urcorner$}}

\newcommand{\names}[1]
 {\mbox{$\ulcorner$} #1 \mbox{$\urcorner$}}

\newcommand{\nml}
 {\triangleleft}

\newcommand{\ob}
 {{\rm ob}}

\newcommand{\op}
 {^{\rm op}}

\newcommand{\palrr}[4]{ 
  \def\labelstyle{\scriptstyle} 
  \xymatrix{ {#1} \ar@<0.5ex>[r]^{#2} \ar@<-0.5ex>[r]_{#3} & {#4} } } 
  
\newcommand{\palrl}[4]{ 
  \def\labelstyle{\scriptstyle} 
  \xymatrix{ {#1} \ar@<0.5ex>[r]^{#2}  &  \ar@<0.5ex>[l]^{#3} {#4} } } 
  
\newcommand{\pepi}
 {\rightharpoondown\kern-0.9em\rightharpoondown}

\newcommand{\pmap}
 {\rightharpoondown}

\newcommand{\Pos}
 {{\bf Pos}}

\newcommand{\prarr}
 {\rightrightarrows}

\newcommand{\princfil}[1]
 {\mbox{$\uparrow\!(#1)$}}

\newcommand{\princid}[1]
 {\mbox{$\downarrow\!(#1)$}}

\newcommand{\prstp}[3]
 {\mbox{$#1\! : #2 \prarr #3$}}

\newcommand{\pstp}[3]
 {\mbox{$#1\! : #2 \pmap #3$}}

\newcommand{\relarr}
 {\looparrowright}

\newcommand{\rlim}
 {{\mbox{$\lower.95ex\hbox{{\rm lim}}$}\atop{\scriptstyle
{\rightarrow}}}{}}

\newcommand{\rlimd}
 {\lower0.37ex\hbox{$\pile{\lim \\ {\scriptstyle
\rightarrow}}$}{}}

\newcommand{\rstp}[3]
 {\mbox{$#1\! : #2 \relarr #3$}}

\newcommand{\scn}
 {{\bf scn}}

\newcommand{\scnS}
 {{\bf scn}_{\cal S}}

\newcommand{\semid}
 {\rtimes}

\newcommand{\Sep}
 {{\bf Sep}}

\newcommand{\sep}
 {{\bf sep}}

\newcommand{\Set}
 {{\bf Set }}

\newcommand{\Sh}
 {{\bf Sh}}

\newcommand{\ShE}
 {{\bf Sh}_{\cal E}}

\newcommand{\ShS}
 {{\bf Sh}_{\cal S}}

\newcommand{\sh}
 {{\bf sh}}

\newcommand{\Simp}
 {{\bf \Delta}}

\newcommand{\snd}
 {{\sf snd}}

\newcommand{\stg}[1]
 {\vec{#1}}

\newcommand{\stp}[3]
 {\mbox{$#1\! : #2 \to #3$}}

\newcommand{\Sub}
 {{\rm Sub}}

\newcommand{\SUB}
 {{\mathbb S}{\sf ub}}

\newcommand{\tbel}
 {\prec\!\prec}

\newcommand{\tic}[2]
 {\mbox{$#1\!.\!#2$}}

\newcommand{\tp}
 {\!:}

\newcommand{\tps}
 {:}

\newcommand{\tsub}
 {\pile{\lower0.5ex\hbox{.} \\ -}}

\newcommand{\wavy}
 {\leadsto}

\newcommand{\wavydown}
 {\,{\mbox{\raise.2ex\hbox{\hbox{$\wr$}
\kern-.73em{\lower.5ex\hbox{$\scriptstyle{\vee}$}}}}}\,}

\newcommand{\wbel}
 {\lt\!\lt}

\newcommand{\wstp}[3]
 {\mbox{$#1\!: #2 \wavy #3$}}

\newcommand{\fu}[2]
{[#1,#2]}

\newcommand{\st}[2]
 {\mbox{$#1 \to #2$}}
\section{Introduction}

This paper provides a unification of the theory of elementary toposes with geometric logic, by passing through the theory of Grothendieck toposes.\\
The main ingredient of the paper is the duality theorem proved in section \ref{section_dualita}, which asserts the existence of a bijection between the subtoposes of the classifying topos of a given geometric theory $\mathbb T$ and the closed `quotients' of $\mathbb T$. In fact, the theorem allows us to interpret many concepts of elementary topos theory which apply to the lattice of subtoposes of a given topos at the level of geometric theories.\\
Notions that will be analyzed in the course of the paper include the coHeyting algebra structure on the lattice of subtoposes of a given topos, open, closed, quasi-closed subtoposes, the dense-closed factorization of a geometric inclusion, coherent subtoposes, subtoposes with enough points, the surjection-inclusion factorization of a geometric morphism, skeletal inclusions, atoms in the lattice of subtoposes of a given topos, Booleanization and DeMorganization of a topos.\\
Many results are established on the way. Specifically, section \ref{section_prooftheory} contains a proof-theoretic analysis of the notion of Grothendieck topology in view of the duality theorem, while section \ref{latticestructure} contains explicit descriptions of the Heyting operation between Grothendieck topologies on a given category and of the Grothendieck topology generated by a given collection of sieves; also, a number of results about the problem of `relativizing' a local operator with respect to a given subtopos are derived in section \ref{relat}.\\
We also provide applications of the duality theorem in various contexts; in particular, we discuss how the theorem can be used to shed light on axiomatization problems for geometric theories, and we prove a deduction theorem for geometric logic.\\
The final part of the paper is devoted to discussing the problem of characterizing the classifying toposes of theories presented as quotients of theories of presheaf type; here, we unify the `semantic' point of view of homogeneous models with respect to a given Grothendieck topology introduced in \cite{OC1} with the syntactic perspective provided by the duality theorem. In this context, we also derive a syntactic description of the finitely presented models of a cartesian theory.

\newpage
\section{Preliminary facts}

In this section we present some basic facts which will be useful for our analysis. All the terminology used in the course of the paper is borrowed from \cite{El} and \cite{El2}, if not otherwise stated.
\subsection{A 2-dimensional Yoneda Lemma}
An essential role in the present paper is played by a $2$-dimensional version of the Yoneda Lemma.\\
Recall that there are a number of $2$-categories which naturally play a role in topos theory; among them, there are certainly the $2$-category $\mathfrak{Cat}$ of small categories, functors and natural transformations between them and the $2$-category $\mathfrak{BTop}$ of Grothendieck toposes, geometric morphisms and geometric transformations between them. Also, we have all the $2$-categories arising from notable fragments of geometric logic, namely the $2$-category $\mathfrak{Cart}$ of cartesian categories, cartesian functors and natural transformations between them, the $2$-category $\mathfrak{Reg}$ of regular categories, regular functors and natural transformations between them, the $2$-category $\mathfrak{Coh}$ of coherent categories, coherent functors and natural transformations between them, and the $2$-category $\mathfrak{Geom}$ of geometric categories, geometric functors and natural transformations between them.\\
Given a strict $2$-category $\cal R$ and two $0$-cells $a$ and $b$ in $\cal R$, we say that $a$ and $b$ are \emph{equivalent} if there exists $1$-cells $f:a\to b$ and $g:b\to a$ and invertible $2$-cells $\alpha:f\circ g\imp 1_{b}$ and $\beta:g\circ f \imp 1_{a}$. 
Given a $2$-category $\cal R$, we have an obvious $2$-functor $Y:{\cal R} \to [{\cal R}^{\textrm{op}}, \mathfrak{Cat}]$ (where  and $[{\cal R}^{\textrm{op}}, \mathfrak{Cat}]$ is the $2$-category of $2$-functors ${\cal R}^{\textrm{op}} \to \mathfrak{Cat}$), which sends a $0$-cell $a$  to the (obviously defined) $2$-functor $Y(a):={\cal R}(-,a):{\cal R}^{\textrm{op}} \to \mathfrak{Cat}$. Notice that this notion of equivalence specializes in $\textbf{Cat}$ to the well-known notion of natural equivalence between small categories.\\
The following result is essentially the $2$-categorical equivalent of the fact that the Yoneda functor in $1$-category theory is faithful; it is probably folklore, but we present a proof for the reader's convenience.
\begin{proposition}
With the notation above, for any $a, b\in {\cal R}$, the functors $Y(a)$ and $Y(b)$ are equivalent (as $0$-cells in the $2$-category $[{\cal R}^{\textrm{op}}, \mathfrak{Cat}]$) (if and) only if $a$ and $b$ are equivalent (as $0$-cells in $\cal R$).\\
\end{proposition}   
\begin{proofs}
It is easy to see that two $2$-functors $F,G:{\cal R}^{\textrm{op}} \to \mathfrak{Cat}$ are equivalent if and only if for each $c\in {\cal R}$, the categories $F(c)$ and $G(c)$ are naturally equivalent via functors $K(c):F(c)\to G(c)$ and $L(c):G(c)\to F(c)$, naturally in $c\in {\cal R}$, i.e. for any $1$-cell $f:c\to d$ in $\cal R$ the obvious naturality squares for both $K$ and $L$ commute up to an invertible natural transformation.\\
Now suppose that for $a, b \in {\cal R}$ we have that $Y(a)$ and $Y(b)$ are equivalent via transformations $K:Y(a)\imp Y(b)$ and $L:Y(b)\imp Y(a)$ such that $K\circ L \cong Y(b)$ and $L\circ K \cong Y(a)$. Then we have $K(a):{\cal R}(a,a)\to {\cal R}(b,a)$ and $L(b):{\cal R}(b,b)\to {\cal R}(a,b)$; let us put $f:=K(a)(1_{a}):a\to b$ and $g:=L(b)(1_{b}):b\to a$. We want to prove that $g\circ f \simeq 1_{a}$ and $f\circ g \simeq 1_{b}$.\\
Consider the naturality square for $K$ corresponding to the arrow $g:b\to a$:
\[  
\xymatrix {
{\cal R}(a,a) \ar[r]^{K(a)} \ar[d]^{-\circ g}  & {\cal R}(a,b) \ar[d]^{-\circ g} \\
{\cal R}(b,a) \ar[r]^{K(b)} & {\cal R}(b,b) }
\] 
This square by our hypothesis commutes up to an invertible natural transformation, so $f\circ g=K(a)(1_{a})\circ g \cong K(b)(g)\cong K(b)(L(b)(1_{b}))\cong 1_{b}$. Dually, or more explicitly by replacing $K$ with $L$ and $f$ with $g$ in the argument above, one obtains the other isomorphism $g\circ f \cong 1_{a}$. So the $1$-cells $f$ and $g$ give an equivalence between $a$ and $b$, as required.
\end{proofs}

\subsection{An alternative view of Grothendieck topologies}
To begin, let us recall from \cite{MM} the definition of Grothendieck topology.\\
A Grothendieck topology on a category $\cal C$ is a function $J$ which assigns to each object $c$ of $\cal C$ a collection $J(c)$ of sieves on $c$ in such a way that\\
(maximality axiom) the maximal sieve $M_{c}=\{f \textrm{ | } cod(f)=c \}$ is in $J(c)$;\\
(stability axiom) if $S\in J(c)$, then $f^{\ast}(S)\in J(d)$ for any arrow $f:d\to c$;\\
(transitivity axiom) if $S\in J(c)$ and $R$ is any sieve on $c$ such that $f^{\ast}(R)\in J(d)$ for all $f:d\to c$ in $S$, then $R\in J(c)$.\\   
In a category $\cal C$ we call a collection of arrows in $\cal C$ with common codomain a \emph{presieve}; given a presieve $P$ on $c\in {\cal C}$, we define the sieve $\overline{P}$ generated by $P$ as the collection of all the arrows in $\cal C$ with codomain $c$ which factor through an arrow in $P$.\\
Given a collection $\cal U$ of presieves on $\cal C$, we define the Grothendieck topology generated by $\cal U$ to be the smallest Grothendieck topology $J$ on $\cal C$ such that all the sieves generated by the presieves in $\cal U$ are $J$-covering.\\
Given two Grothendieck topologies $J$ and $J'$ on a category $\cal C$ such that $J'\supseteq J$, we say that $J'$ is generated over $J$ by a collection $\cal U$ of sieves in $\cal C$ if $J'$ is generated by the collection of all the sieves on $\cal C$ which are either $J$-covering or belonging to $\cal U$.

\begin{rmk}\label{presieves}
\emph{Given a functor $F:{\cal C} \to {\cal E}$, where $\cal E$ is a Grothendieck topos, and a presieve $P$ in $\cal C$, $F$ sends $P$ to an epimorphic family if and only if it sends $\overline{P}$ to an epimorphic family; this remark will be useful below in connection with Diaconescu's theorem.}\\
\end{rmk}
We note that the definition of Grothendieck topology can also be put in the following alternative form.\\ 

\begin{definition} \label{def2}
A Grothendieck topology on a category $\cal C$ is a function $J$ which assigns to each object $c$ of $\cal C$ a collection $J(c)$ of sieves on $c$ in such a way that\\
(i) the maximal sieve $M_{c}$ belongs to $J(c)$;\\
(ii) for each pair of sieves $S$ and $T$ on $c$ such that $T\in J(c)$ and $S\supseteq T$, $S\in J(c)$;\\
(iii) if $R\in J(c)$ then for any arrow $g:d\rightarrow c$ there exists a sieve $S\in J(d)$ such that for each arrow $f$ in $S$, $g\circ f\in R$;\\
(iv) if the sieve $S$ generated by a presieve $\{f_{i}:c_{i}\rightarrow c \textrm{ | } i\in I \}$ belongs to $J(c)$ and for each $i\in I$ we have a presieve $\{g_{ij}:d_{ij}\rightarrow c_{i} \textrm{ | } j\in I_{i} \}$ such that the sieve $T_{i}$ generated by it belongs to $J(c_{i})$, then the sieve $R$ generated by the family of composites $\{f_{i}\circ g_{ij}:d_{ij}\rightarrow c \textrm{ | } i\in I, j\in I_{i} \}$ belongs to $J(c)$.\\
\end{definition}

In this definition, the sieve $R$ will be called the composite of the sieve $S$ with the sieves $T_{i}$ for $i\in I$ and denoted by $S\ast  \{T_{i} \textrm{ | } i\in I\}$.\\ 

Let us prove the equivalence of the two definitions. Let us assume the first definition and derive the second. To prove property (ii) let us assume that $S\supseteq T$ with $T\in J(c)$; then for every arrow $f$ in $T$ we have $f^{\ast}(S)\supseteq f^{\ast}(T)=M_{c}\in J(c)$ so by the transitivity axiom $S\in J(c)$, as required. Property (iii) immediately follows from the stability axiom. Property (iv) follows from the transitivity axiom for Grothendieck topologies by observing that for all arrows $f$ in $S$, $f^{\ast}(R)$ is $J$-covering. Indeed, if $f\in S$ then $f=f_{i}\circ h$ for some $i\in I$ and arrow $h$; so $f^{\ast}(R)=h^{\ast}(f_{i}^{\ast}(R))\supseteq h^{\ast}(T_{i})\in J(dom (h))$ and hence $f^{\ast}(R)\in J(dom(f))$ by property (ii) and the stability axiom.\\
Conversely, let us assume the second definition and derive the first. The stability axiom easily follows from (ii) and (iii); indeed, if $R\in J(c)$ and $g:d\rightarrow c$ is an arrow with codomain $c$, then $h^{\ast}(R)$ contains the sieve $S$ given by property (iii) and hence it is $J$-covering by property (ii). To prove the transitivity axiom we observe that, given a sieve $R$ on $c$ and a sieve $S\in J(c)$ such that for all arrows $h$ in $S$, $h^{\ast}(R)$ is $J$-covering, $R$ contains the composite of the sieve $S$ with the sieves of the form $h^{\ast}(R)$ for $h$ in $S$.\\
Note that, in Definition \ref{def2}, one can equivalently require in property (iv) that the presieves $\{f_{i}:c_{i}\rightarrow c \textrm{ | } i\in I \}$ and $\{g_{ij}:d_{ij}\rightarrow c_{i} \textrm{ | } j\in I_{i} \}$ are sieves; indeed, it is clear from the proof above that both versions of the condition are equivalent, under properties (i), (ii) and (iii), to the transitivity axiom.\\
Notice that, in the case the category $\cal C$ has pullbacks, property (iii) (equivalently, the stability axiom) may be replaced by the following condition: if (the sieve generated by) $\{f_{i}:c_{i}\rightarrow c \textrm{ | } i\in I \}$ belongs to $J(c)$ then for any arrow $g:d\rightarrow c$ the sieve generated by the family of pullbacks $\{ \textrm{p.b.($f_{i}$, $g$)} \rightarrow d \textrm{ | } i\in I \}$ belongs to $J(d)$.\\
\begin{rmk}\label{def2closure}
\emph{The operation of composition of sieves in a category $\cal C$ defined above behaves naturally with respect to the operator $\overline{(-)}^{J}$ of $J$-closure of sieves for a Grothendieck topology $J$ on $\cal C$; that is, with the notation above, we have $\overline{S\ast  \{T_{i} \textrm{ | } i\in I\}}^{J}=\overline{S\ast  \{\overline{T_{i}}^{J} \textrm{ | } i\in I\}}^{J}$. To verify this equality, it clearly suffices to prove that $S\ast  \{\overline{T_{i}}^{J}  \textrm{ | } i\in I\} \subseteq \overline{S\ast  \{T_{i} \textrm{ | } i\in I\}}^{J}$, and this easily follows from property (ii) in Definition \ref{def2}}.  
\end{rmk}

\subsection{Generators for Grothendieck topologies}
If $\cal C$ is a regular category, we may define the regular topology $J^{\textrm{reg}}_{\cal C}$ on $\cal C$ as the Grothendieck topology on $\cal C$ having as sieves exactly those which contain a cover. If $\cal C$ is a geometric category, we may define the geometric topology $J^{\textrm{geom}}_{\cal C}$ on it as the Grothendieck topology on $\cal C$ having as sieves exactly those which contain a small covering family. Notice that if ${\cal C}_{\mathbb T}$ is the geometric syntactic category of a geometric theory $\mathbb T$, then the geometric topology on ${\cal C}_{\mathbb T}$ concides with the syntactic topology $J_{\mathbb T}$ on ${\cal C}_{\mathbb T}$ (cfr. section \ref{section_dualita}).\\ The following result about these topologies hold. Below, by a principal sieve we mean a sieve which is generated by a single arrow. 

\begin{proposition}\label{gengrot}
Let $\cal C$ be a category and $J$ a Grothendieck topology on it. Then\\
(i) if $\cal C$ is regular and $J\supseteq J^{\textrm{reg}}_{\cal C}$ then $J$ is generated over $J^{\textrm{reg}}_{\cal C}$ by a collection of sieves generated by monomorphisms;\\
(ii) if $\cal C$ is geometric and $J\supseteq J^{\textrm{geom}}_{\cal C}$ then $J$ is generated over $J^{\textrm{geom}}_{\cal C}$  by a collection of principal sieves generated by a monomorphism.\\ 
\end{proposition}

\begin{proofs}
(i) Given an object $c\in {\cal C}$ and a sieve $R$ on $c$ in $\cal C$, let us denote, for each arrow $r$ in $R$, by $dom(r) \stackrel{r''}{\epi} x \stackrel{r'}{\mono} c$ its cover-mono factorization in $\cal C$ and by $R'$ the sieve in ${\cal C}$ generated by the arrows $r'$ (for $r$ in $R$). Clearly, it is enough to prove that $R\in J(c)$ if and only if $R'\in J(c)$. The `only if' part follows from property (ii) in Definition \ref{def2} since $R'\supseteq R$, while the `if' part follows from property (iv) in Definition \ref{def2} by using that, since $J\supseteq J^{\textrm{reg}}_{\cal C}$, all the sieves generated by the single arrows $r''$ (for $r\in R$) are $J$-covering.\\
(ii) Given an object $c\in {\cal C}$ and a sieve $R$ on $c\in {\cal C}$, let $r$ be the subobject of $c$ given by the union in $\Sub_{\cal C}(c)$ of all the images in $\cal C$ of the morphisms in $R$ (this union exists because, $\cal C$ being well-powered, there is, up to isomorphism, only a \emph{set} of monomorphisms with a given codomain). Clearly, it is enough to prove that $R\in J(c)$ if and only if $(r)\in J(c)$ (where $(r)$ denotes the sieve generated by the arrow $r$ in $\cal C$). The `only if' part follows from property (ii) in Definition \ref{def2} since $(r)\supseteq R$, while the `if' part follows from property (iv) in Definition \ref{def2}  by using that, since $J\supseteq J^{\textrm{geom}}_{\cal C}$, the sieve generated by the inclusions into $r$ of the images of the morphisms in $R$ is $J$-covering.\\   
\end{proofs}

Let us note that, given a sieve $R$ on a regular category $\cal C$, it is natural to consider the sieve $R^{\textrm{reg}}$ generated by the images of all the morphisms in $R$; similarly, if $\cal C$ is a geometric category, it is natural to consider the sieve $R^{\textrm{geom}}$ generated by the union (in the appropriate subobject lattice) of all the images of morphisms in $R$. In fact, these notions played an essential role in \cite{OC3}. The following result provides a link between these latter concepts and the notions of regular and geometric topology.\\
Regarding notation, given a small category $\cal C$ with a Grothendieck topology $J$ on it and a sieve $R$ in $\cal C$, we denote by $\overline{R}^{J}$ the $J$-closure of $R$, that is the sieve $\overline{R}^{J}:=\{f:d\to c \textrm{ | } f^{\ast}(R)\in J(d)\}$; recall that, via the identification of sieves on an object $c$ with subobjects in $[{\cal C}^{\textrm{op}}, \Set]$ of ${\cal C}(-,c)$, $\overline{R}^{J}$ corresponds to the closure of $R\mono {\cal C}(-,c)$ with respect to the universal closure operator on $[{\cal C}^{\textrm{op}}, \Set]$ corresponding to the (local operator associated to the) Grothendieck topology $J$ on $\cal C$. The monic part of the cover-mono factorization of an arrow $f$ in a regular category will be denoted by $Im(f)$.  

\begin{proposition}\label{regulargeometric}
Let $R$ be a sieve on a category $\cal C$. Then\\
(i) If $\cal C$ is a regular category then $R^{\textrm{reg}}=\overline{R}^{J^{\textrm{reg}}_{\cal C}}$;\\
(ii) If $\cal C$ is a geometric category then $R^{\textrm{geom}}=\overline{R}^{J^{\textrm{geom}}_{\cal C}}$.
\end{proposition}     
  
\begin{proofs}
(i) Let us begin by proving that the inclusion $R^{\textrm{reg}}\subseteq \overline{R}^{J^{\textrm{reg}}_{\cal C}}$ holds. Clearly, it suffices to show that for any $f$ in $R$, $Im(f)\in \overline{R}^{J^{\textrm{reg}}_{\cal C}}$; now, if $a\stackrel{f'}{\epi} a' \stackrel{Im(f)}{\mono} b$ is the cover-mono factorization of $f$ then $(f')\subseteq Im(f)^{\ast}(R)$; but $f'$ is a cover so $(f')\in J^{\textrm{reg}}_{\cal C}(a)$ and hence $Im(f)^{\ast}(R)$ is $J^{\textrm{reg}}_{\cal C}$-covering by property (ii) in Definition \ref{def2}. It remains to prove the other inclusion. If $f\in \overline{R}^{J^{\textrm{reg}}_{\cal C}}$ then $f^{\ast}(R)$ contains a cover, call it $h$. Since composition of covers is a cover then $f$ factors through $Im(f\circ h)$ and hence $f\in R^{\textrm{reg}}$, as required.\\
(ii) Let $R$ be a sieve $\{r_{i} \textrm{ | } i\in I \}$ on an object $c\in {\cal C}$ (for our purposes we can suppose $I$ to be a set without loss of generality, every geometric category being well-powered). Let us denote by $r:d\mono c$ the union in $\Sub_{\cal C}(c)$ of the $Im(r_{i})$ as $i$ varies in $I$ and by $h_{i}$ the (unique) factorization of $r_{i}$ through $r$ (for each $i\in I$). To prove the inclusion $R^{\textrm{geom}}\subseteq \overline{R}^{J^{\textrm{geom}}_{\cal C}}$, it is enough to show that $r\in \overline{R}^{J^{\textrm{geom}}_{\cal C}}$. Now, $r = \mathbin{\mathop{\textrm{\huge $\cup$}}\limits_{i\in I}} Im(r_{i})$ so $1_{d}=r^{\ast}(\mathbin{\mathop{\textrm{\huge $\cup$}}\limits_{i\in I}} Im(r_{i}))=\mathbin{\mathop{\textrm{\huge $\cup$}}\limits_{i\in I}}r^{\ast}(Im(r\circ h_{i}))=\mathbin{\mathop{\textrm{\huge $\cup$}}\limits_{i\in I}}Im(r^{\ast}(r\circ h_{i}))=\mathbin{\mathop{\textrm{\huge $\cup$}}\limits_{i\in I}}Im(h_{i})$, where the second and third equalities follows from the fact that in any geometric category cover-mono factorizations and small unions of subobjects are stable under pullback and the last equality follows from the fact that $r$ is monic. So we obtain that $\{h_{i} \textrm{ | } i\in I\}$ is a small covering family contained in $r^{\ast}(R)$ and hence $r^{\ast}(R)$ is $J^{\textrm{geom}}_{\cal C}$-covering, as required. Conversely, let us suppose that, given $f:d\to c$, $f^{\ast}(R)$ contains a small covering family $\{h_{j} \textrm{ | } j\in J\}$. We want to prove that $f$ factors through $r$. Since $r$ is monic, this condition is clearly equivalent to requiring that $f^{\ast}(r)=1_{d}$. Now, $f^{\ast}(r)=f^{\ast}(\mathbin{\mathop{\textrm{\huge $\cup$}}\limits_{i\in I}} Im(r_{i}))=\mathbin{\mathop{\textrm{\huge $\cup$}}\limits_{i\in I}}Im(f^{\ast}(r_{i}))$. For each $j\in J$ there exists $i\in I$ such that $f\circ h_{j}=r_{i}$ and hence $h_{j}$ factors through $f^{\ast}(r_{i})$; this in turn clearly implies that $Im(h_{j})$ factors through $Im(f^{\ast}(r_{i}))$, so that $\mathbin{\mathop{\textrm{\huge $\cup$}}\limits_{i\in I}}Im(f^{\ast}(r_{i}))\supseteq \mathbin{\mathop{\textrm{\huge $\cup$}}\limits_{j\in J}} Im(h_{j})=1_{d}$. Therefore $f^{\ast}(r)=1_{d}$, as required.  \end{proofs}  
\begin{rmk}\label{pbstable}
\emph{As a consequence of our proposition we may deduce that if $\cal C$ is regular (resp. geometric) then for any sieve $R$ on $c$ and any arrow $f:d\to c$, $f^{\ast}(R^{\textrm{reg}})=f^{\ast}(R)^{\textrm{reg}}$ (resp. $f^{\ast}(R^{\textrm{geom}})=f^{\ast}(R)^{\textrm{geom}}$); indeed, universal closure operators always commute with pullbacks.}
\end{rmk}

\subsection{Categories with logical structure as syntactic categories}

We recall from \cite{El2} that if $\mathbb T$ is a cartesian (resp. regular, coherent, geometric) theory over a signature $\Sigma$, one may construct the cartesian (resp. regular, coherent, geometric) syntactic category ${\cal C}_{\mathbb T}^{\textrm{cart}}$ (resp. ${\cal C}_{\mathbb T}^{\textrm{reg}}$, ${\cal C}_{\mathbb T}^{\textrm{coh}}$, ${\cal C}_{\mathbb T}^{\textrm{geom}}$) of $\mathbb T$. By Lemma D1.4.10 \cite{El}, this category is cartesian (resp. regular, coherent, geometric) and satisfies the property that the category of cartesian (resp. regular, coherent, geometric) functors from it to any cartesian (resp. regular, coherent, geometric) category $\cal D$ is naturally equivalent to the category of models of the theory $\mathbb T$ in $\cal D$, the equivalence sending each model $M\in {\mathbb T}\textrm{-mod}({\cal E})$ to the functor $F_{M}:{\cal C}_{\mathbb T}\rightarrow {\cal E}$ assigning to a formula $\phi(\vec{x})$ its interpretation $[[\phi(\vec{x})]]_{M}$ in $M$. Let us now show that, conversely, any cartesian (resp. regular, coherent, geometric) category can be regarded as (that is, it is naturally equivalent to) the syntactic category of a cartesian (resp. regular, coherent, geometric) theory. The ingredients for this result are all in \cite{El}, the main one being the construction of the canonical signature $\Sigma_{\cal C}$ of a category $\cal C$ with at least finite limits described at p. 837. This signature has one sort $\name{A}$ for each object $A$ of $\cal C$, one function symbol $\name{f}:\name{A_{1}},\cdots, \name{A_{n}}\to \name{B}$ for each arrow $f:A_{1}\times \cdots \times A_{n} \to B$ in $\cal C$, and one relation symbol $\name{R}\mono \name{A_{1}},\cdots, \name{A_{n}}$ for each subobject $R\mono A_{1}\times \cdots \times A_{n}$. Now, let ${\mathbb T}^{\cal C}$ be the theory formed by the following cartesian sequents over $\Sigma_{\cal C}$:
\[
(\top \vdash_{x} (\name{f}(x)=x))
\]
for any identity arrow $f$ in $\cal C$;
\[
(\top \vdash_{x} (\name{f}(x)=\name{h}(\name{g}(x))))
\]
for any triple of arrows $f,g,h$ of $\cal C$ such that $f$ is equal to the composite $h\circ g$;

\[
(\top \vdash_{[]} (\exists x)\top) \textrm{ and } (\top \vdash_{x,x'} (x=x'))
\]
where $x$ and $x'$ are of sort $\name{1}$, $1$ being the terminal object of $\cal C$;

\[
\begin{array}{c}
(\top \vdash_{x} (\name{h}(\name{f}(x))=\name{k}(\name{g}(x)))),\\
((\name{f}(x)=\name{f'}(x'))\wedge (\name{g}(x)=\name{g'}(x'))\vdash_{x,x'} (x=x')), \textrm{and}\\
((\name{h}(y)=\name{k}(z)) \vdash_{y,z} (\exists x)((\name{f}(x)=y)\wedge (\name{g}(x)=z)))
\end{array}
\] 
for any pullback square
\[  
\xymatrix {
a \ar[r]^{f} \ar[d]^{g}  & b \ar[d]^{h} \\
c \ar[r]^{k} & d }
\] 
in $\cal C$.\\
It is an immediate consequence of Lemma D1.3.11 \cite{El2} that for any cartesian category $\cal D$, the ${\mathbb T}^{\cal C}$-models are the same thing as functors ${\cal C}\to {\cal D}$ i.e. cartesian functors (cfr. Example D1.4.8 \cite{El2}). So we have an equivalence of categories ${\mathbb T}^{\cal C}\textrm{-mod}({\cal D})\simeq \mathfrak{Cart}({\cal C}, {\cal D})$ natural in ${\cal D}\in \mathfrak{Cart}$. Since we also have an equivalence $\mathfrak{Cart}({\cal C}^{\textrm{cart}}_{{\mathbb T}^{\cal C}}, {\cal D}) \simeq {\mathbb T}^{\cal C}\textrm{-mod}({\cal D})$ natural in ${\cal D}\in \mathfrak{Cart}$ (by definition of syntactic category), by composing the two we find an equivalence $\mathfrak{Cart}({\cal C}, {\cal D}) \simeq \mathfrak{Cart}({\cal C}^{\textrm{cart}}_{{\mathbb T}^{\cal C}}, {\cal D})$ natural in ${\cal D}\in \mathfrak{Cart}$ and hence, by the $2$-dimensional Yoneda Lemma, a natural equivalence of categories ${\cal C}^{\textrm{cart}}_{{\mathbb T}^{\cal C}} \simeq {\cal C}$, one half of which sends a formula $\phi(\vec{x})$ to (the domain of) its interpretation $[[\phi(\vec{x})]]$ in the canonical $\Sigma_{\cal C}$-structure in $\cal C$.\\
One can easily extend this result to more general fragments of geometric logic. Indeed, given a Grothendieck topology $J$ on a category $\cal C$, recall from \cite{El2} (Remark D3.1.13) that the cartesian and $J$-cover-preserving (i.e. which send every $J$-covering sieve to a covering family) functors on $\cal C$ correspond exactly to the models of the theory ${\mathbb T}^{\cal C}$ which satisfy the additional axioms
\[
(\top \vdash_{x} \mathbin{\mathop{\textrm{\huge $\vee$}}\limits_{i\in I}}(\exists y_{i})(\name{f_{i}}(y_{i})=x))  
\]
for each $J$-covering family $(f_{i}:B_{i}\to A \textrm{ | } i\in I)$. Let us call ${\mathbb T}^{\cal C}_{J}$ the theory obtained from ${\mathbb T}^{\cal C}$ by adding the axioms above. Now, it is easy to verify that if $\cal C$ is a regular (resp. coherent, geometric) category then for any regular (resp. coherent, geometric) category $\cal D$, the regular (resp. coherent, geometric) functors ${\cal C}\to {\cal D}$ are exactly the cartesian functors on $\cal C$ which are $J$-cover-preserving, where $J$ is the regular (resp. coherent, geometric) coverage on $\cal C$. So we conclude as above that if $\cal C$ is a regular (resp. coherent, geometric) category then there is an equivalence of categories ${\cal C}_{{\mathbb T}^{\cal C}_{J}}^{\textrm{reg}} \simeq {\cal C}$ (resp. ${\cal C}_{{\mathbb T}^{\cal C}_{J}}^{\textrm{coh}}\simeq {\cal C}$, ${\cal C}_{{\mathbb T}^{\cal C}_{J}}^{\textrm{geom}}\simeq {\cal C}$) one half of which sends a formula $\phi(\vec{x})$ to (the domain of) its interpretation $[[\phi(\vec{x})]]$ in the canonical $\Sigma_{\cal C}$-structure in $\cal C$.\\
Hence we have arrived at the following result
\begin{proposition}\label{loccat}
The cartesian (resp. regular, coherent, geometric) categories are, up to natural equivalence, exactly the syntactic categories of cartesian (resp. regular, coherent, geometric) theories.
\end{proposition}\qed
We note that the fact that every cartesian (resp. regular, coherent, geometric) category $\cal C$ is naturally equivalent to the syntactic category of a theory $\mathbb T$ enables us to interpret categorical constructions on $\cal C$ as logical operations involving $\mathbb T$.  

\newpage
\section{The duality theorem}\label{section_dualita}
In this section we prove our main theorem, which asserts the existence of a bijection between the subtoposes of the classifying topos of a geometric theory $\mathbb{T}$ over $\Sigma$ and the closed geometric theories over $\Sigma$ which are `quotients' of the theory $\mathbb{T}$.\\
Let us start with an easy remark: every subtopos of a Grothendieck topos is a Grothendieck topos. This can be proved in (at least) two different ways, as follows.\\
We recall that a subtopos of a topos $\cal E$ is a geometric inclusion of the form $\sh_{j}({\cal E})\hookrightarrow {\cal E}$ for a local operator $j$ on $\cal E$, equivalently an equivalence class of geometric inclusions to the topos $\cal E$. It is well-known that the subtoposes of a presheaf topos $[{\cal C}^{\textrm{op}}, \Set]$ are in bijection with the Grothendieck topologies $J$ on the category $\cal C$, i.e. every geometric inclusion to $[{\cal C}^{\textrm{op}}, \Set]$ is, up to equivalence, of the form $\Sh({\cal C}, J)\hookrightarrow [{\cal C}^{\textrm{op}}, \Set]$ for a unique Grothendieck topology $J$ on $\cal C$; moreover, a geometric inclusion $\Sh({\cal C}, J)\hookrightarrow [{\cal C}^{\textrm{op}}, \Set]$ factors through another geometric inclusion $\Sh({\cal C}, J')\hookrightarrow [{\cal C}^{\textrm{op}}, \Set]$ of the same form if and only if $J'\subseteq J$ (i.e. every $J'$-covering sieve is a $J$-covering sieve). Now, the geometric inclusions to a Grothendieck topos $\Sh({\cal C}, J)$ can be clearly identified with the geometric inclusions to $[{\cal C}^{\textrm{op}}, \Set]$ which factors through $\Sh({\cal C}, J)\hookrightarrow [{\cal C}^{\textrm{op}}, \Set]$ and hence the subtoposes of $\Sh({\cal C}, J)$ correspond precisely to the Grothendieck topologies $J'$ on $\cal C$ such that $J'\supseteq J$. This provides us with the first proof of our claim. Alternatively, we can argue as follows. By Theorem C2.2.8 \cite{El2}, an elementary topos $\cal E$ is a Grothendieck topos if and only if there exists a bounded geometric morphism ${\cal E}\to \Set$ (cfr. B3.1.7 \cite{El}). Now, a geometric inclusion is always a localic morphism (cfr. Example A4.6.2(a) \cite{El}), and hence a bounded morphism (cfr. Example B3.1.8 \cite{El}); but a composite of bounded morphism is a bounded morphism (by Lemma B3.1.10(i)), so that our thesis follows from the above-mentioned characterization.\\            
Our remark is fundamental for our purposes for the following reason. For each elementary topos $\cal E$, the collection of subtoposes of $\cal E$ has the structure of a coHeyting algebra (cfr. Example A4.5.13(f) \cite{El}), and there are many important concepts in topos theory that apply to this context (cfr. section A4 \cite{El}); so we are naturally led to investigating their meaning in the context of Grothendieck toposes.  In fact, thanks to the duality theorem established below, we will also be able to interpret all these concepts in the context of geometric theories. All of this will be carried out in the following sections of the paper.\\ 
Before we can state our duality theorem, which describes how the relationship between Grothendieck toposes and geometric theories given by the theory of classifying toposes `restricts' to the context of all the subtoposes of a given Grothendieck topos, we need to introduce some definitions. Regarding terminology, we use the term theory to mean a \emph{presentation} of a theory, that is a collection of axioms of the theory, and accordingly we consider two theories over a given signature equal when they have exactly the same axioms.

\begin{definition}
Let $\mathbb T$ be a geometric theory over a signature $\Sigma$ and $\sigma, \sigma'$ two geometric sequents over $\Sigma$. Then $\sigma$ and $\sigma'$ are said to be $\mathbb T$-equivalent if $\sigma$ is provable in ${\mathbb T}\cup \{\sigma'\}$ and $\sigma'$ is provable in ${\mathbb T}\cup \{\sigma\}$.
\end{definition} 

\begin{definition}
Let $\mathbb{T}$ be a geometric theory over a signature $\Sigma$. A quotient of $\mathbb{T}$ is a geometric theory $\mathbb{T}'$ over $\Sigma$ such that every axiom of $\mathbb{T}$ is provable in $\mathbb{T}'$.
\end{definition}

\begin{rmk}     
\emph{The notion of provability in geometric logic to which we refer here (and below) is that defined p. 832 \cite{El2}; that system is essentially constructive, but, by Proposition D3.1.16 \cite{El2}, we may add the law of excluded middle to it (thus making it classical) without affecting the correponding notion of provability.}
\end{rmk} 
\begin{definition}
Let $\mathbb{T}$ and $\mathbb{T}$ be geometric theories over a signature $\Sigma$. We say that $\mathbb{T}$ and $\mathbb{T}$ are syntactically equivalent, and we write $\mathbb{T} \equiv_{s} \mathbb{T}$, if for every geometric sequent $\sigma$ over $\Sigma$, $\sigma$ is provable in $\mathbb{T}$ if and only if $\sigma$ is provable in $\mathbb{T}'$.
\end{definition}     
We note that we can take a canonical representative for each of the $\equiv_{s}$-equivalence classes, namely the theory having as axioms exactly the geometric sequents over $\Sigma$ which are provable in one (equivalently, all) of the theories belonging to that equivalence class.\\
Borrowing a term from classical model theory, we will say that a geometric theory $\mathbb{T}$ over a signature $\Sigma$ is \emph{closed} if all the geometric sequents over the signature of $\mathbb{T}$ which are provable in $\mathbb{T}$ already belong to $\mathbb{T}$. Thus, there is exactly one closed theory in every $\equiv_{s}$-equivalence class, which is in fact our canonical representative. Accordingly, we define the \emph{closure} of a geometric theory over a given signature as the unique closed theory in its $\equiv_{s}$-equivalence class.\\
Let us recall the following definition.
\begin{definition}
Let $\mathbb{T}$ and $\mathbb{T}'$ be geometric theories. We say that $\mathbb{T}$ and $\mathbb{T}'$ are Morita-equivalent if they have equivalent classifying toposes (equivalently, by the $2$-dimensional Yoneda Lemma, if they have equivalent categories of models in every Grothendieck topos $\cal E$, naturally in ${\cal E}\in \mathfrak{BTop}$).
\end{definition}  
We are now ready to state our duality theorem.     
Concerning notation, given two Grothendieck toposes $\cal E$ and $\cal F$ and a Grothendieck topology $J$ on a small category $\cal C$, we denote by ${\bf Geom}({\cal E}, {\cal F})$ the category of geometric morphisms from $\cal E$ to $\cal F$ and by ${\bf Flat}_{J}({\cal C}, {\cal E})$ the category of $J$-continuous flat functors from $\cal C$ to $\cal E$.
\begin{theorem}\label{dualita}
Let $\mathbb{T}$ be a geometric theory over a signature $\Sigma$. Then the assignment sending a quotient of $\mathbb{T}$ to its classifying topos defines a bijection between the $\equiv_{s}$-equivalence classes of quotients of $\mathbb T$ (equivalently, the closed quotients of $\mathbb T$) and the subtoposes of the classifying topos $\Set[\mathbb{T}]$ of $\mathbb{T}$.   
\end{theorem}

\begin{proofs}
First, we note that two syntactically equivalent theories are Morita-equivalent; indeed, by the soundness theorem for geometric logic, they have the same (categories of) models in every Grothendieck topos.
Let us recall from \cite{El2} that the classifying topos $\Set[\mathbb T]$ of $\mathbb T$ can be represented as $\Sh({\cal C}_{\mathbb T}, J_{\mathbb T})$, where ${\cal C}_{\mathbb T}$ is the geometric syntactic category of $T$ and $J_{\mathbb T}$ is the canonical topology on ${\cal C}_{\mathbb{T}}$ (i.e. the Grothendieck topology on ${\cal C}_{\mathbb T}$ having as covering sieves exactly those which contain small covering families), and that we have an equivalence of categories ${\mathbb T}\textrm{-mod}({\cal E})\simeq {\bf Flat}_{J_{\mathbb T}}({\cal C}_{\mathbb T}, {\cal E})$ (natural in ${\cal E}\in \mathfrak{BTop}$) which sends each model $M\in {\mathbb T}\textrm{-mod}({\cal E})$ to the functor $F_{M}:{\cal C}_{\mathbb T}\rightarrow {\cal E}$ assigning to a formula $\{\vec{x}. \phi\}$ (the domain of) its interpretation $[[\phi(\vec{x})]]_{M}$ in $M$.\\
We note that, although not small, ${\cal C}_{\mathbb T}$ is an essentially small category i.e. it is equivalent to a small category (by the results in Part D \cite{El2}); hence all the results valid for small Grothendieck sites naturally extend to sites involving the category ${\cal C}_{\mathbb T}$.\\
Let us recall the construction of pullbacks in ${\cal C}_{\mathbb{T}}$. Given two morphisms
\[  
\xymatrix {
\{\vec{x}.\phi\} \ar[r]^{[\theta]} & \{\vec{y}.\psi\}}
\]
and
\[  
\xymatrix {
\{\vec{x'}.\phi'\} \ar[r]^{[\theta']} & \{\vec{y}.\psi\}}
\]
in ${\cal C}_{\mathbb{T}}$ with common codomain, we have the following pullback in ${\cal C}_{\mathbb{T}}$:
\[  
\xymatrix {
\{\vec{\underline{x}},\vec{\underline{x'}}. (\exists \vec{y})(\theta[\vec{\underline{x}}\slash \vec{x}]\wedge \theta'[\vec{\underline{x'}}\slash \vec{x'}])\} \ar[rrrrr]^{[(\exists \vec{y})(\theta \wedge \theta' \wedge \vec{\underline{x'}}=\vec{x'})]} \ar[d]^{[(\exists \vec{y})(\theta \wedge \theta' \wedge \vec{\underline{x}}=\vec{x})]} & & & & & \{\vec{x'}.\phi'\} \ar[d]^{[\theta']} \\
\{\vec{x}. \phi\} \ar[rrrrr]^{[\theta]} & & & & & \{\vec{y}. \psi\} }
\]\\ 
Let us note that the sequent $\phi' \vdash_{\vec{x'}} (\exists \vec{\underline{x}},\vec{\underline{x'}})((\exists \vec{y})(\theta \wedge \theta' \wedge \vec{\underline{x'}}=\vec{x'}))$ is provable in $\mathbb T$ from the sequent $\psi \vdash_{\vec{y}} (\exists \vec{x})\theta$. Indeed, it is clearly equivalent in geometric logic to the sequent $\phi' \vdash_{\vec{x'}} (\exists \vec{y})((\exists \vec{\underline{x}})\theta \wedge \theta')$, and the sequents $\phi' \vdash_{\vec{x'}} (\exists \vec{y})\theta'$ and $\theta' \vdash_{\vec{x'}, \vec{y}} \psi$ are provable in $\mathbb T$ since $[\theta']$ is a morphisms in the syntactic category ${\cal C}_{\mathbb{T}}$.\\

Next, we observe that, given a $\mathbb T$-model $M$ in a Grothendieck topos $\cal E$, $F_{M}:{\cal C}_{\mathbb T}\rightarrow {\cal E}$ sends a small family $\{\theta_{i} \textrm{ | } i\in I\}$ of morphisms
\[  
\xymatrix {
\{\vec{x_{i}}.\phi_{i}\} \ar[r]^{[\theta_{i}]} & \{\vec{y}.\psi\}}
\]
in ${\cal C}_{\mathbb T}$ with common codomain to an epimorphic family in $\cal E$ if and only if $[[\vec{y}.\psi]]_{M}=[[\vec{y}. \mathbin{\mathop{\textrm{\huge $\vee$}}\limits_{i\in I}}(\exists \vec{x_{i}})\theta_{i}]]_{M}$, equivalently if and only if the sequent $\psi \vdash_{\vec{y}} \mathbin{\mathop{\textrm{\huge $\vee$}}\limits_{i\in I}}(\exists \vec{x_{i}})\theta_{i}$ is satisfied in $M$.\\ 

This remark shows, by the soundess theorem for geometric logic, that for any small presieve $R$ in ${\cal C}_{\mathbb{T}}$, the $J_{\mathbb T}$-continuous flat functors on ${\cal C}_{\mathbb{T}}$ sending $R$ to an epimorphic family also send all the pullbacks of $R$ along arrows in ${\cal C}_{\mathbb{T}}$ to epimorphic families. This implies, by Remark \ref{presieves} and Lemma $3$ \cite{blasce}, that the $J_{\mathbb T}$-continuous flat functors on ${\cal C}_{\mathbb{T}}$ which send each of the small presieves in a given collection $\cal F$ to an epimorphic family coincide with the $J_{\mathbb T}$-continuous flat functors on ${\cal C}_{\mathbb{T}}$ which are $J_{\cal F}$-continuous, where $J_{\cal F}$ is the Grothendieck topology on ${\cal C}_{\mathbb{T}}$ generated over $J_{\mathbb T}$ by the sieves generated by presieves in $\cal F$.\\ 
  
Given a quotient $\mathbb{T}'$ of $\mathbb{T}$, we may construct its classifying topos as follows. Let $\mathbb{T}'$ be obtained from $\mathbb T$ by adding a number of axioms of the form $\phi \vdash_{\vec{x}} \psi$, where $\phi$ and $\psi$ are geometric formulae over $\Sigma$; of course, up to syntactic equivalence, there are many possible ways of presenting $\mathbb{T}'$ in such form (for example one may take as axioms all the axioms of $\mathbb{T}'$ or, more economically, all the axioms of ${\mathbb T}'$ which are not provable in $\mathbb T$), but we will show that our construction is independent from any particular presentation. For each of these axioms $\phi \vdash_{\vec{x}} \psi$, consider the corresponding morphism
\[  
\xymatrix {
\{\vec{x'}.\phi\wedge \psi\} \ar[rrr]^{[(\phi\wedge \psi\wedge \vec{x'}=\vec{x})]} & & & \{\vec{x}.\phi\}}
\]
in the geometric syntactic category ${\cal C}_{\mathbb{T}}$ of $\mathbb T$.\\ 
It is clear that, given a $\mathbb T$-model $M$ in a Grothendieck topos $\cal E$, $F_{M}:{\cal C}_{\mathbb T}\rightarrow {\cal E}$ sends the morphism
\[  
\xymatrix {
\{\vec{x'}.\phi\wedge \psi\} \ar[rrr]^{[(\phi\wedge \psi\wedge \vec{x'}=\vec{x})]} & & & \{\vec{x}.\phi\}}
\]
to an epimorphism if and only if $[[\vec{x}.\phi]]_{M}\leq [[\vec{x}.\psi]]_{M}$ i.e. if and only if the sequent $\phi \vdash_{\vec{x}} \psi$ holds in $M$.\\  

So the $J_{\mathbb T}$-continuous flat functors on ${\cal C}_{\mathbb{T}}$ which send each of the morphism corresponding to the axioms of ${\mathbb T}'$ to an epimorphism classify the models of ${\mathbb T}'$. Therefore, from the discussion above we deduce that if $J_{\mathbb T'}^{\mathbb T}$ is the smallest Grothendieck topology on ${\cal C}_{\mathbb{T}}$ for which all the $J_{\mathbb T}$-covering sieves and the sieves containing a morphism corresponding to an axiom of ${\mathbb T}'$ are $J_{\mathbb T'}^{\mathbb T}$-covering then, by Diaconescu's theorem, the topos $\Sh({\cal C}_{\mathbb{T}},J_{\mathbb T'}^{\mathbb T})$ classifies the theory $\mathbb T'$; moreover, the canonical geometric inclusion $\Sh({\cal C}_{\mathbb{T}},J_{\mathbb T'}^{\mathbb T}) \hookrightarrow \Sh({\cal C}_{\mathbb{T}},J_{\mathbb T})$ corresponding to the inclusion $J_{\mathbb T} \subseteq  J_{\mathbb T'}^{\mathbb T}$ makes $\Sh({\cal C}_{\mathbb{T}},J_{\mathbb T'}^{\mathbb T})$ into a subtopos of $\Set[\mathbb T]$.\\ 
Now, to have a well-defined assignment from the $\equiv_{s}$-equivalence classes of quotients of $\mathbb T$ to the subtoposes of $\Set[\mathbb T]$, it remains to verify that the topology $J_{\mathbb T'}^{\mathbb T}$ defined above does not depend on the particular choice of axioms for ${\mathbb T}'$, i.e. it is the same for all the quotients in a given $\equiv_{s}$-equivalence class.\\  
Let ${\mathbb T}_{1}$ and ${\mathbb T}_{2}$ be quotients of $\mathbb T$ such that ${\mathbb T}_{1}\equiv_{s}{\mathbb T}_{2}$; we want to prove that $J_{{\mathbb T}_{1}}^{\mathbb T}=J_{{\mathbb T}_{2}}^{\mathbb T}$. We will prove the existence of a geometric equivalence $\tau:\Sh({\cal C}_{\mathbb{T}},J_{{\mathbb T}_{1}}^{\mathbb T})\rightarrow \Sh({\cal C}_{\mathbb{T}},J_{{\mathbb T}_{2}}^{\mathbb T})$ such that the diagram in $\mathfrak{BTop}$  
\[  
\xymatrix {
\Sh({\cal C}_{\mathbb{T}},J_{{\mathbb T}_{1}}^{\mathbb T}) \ar_{i_{1}}[dr] \ar^{\tau}[rr]  &  & \Sh({\cal C}_{\mathbb{T}},J_{{\mathbb T}_{2}}^{\mathbb T}) \ar^{i_{2}}[dl]  \\
& [{{\cal C}_{\mathbb T}}^{\textrm{op}}, \Set ] &}
\]
where the geometric inclusions $\Sh({\cal C}_{\mathbb{T}},J_{{\mathbb T}_{1}}^{\mathbb T})\to [{{\cal C}_{\mathbb T}}^{\textrm{op}}, \Set ]$ and $\Sh({\cal C}_{\mathbb{T}},J_{{\mathbb T}_{2}}^{\mathbb T})\to [{{\cal C}_{\mathbb T}}^{\textrm{op}}, \Set ]$ are the canonical ones, commutes up to isomorphism.\\

From the identification of equivalence classes of geometric inclusions to a given topos with local operators on that topos (given by the theory of elementary toposes) it will then follow the equality of the two topologies $J_{{\mathbb T}_{1}}^{\mathbb T}$ and $J_{{\mathbb T}_{2}}^{\mathbb T}$. By the $2$-dimensional Yoneda Lemma, it is equivalent to prove the existence of an equivalence of categories $l_{\cal E}:{\bf Geom}({\cal E},\Sh({\cal C}_{\mathbb{T}},J_{{\mathbb T}_{1}}^{\mathbb T})) \rightarrow {\bf Geom}({\cal E},\Sh({\cal C}_{\mathbb{T}},J_{{\mathbb T}_{2}}^{\mathbb T}))$ natural in ${\cal E}\in \mathfrak{BTop}$ such that $(i_{1}\circ -) \circ l_{\cal E}\cong (i_{2}\circ -)$ for each ${\cal E}\in \mathfrak{BTop}$. Since ${\mathbb T}_{1}\equiv_{s}{\mathbb T}_{2}$, ${\mathbb T}_{1}$ and ${\mathbb T}_{2}$ have the same models (in every Grothendieck topos), and hence we may obtain such an equivalence by composing ${\bf Geom}({\cal E},\Sh({\cal C}_{\mathbb{T}},J_{{\mathbb T}_{1}}^{\mathbb T})) \simeq {\bf Flat}_{J_{{\mathbb T}_{1}}}({\cal C}_{\mathbb T}, {\cal E}) \simeq {\mathbb T}_{1}\textrm{-mod}({\cal E}) = {\mathbb T}_{2}\textrm{-mod}({\cal E}) \simeq {\bf Flat}_{J_{{\mathbb T}_{2}}}({\cal C}_{\mathbb T}, {\cal E}) \simeq {\bf Geom}({\cal E},\Sh({\cal C}_{\mathbb{T}},J_{{\mathbb T}_{2}}^{\mathbb T}))$, where the first and last equivalences are given by Diaconescu's theorem.\\
Conversely, suppose starting with a subtopos $\cal E$ of $\Set[\mathbb T]$; then $\cal E$ has the form $\Sh({\cal C}_{\mathbb T}, J)$ for a unique Grothendieck topology $J$ such that $J\supseteq J_{\mathbb T}$. Let us prove that there exists a quotient ${\mathbb T}^{J}$ of $\mathbb T$ such that $\cal E$ is its classifying topos. Let us define ${\mathbb T}^{J}$ to consist of all the axioms over $\Sigma$ of the form $\psi \vdash_{\vec{y}} (\exists \vec{x})\theta$, where $[\theta]$ is any monomorphism
\[  
\xymatrix {
\{\vec{x}.\phi\} \ar[r]^{[\theta]} & \{\vec{y}.\psi\}}
\]
in ${\cal C}_{\mathbb T}$ generating a $J$-covering sieve.\\
Since, for any $\mathbb T$-model $M$ in a Grothendieck topos $\cal E$, $F_{M}$ sends $[\theta]$ to an epimorphism if and only if the sequent $\psi \vdash_{\vec{y}} (\exists \vec{x})\theta$ holds in $M$, it follows from Remark \ref{presieves}, Proposition \ref{gengrot} and Lemma $3$ \cite{blasce} that the equivalence ${\mathbb T}\textrm{-mod}({\cal E}) \simeq {\bf Flat}_{J_{{\mathbb T}}}({\cal C}_{\mathbb T}, {\cal E})$ restricts to an equivalence ${\mathbb T}^{J}\textrm{-mod}({\cal E}) \simeq {\bf Flat}_{J}({\cal C}_{\mathbb T}, {\cal E})$ (naturally in ${\cal E} \in \mathfrak{BTop}$) and hence that ${\cal E}=\Sh({\cal C}_{\mathbb T}, J)$ classifies the theory ${\mathbb T}^{J}$.\\
To conclude the proof of the theorem it remains to show that the two assignments ${\mathbb T}'\to J_{{\mathbb T}'}^{\mathbb T}$ and $J\to {\mathbb T}^{J}$ are bijections inverse to each other between the $\equiv_{s}$-equivalence classes of quotients of $\mathbb T$ and the subtoposes of the classifying topos $\Set[\mathbb{T}]$ of $\mathbb{T}$.\\

To prove that for any quotient ${\mathbb T}'$ of $\mathbb T$ we have ${\mathbb T}' \equiv_{s} {\mathbb T}^{J_{{\mathbb T}'}^{\mathbb T}}$ we argue as follows. First, we observe that for any $\mathbb T$-model $M$ in a Grothendieck topos $\cal E$, $M$ is a ${\mathbb T}'$-model if and only if it is a ${\mathbb T}^{J_{{\mathbb T}'}^{\mathbb T}}$-model; indeed, by definition of $J_{{\mathbb T}'}^{\mathbb T}$ and of ${\mathbb T}^{J_{{\mathbb T}'}^{\mathbb T}}$, both ${\mathbb T}'$-models and ${\mathbb T}^{J_{{\mathbb T}'}^{\mathbb T}}$-models in $\cal E$ correspond to functors in ${\bf Flat}_{J_{{\mathbb T}'}^{\mathbb T}}({\cal C}_{{\mathbb T}'}, {\cal E})$ via the equivalence ${\bf Flat}_{J_{\mathbb T}}({\cal C}_{\mathbb T}, {\cal E}) \simeq {\mathbb T}\textrm{-mod}({\cal E})$.\\

Now, let us denote by $U_{{\mathbb T}'}^{\mathbb T}$ the image of $a_{J_{{\mathbb T}'}^{\mathbb T}} \circ y^{{\mathbb T}}$ in ${\mathbb T}'\textrm{-mod}({\cal E})$ through the equivalence ${\bf Flat}_{J_{{\mathbb T}'}^{\mathbb T}}({\cal C}_{\mathbb T}, {\cal E}) \simeq {\mathbb T}'\textrm{-mod}({\cal E})$, where $y^{{\mathbb T}}:{\cal C}_{\mathbb T}\to [{{\cal C}_{\mathbb T}}^{\textrm{op}}, \Set]$ is the Yoneda embedding and $a_{J_{{\mathbb T}'}^{\mathbb T}}:[{{\cal C}_{\mathbb T}}^{\textrm{op}}, \Set] \to \Sh({\cal C}_{\mathbb{T}},J_{{\mathbb T}'}^{\mathbb T})$ is the associated sheaf functor. By Diaconescu's theorem and the naturality in ${\cal E}\in \mathfrak{BTop}$ of the equivalences ${\mathbb T}'\textrm{-mod}({\cal E}) = {\mathbb T}^{J_{{\mathbb T}'}^{\mathbb T}}\textrm{-mod}({\cal E}) \simeq {\bf Flat}_{J_{{\mathbb T}'}^{\mathbb T}}({\cal C}_{\mathbb T}, {\cal E}) \simeq {\bf Geom}({\cal E},\Sh({\cal C}_{\mathbb{T}},J_{{\mathbb T}'}^{\mathbb T}))$, the $\Sigma$-structure $U_{{\mathbb T}'}^{\mathbb T}$ is a universal model for both ${\mathbb T}'$ and ${\mathbb T}^{J_{{\mathbb T}'}^{\mathbb T}}$ (i.e. every ${\mathbb T}'$-model $M$ in a Grothendieck topos $\cal G$ is the image $g^{\ast}(U_{{\mathbb T}'}^{\mathbb T})$ for a unique up to isomorphism geometric morphism $g:{\cal G} \to \Sh({\cal C}_{\mathbb{T}},J_{{\mathbb T}'}^{\mathbb T})$); in particular, it is conservative both as a ${\mathbb T}'$-model and as a ${\mathbb T}^{J_{{\mathbb T}'}^{\mathbb T}}$-model (since for every geometric theory $\mathbb Z$ over a signature $\Sigma'$ its classifying topos $\Sh({\cal C}_{\mathbb Z}, J_{\mathbb Z})$ contains a conservative $\mathbb Z$-model, cfr. the discussion preceding Proposition D3.1.12 \cite{El2}). From this it clearly follows that ${\mathbb T}' \equiv_{s} {\mathbb T}^{J_{{\mathbb T}'}^{\mathbb T}}$, as required.\\

On the other hand, the fact that $J=J_{{\mathbb T}^{J}}^{\mathbb T}$ directly follows from the definition of the assigmnent ${\mathbb T}'\to J_{{\mathbb T}'}^{\mathbb T}$.

\end{proofs}

With the above notation, we will refer to the topology $J_{{\mathbb T}'}^{\mathbb T}$ as the \emph{associated $\mathbb T$-topology of ${\mathbb T}'$}, and to the ($\equiv_{s}$-equivalence class of the) quotient ${\mathbb T}^{J}$ as the \emph{associated $\mathbb T$-quotient of $J$}.\\ 

For each Grothendieck topos $\cal E$, we denote by $\tau^{\cal E}:{\mathbb T}\textrm{-mod}({\cal E}) \simeq {\bf Geom}({\cal E}, \Sh({\cal C}_{\mathbb{T}},J_{\mathbb T}))$ the composite of the equivalence ${\mathbb T}\textrm{-mod}({\cal E})\simeq {\bf Flat}_{J_{\mathbb T}}({\cal C}_{\mathbb T}, {\cal E})$ considered in the proof of the theorem with Diaconescu's equivalence ${\bf Flat}_{J_{\mathbb T}}({\cal C}_{\mathbb T}, {\cal E}) \simeq {\bf Geom}({\cal E},\Sh({\cal C}_{\mathbb{T}},J_{{\mathbb T}}))$; given a quotient ${\mathbb T}'$ of a theory $\mathbb T$, we denote by $i^{\cal E}_{{\mathbb T}'}:{{\mathbb T}'}\textrm{-mod}({\cal E})\hookrightarrow {\mathbb T}\textrm{-mod}({\cal E})$ the inclusion into ${\mathbb T}\textrm{-mod}({\cal E})$ of the full subcategory ${{\mathbb T}'}\textrm{-mod}({\cal E})$ on the ${\mathbb T}'$-models in $\cal E$.   

\begin{rmk}\label{diagrammafond}
\emph{With the notation above, we note that, given a Grothendieck topology $J$ on ${\cal C}_{\mathbb T}$ such that $J\supseteq J_{\mathbb T}$ with corresponding canonical geometric inclusion $i_{J}:\Sh({\cal C}_{\mathbb{T}},J)\hookrightarrow \Sh({\cal C}_{\mathbb{T}},J_{\mathbb T})$, the duality theorem asserts in particular that there exists exactly one quotient ${\mathbb T}'$ of ${\mathbb T}$, up to syntactic equivalence, such that the diagram in $\mathfrak{Cat}$} 

\[  
\xymatrix {
{{\mathbb T}'}\emph{\textrm{-mod}}({\cal E}) \ar[rr]^{\simeq} \ar[d]^{i^{\cal E}_{{\mathbb T}'}} & & {\bf Geom}({\cal E},\Sh({\cal C}_{\mathbb{T}},J)) \ar[d]^{i_{J}\circ -} \\
{\mathbb T}\emph{\textrm{-mod}}({\cal E}) \ar[rr]^{\simeq}_{\tau^{\cal E}} & &  {\bf Geom}({\cal E},\Sh({\cal C}_{\mathbb{T}},J_{\mathbb T}))}
\]

\emph{commutes (up to invertible natural equivalence) naturally in ${\cal E}\in \mathfrak{BTop}$.}
\end{rmk}

We remark that our method of constructing the $\mathbb T$-topology associated to a given quotient of $\mathbb T$ has points in common with the `forcing' method summarized by Proposition D3.1.10 \cite{El2}. In fact, our arguments show that, more generally, it is always possible to construct the classifying topos of a quotient ${\mathbb T}'$ of a given theory $\mathbb T$ as a category of sheaves on the cartesian (resp. regular, coherent or geometric) syntactic category of $\mathbb T$ starting from a way of expressing ${\mathbb T}'$ as a theory obtained from $\mathbb T$ by adding axioms of the form $\psi \vdash_{\vec{y}} \mathbin{\mathop{\textrm{\huge $\vee$}}\limits_{i\in I}}(\exists \vec{x_{i}})\theta_{i}$ where $\psi$ and the $\theta_{i}$ are cartesian (resp. regular, coherent or geometric) formulae.\\ 

Finally, consider the following question: given a Grothendieck topos $\cal E$ and a signature $\Sigma$, when is it the case that there exists a geometric theory $\mathbb T$ over $\Sigma$ such that $\cal E$ is a classifying topos for $\mathbb T$? Our duality theorem gives us an answer to this question: the condition on $\cal E$ and $\Sigma$ amounts precisely to requiring that there should exist a geometric inclusion from $\cal E$ to the classifying topos for the empty (geometric) theory over $\Sigma$. 

\newpage
\section{The proof-theoretic interpretation} \label{section_prooftheory}

In this section, we provide an alternative, syntactic, proof of our duality theorem. This will be based on a logical interpretation of the notion of Grothendieck topology. Specifically, given a collection $\cal A$ of sieves on a given category $\cal C$, the notion of Grothendieck topology on $\cal C$ gives naturally rise to a proof system ${\cal T}_{\cal C}^{\cal A}$, as follows: the axioms of ${\cal T}_{\cal C}^{\cal A}$ are the sieves in $\cal A$ together with all the maximal sieves, while the inference rules of ${\cal T}_{\cal C}^{\cal A}$ are the proof-theoretic versions of the well-known axioms for Grothendieck topologies, i.e. the rules:\\
\emph{Stability rule:}
\[
\begin{array}{c}
R\\
\hline
f^{\ast}(R)  
\end{array}   
\] 
where $R$ is any sieve on an object $c$ in ${\cal C}$ and $f$ is any arrow in ${\cal C}$ with codomain $c$.\\
\emph{Transitivity rule:}
\[
\begin{array}{c}
Z \textrm{     } \{f^{\ast}(R) \textrm{ | } f\in Z \}\\
\hline
R  
\end{array}   
\] 
where $R$ and $Z$ are sieves in ${\cal C}$ on a given object of $\cal C$.\\ 
Notice that the `closed theories' of this proof system are precisely the Grothendieck topologies on $\cal C$ which contain the sieves in $\cal A$ as their covering sieves, and the closure of a theory in ${\cal T}_{\cal C}^{\cal A}$ i.e. of a collection $\cal U$ of sieves in $\cal C$, is exactly the Grothendieck topology on $\cal C$ generated by $\cal A$ and $\cal U$.\\ 
Our Theorem \ref{prooftheory} below can be interpreted as giving a `proof-theoretic equivalence' between the system of geometric logic over a given geometric theory $\mathbb T$ and the system ${\cal T}_{{\cal C}_{\mathbb T}}^{J_{\mathbb T}}$.\\ 
 
Given a geometric theory $\mathbb T$ over a signature $\Sigma$, let $\cal S$ be the collection of geometric sequents over $\Sigma$, $\tilde{S}$ the quotient of $\cal S$ by the relation of $\mathbb T$-equivalence, and $Sieves({\cal C}_{\mathbb T})$ the collection of sieves on the geometric syntactic category ${\cal C}_{\mathbb T}$ of $\mathbb T$.\\ 
Motivated by the proof of the duality theorem in section \ref{section_dualita}, let us define two corresponences ${\cal F}:{\cal S}\to Sieves({\cal C}_{\mathbb T})$ and ${\cal G}:Sieves({\cal C}_{\mathbb T}):\to \tilde{\cal S}$, as follows. 

Given a geometric sequent $\phi \vdash_{\vec{x}} \psi$ over $\Sigma$, we put ${\cal F}(\sigma)$ equal to principal sieve in ${\cal C}_{\mathbb T}$ generated by the monomorphism
\[  
\xymatrix {
\{\vec{x'}.\phi\wedge \psi\} \ar[rrr]^{[(\phi\wedge \psi\wedge \vec{x'}=\vec{x})]} & & & \{\vec{x}.\phi\}}
\]

Conversely, given a sieve $R$ in ${\cal C}_{\mathbb T}$, we put ${\cal G}(R)$ equal to the $\mathbb T$-equivalence class of any geometric sequent $ \psi \vdash_{\vec{y}} (\exists \vec{x})\theta$ such that $[\theta]$ is a monomorphism
\[  
\xymatrix {
\{\vec{x}.\phi\} \ar[r]^{[\theta]} & \{\vec{y}.\psi\}}
\]
in ${\cal C}_{\mathbb T}$ generating the principal sieve $\overline{R}^{J_{\mathbb T}}$ (cfr. Proposition \ref{regulargeometric}).\\
Applying the powerset functor to $\cal F$ and $\cal G$, we obtain maps of posets ${\mathscr{P}}({\cal F}):{\mathscr{P}}({\cal S}) \to {\mathscr{P}}(Sieves({\cal C}_{\mathbb T}))$ and ${\mathscr{P}}({\cal G}):{\mathscr{P}}(Sieves({\cal C}_{\mathbb T}))\to {\mathscr{P}}(\tilde{\cal S})$ (where the partial order on these sets is given by the inclusion). Concerning notation, we will write ${\cal F}(U)$ for ${\mathscr{P}}({\cal F})(U)$ and ${\cal G}(V)$ for ${\mathscr{P}}({\cal G})(V)$.\\
We have closure operators $\overline{(-)}^{\mathbb T}:{\mathscr{P}}({\cal S}) \to {\mathscr{P}}({\cal S})$ and $\overline{(-)}^{\cal T}:{\mathscr{P}}(Sieves({\cal C}_{\mathbb T})) \to {\mathscr{P}}(Sieves({\cal C}_{\mathbb T}))$ defined as follows: for a collection $U$ of geometric sequents over $\Sigma$, $\overline{U}^{\mathbb T}$ is the collection of geometric sequents $\sigma$ which are provable in ${\mathbb T}\cup U$ using geometric logic, while, for a collection $V$ of sieves in ${\cal C}_{\mathbb T}$, $\overline{V}^{\cal T}$ is the Grothendieck topology in ${\cal C}_{\mathbb T}$ generated by $J_{\mathbb T}$ and $V$ (i.e. the smallest Grothendieck topology $J$ on ${\cal C}_{\mathbb T}$ such that all the $J_{\mathbb T}$-covering sieves and the sieves in $V$ are $J$-covering); note that the relation of $\mathbb T$-equivalence on $\cal S$ is compatible with the closure operator $\overline{(-)}^{\mathbb T}$, that is we have a factorization $\overline{(-)}_{\tilde{S}}^{\mathbb T}:{\mathscr{P}}(\tilde{{\cal S}}) \to {\mathscr{P}}({\cal S})$ of $\overline{(-)}^{\mathbb T}:{\mathscr{P}}({\cal S}) \to {\mathscr{P}}({\cal S})$ through the image ${\mathscr{P}}({\cal S}) \to {\mathscr{P}}(\tilde{\cal S})$ via ${\mathscr{P}}$ of natural projection map ${\cal S} \to \tilde{\cal S}$.\\    
We note that the closed points with respect to these closure operators are respectively the closed quotients of $\mathbb T$ and the Grothendieck topologies $J$ on ${\cal C}_{\mathbb T}$ such that $J\supseteq J_{\mathbb T}$.\\
Let us define $F:{\mathscr{P}}({\cal S}) \to {\mathscr{P}}(Sieves({\cal C}_{\mathbb T}))$ as the composite $\overline{(-)}^{\cal T}\circ {\mathscr{P}}({\cal F})$ and $G:{\mathscr{P}}(Sieves({\cal C}_{\mathbb T})) \to {\mathscr{P}}({\cal S})$ as the composite $\overline{(-)}_{\tilde{\cal S}}^{\mathbb T}\circ {\mathscr{P}}({\cal G})$.\\

Given a collection $U$ of geometric sequents over $\Sigma$, we define ${\mathbb T}^{U}$ to be the collection of all the geometric sequents $\sigma$ over $\Sigma$ such that ${\cal F}(\sigma)$ belongs to $\overline{{\cal F}(U)}^{\cal T}$. Similarly, given a collection $V$ of sieves on ${\cal C}_{\mathbb T}$, we define $J_{V}$ to be the collection of sieves $R$ in $\cal C$ such that any sequent in ${\cal G}(R)$ is provable in ${\mathbb T}\cup {\cal G}(V)$ using geometric logic.\\
The following result shows that our maps ${\mathscr{P}}({\cal F})$ and ${\mathscr{P}}({\cal G})$ are compatible with respect to these closure operators, and that $F$ and $G$ are inverse to each other on the subsets of closed points, that is between the collection of closed quotients of $\mathbb T$ and the collection of Grothendieck topologies on ${\cal C}_{\mathbb T}$ which contain $J_{\mathbb T}$. In fact, given a quotient ${\mathbb T}'$ of $\mathbb T$, $F({\mathbb T}')=J_{{\mathbb T}'}^{\mathbb T}$ while for a Grothendieck topology $J\supseteq J_{\mathbb T}$, $G(J)={\mathbb T}^{J}$ (where the notations here are those of section \ref{section_dualita}). Thus this approach provides a different, entirely syntactic, way to arrive at the duality of Theorem \ref{dualita}.

\begin{theorem}\label{prooftheory}
With the above notation:\\
(i) For any $U\in {\mathscr{P}}({\cal S})$, ${\cal F}(\overline{U}^{\mathbb T})\subseteq \overline{{\cal F}(U)}^{\cal T}$;\\
(ii) For any $V \in {\mathscr{P}}(Sieves({\cal C}_{\mathbb T}))$, ${\cal G}(\overline{V}^{\cal T})\subseteq \overline{{\cal G}(V)}^{\mathbb T}$;\\
(iii) For any $U\in {\mathscr{P}}({\cal S})$, $G(F(U))=\overline{U}^{\mathbb T}={\mathbb T}^{U}$;\\
(iv) For any $V \in {\mathscr{P}}(Sieves({\cal C}_{\mathbb T}))$, $F(G(V))=\overline{V}^{\cal T}=J_{V}$.
\end{theorem}

\begin{proofs}
(i) We have to prove that, given $U:=\{\sigma_{i} \textrm{ | } i\in I\}\in {\mathscr{P}}({\cal S})$, if a geometric sequent $\sigma$ is provable in ${\mathbb T}\cup U$ using geometric logic, then ${\cal F}(\{\sigma\})$ belongs to $\overline{{\cal F}(U)}^{\cal T}$. Let us show this by induction on the complexity of a proof of $\sigma \equiv \phi \vdash_{\vec{x}} \psi$ in ${\mathbb T}\cup U$.\\
If $\sigma\in U$ then the thesis is clear.\\
If $\sigma$ belongs to $\mathbb T$ or, more generally, is provable in $\mathbb T$, then the morphism 
\[  
\xymatrix {
\{\vec{x'}.\phi\wedge \psi\} \ar[rrr]^{[(\phi\wedge \psi\wedge \vec{x'}=\vec{x})]} & & & \{\vec{x}.\phi\}}
\]
in ${\cal C}_{\mathbb T}$ is isomorphic to the identity morphism on $\{\vec{x}.\phi\}$, and hence it belongs to $\overline{{\cal F}(U)}^{\cal T}$ by the maximality axiom for Grothendieck topologies. Notice in particular that if $\sigma$ is an axiom of geometric logic then ${\cal F}(\sigma)$ belongs to $\overline{{\cal F}(U)}^{\cal T}$.\\
Now, let us verify that all the inference rules for geometric logic (described p. 830 \cite{El}) are `sound' with respect to the operation $\cal F$, that is if each of the premises $\sigma$ of an inference rule satisfies `${\cal F}(\sigma)$ belongs to $\overline{{\cal F}(U)}^{\cal T}$` then the conclusion $\sigma'$ of the rule also satisfies `${\cal F}(\sigma')$ belongs to $\overline{{\cal F}(U)}^{\cal T}$`.\\

\emph{Substitution rule}:
\[
\begin{array}{c}
(\phi \vdash_{\vec{x}} \psi)\\
\hline
(\phi[\vec{s}\slash \vec{x}] \vdash_{\vec{y}} \psi[\vec{s}\slash \vec{x}])  
\end{array}   
\]
where $\vec{y}$ is any string of variables including all the variables occurring in the string of terms $\vec{s}$.\\
We have to prove that if the sieve in ${\cal C}_{\mathbb T}$ generated by the single morphism 
\[  
\xymatrix {
\{\vec{x'}.\phi\wedge \psi\} \ar[rrr]^{[(\phi\wedge \psi\wedge \vec{x'}=\vec{x})]} & & & \{\vec{x}.\phi\}}
\]
is $\overline{{\cal F}(U)}^{\cal T}$-covering then the sieve generated by the single morphism
\[  
\xymatrix {
\{\vec{y'}.\phi[\vec{s}\slash \vec{x}]\wedge \psi[\vec{s}\slash \vec{x}]\} \ar[rrrr]^{[(\phi[\vec{s}\slash \vec{x}]\wedge \psi[\vec{s}\slash \vec{x}]\wedge \vec{y'}=\vec{y})]} & & & & \{\vec{y}.\phi[\vec{s}\slash \vec{x}]\}}
\]
is also $\overline{{\cal F}(U)}^{\cal T}$-covering.\\ 
For any geometric formula $\phi(\vec{x})$ and a term $s(\vec{y})$ over $\Sigma$, the diagram
\[  
\xymatrix {
\{\vec{y}. \phi[\vec{s}\slash \vec{x}]\} \ar[rrr]^{[(s(\vec{y})=\vec{x}) \wedge \phi]} \ar[d]^{[(\phi[\vec{s}\slash \vec{x}])[\vec{y'}\slash \vec{y}]\wedge \vec{y'}=\vec{y}]}  & & & \{\vec{x'}. \phi[\vec{x'}\slash \vec{x}]\} \ar[d]^{[(\phi\wedge \vec{x'}=\vec{x})]}\\
\{\vec{y}. \top\} \ar[rrr]^{[s(\vec{y})=\vec{x}]} & & & \{\vec{x}. \top\} }
\]\\ 

is a pullback in ${\cal C}_{\mathbb T}$. To prove this, let us first observe that if $\chi=(\exists \vec{y})\xi$ is a geometric formula in a context $\vec{x}$ such that the sequent 
\[
((\xi \wedge \xi[\vec{z}\slash \vec{y}]) \vdash_{\vec{x}, \vec{y}, \vec{z}} (\vec{y}=\vec{z}))
\]
is provable in $\mathbb T$ then the objects $\{\vec{x}. \chi\}$ and $\{\vec{x},\vec{y}. \xi\}$ are isomorphic in ${\cal C}_{\mathbb T}$. Indeed, it is an easy consequence of Lemma D1.4.4(i) \cite{El2} that the arrow
\[  
\xymatrix {
\{\vec{x}, \vec{y}.\xi\} \ar[rr]^{[(\xi \wedge (\vec{x'}=\vec{x}))]} & & \{\vec{x'}. \chi[\vec{x'}\slash \vec{x}]\}
}
\]
is an isomorphism.\\
Now, it immediately follows from the substitution axiom (and the equality axioms) that the sequent $(\exists \vec{x})((s(\vec{y})=\vec{x}) \wedge \phi(\vec{x})) \vdash_{\vec{y}} \phi[\vec{s}\slash \vec{x}]$ and its converse are provable in geometric logic.\\
So, in view of the construction of pullbacks given in section \ref{section_dualita} above, these two remarks together imply that our square is a pullback in ${\cal C}_{\mathbb T}$, as required.\\
From this we immediately deduce that the morphism

\[  
\xymatrix {
\{\vec{x'}.\phi[\vec{s}\slash \vec{x}]\wedge \psi[\vec{s}\slash \vec{x}]\} \ar[rrrr]^{[(\phi[\vec{s}\slash \vec{x}]\wedge \psi[\vec{s}\slash \vec{x}]\wedge \vec{x'}=\vec{x})]} & & & & \{\vec{x}.\phi[\vec{s}\slash \vec{x}]\}}
\]  

is (isomorphic to) the pullback in ${\cal C}_{\mathbb T}$ along $[(s(\vec{y})=\vec{x}) \wedge \phi]: \{\vec{y}. \phi[\vec{s}\slash \vec{x}]\} \to \{\vec{x'}. \phi[\vec{x'}\slash \vec{x}]\}$ of the morphism 
 
\[  
\xymatrix {
\{\vec{x'}.\phi\wedge \psi\} \ar[rrr]^{[(\phi\wedge \psi\wedge \vec{x'}=\vec{x})]} & & & \{\vec{x}.\phi\}}
\]
Now, for a Grothendieck topology $J$ on a category $\cal C$, it is always true that if the diagram 
\[  
\xymatrix {
a \ar[r] \ar[d]^{f'}  & b \ar[d]^{f} \\
d \ar[r]^{h} & c }
\] 
is a pullback in $\cal C$ then $(f)\in J(c)$ implies $(f')\in J(d)$. Indeed, by the universal property of the pullback, we have $(f')=h^{\ast}((f))$ and hence the thesis follows from the stability axiom for Grothendieck topologies.\\
This concludes the proof that the substitution rule is `sound' for the operation $\cal F$.\\

\emph{Cut rule:} 
\[
\begin{array}{c}
(\phi \vdash_{\vec{x}} \psi) (\psi \vdash_{\vec{x}} \chi)\\
\hline
(\phi \vdash_{\vec{x}} \chi)  
\end{array}   
\] 
We have to prove that if the sieves in ${\cal C}_{\mathbb T}$ respectively generated by the morphisms 
\[  
\xymatrix {
\{\vec{x'}.\phi\wedge \psi\} \ar[rrr]^{[(\phi\wedge \psi\wedge \vec{x'}=\vec{x})]} & & & \{\vec{x}.\phi\}}
\]
and
\[  
\xymatrix {
\{\vec{x'}.\psi\wedge \chi\} \ar[rrr]^{[(\psi\wedge \chi\wedge \vec{x'}=\vec{x})]} & & & \{\vec{x}.\psi\}}
\]
are $\overline{{\cal F}(U)}^{\cal T}$-covering then the sieve generated by the single morphism
\[  
\xymatrix {
\{\vec{x'}.\phi\wedge \chi\} \ar[rrr]^{[(\phi\wedge \chi\wedge \vec{x'}=\vec{x})]} & & & \{\vec{x}.\phi\}}
\]
is also $\overline{{\cal F}(U)}^{\cal T}$-covering.\\ 
The diagrams
\[  
\xymatrix {
\{\vec{x'''}.\phi\wedge \psi \wedge \chi\} \ar[rrr]^{[(\phi\wedge \psi \wedge \chi\wedge \vec{x'''}=\vec{x''})]} \ar[dd]^{[(\phi\wedge \psi \wedge \chi\wedge \vec{x'''}=\vec{x'})]} & & & \{\vec{x''}.\phi\wedge \chi\} \ar[dd]^{[(\phi\wedge \chi\wedge \vec{x''}=\vec{x})]} \\
& & & \\
\{\vec{x'}.\phi\wedge \psi\} \ar[rrr]^{[(\phi\wedge \psi\wedge \vec{x'}=\vec{x})]} & & & \{\vec{x}.\phi\}}
\] 
and 
\[  
\xymatrix {
\{\vec{x'''}.\phi\wedge \psi \wedge \chi\} \ar[rrr]^{[(\phi\wedge \psi \wedge \chi\wedge \vec{x'''}=\vec{x''})]} \ar[dd]^{[(\phi\wedge \psi \wedge \chi\wedge \vec{x'''}=\vec{x'})]} & & & \{\vec{x''}.\psi\wedge \chi\} \ar[dd]^{[(\psi\wedge \chi\wedge \vec{x''}=\vec{x})]} \\
& & & \\
\{\vec{x'}.\phi\wedge \psi\} \ar[rrr]^{[(\phi\wedge \psi\wedge \vec{x'}=\vec{x})]} & & & \{\vec{x}.\psi\}}
\] 
are clearly pullback squares in ${\cal C}_{\mathbb T}$.\\
By the stability axiom for Grothendieck topologies, the sieve generated by the morphism
\[
\xymatrix {
\{\vec{x'''}.\phi\wedge \psi \wedge \chi\} \ar[rrr]^{[(\phi\wedge \psi \wedge \chi\wedge \vec{x'''}=\vec{x'})]} & & & \{\vec{x'}.\phi\wedge \psi\}}
\]
is $\overline{{\cal F}(U)}^{\cal T}$-covering, since it is the pullback of the ($\overline{{\cal F}(U)}^{\cal T}$-covering) sieve generated by the morphism 
\[  
\xymatrix {
\{\vec{x'}.\psi\wedge \chi\} \ar[rrr]^{[(\psi\wedge \chi\wedge \vec{x'}=\vec{x})]} & & & \{\vec{x}.\psi\}}
\]
along the arrow
\[  
\xymatrix {
\{\vec{x'}.\phi\wedge \psi\} \ar[rrr]^{[(\phi\wedge \psi\wedge \vec{x'}=\vec{x})]} & & & \{\vec{x}.\psi\}}
\]  
So, since the sieve generated by the morphism
\[  
\xymatrix {
\{\vec{x'}.\phi\wedge \psi\} \ar[rrr]^{[(\phi\wedge \psi\wedge \vec{x'}=\vec{x})]} & & & \{\vec{x}.\phi\}}
\]
is $\overline{{\cal F}(U)}^{\cal T}$-covering, we conclude, by the transitivity axiom for Grothendieck topologies and the fact that the first square above is a pullback, that the sieve generated by the morphism
\[  
\xymatrix {
\{\vec{x'}.\phi\wedge \chi\} \ar[rrr]^{[(\phi\wedge \chi\wedge \vec{x'}=\vec{x})]} & & & \{\vec{x}.\phi\}}
\]
is $\overline{{\cal F}(U)}^{\cal T}$-covering, as required.\\

\emph{Rule for finite conjunction:}
\[
\begin{array}{c}
(\phi \vdash_{\vec{x}} \psi) (\phi \vdash_{\vec{x}} \chi)\\
\hline
(\phi \vdash_{\vec{x}} (\psi \wedge \chi))  
\end{array}   
\] 
We have to prove that if the sieves in ${\cal C}_{\mathbb T}$ respectively generated by the morphisms 
\[  
\xymatrix {
\{\vec{x'}.\phi\wedge \psi\} \ar[rrr]^{[(\phi\wedge \psi\wedge \vec{x'}=\vec{x})]} & & & \{\vec{x}.\phi\}}
\]
and
\[  
\xymatrix {
\{\vec{x'}.\phi\wedge \chi\} \ar[rrr]^{[(\phi\wedge \chi\wedge \vec{x'}=\vec{x})]} & & & \{\vec{x}.\phi\}}
\]
are $\overline{{\cal F}(U)}^{\cal T}$-covering then the sieve generated by the single morphism
\[  
\xymatrix {
\{\vec{x'}.\phi\wedge (\psi \wedge \chi)\} \ar[rrrr]^{[(\phi\wedge (\psi \wedge \chi)\wedge \vec{x'}=\vec{x})]} & & & & \{\vec{x}.\phi\}}
\]
is also $\overline{{\cal F}(U)}^{\cal T}$-covering.\\
We observed above that the diagram
\[  
\xymatrix {
\{\vec{x'''}.\phi\wedge \psi \wedge \chi\} \ar[rrr]^{[(\phi\wedge \psi \wedge \chi\wedge \vec{x'''}=\vec{x''})]} \ar[dd]^{[(\phi\wedge \psi \wedge \chi\wedge \vec{x'''}=\vec{x'})]} & & & \{\vec{x''}.\phi\wedge \chi\} \ar[dd]^{[(\phi\wedge \chi\wedge \vec{x''}=\vec{x})]} \\
& & & \\
\{\vec{x'}.\phi\wedge \psi\} \ar[rrr]^{[(\phi\wedge \psi\wedge \vec{x'}=\vec{x})]} & & & \{\vec{x}.\phi\}}
\] 
is a pullback in ${\cal C}_{\mathbb T}$. Thus, by the stability axiom for Grothendieck topologies, the sieve generated by the arrow
\[
\xymatrix {
\{\vec{x'''}.\phi\wedge \psi \wedge \chi\} \ar[rrr]^{[(\phi\wedge \psi \wedge \chi\wedge \vec{x'''}=\vec{x'})]} & & & \{\vec{x'}.\phi\wedge \psi\}}
\]
is $\overline{{\cal F}(U)}^{\cal T}$-covering, since it is the pullback of the ($\overline{{\cal F}(U)}^{\cal T}$-covering) sieve generated by the arrow 
\[  
\xymatrix {
\{\vec{x'}.\phi\wedge \chi\} \ar[rrr]^{[(\phi\wedge \chi\wedge \vec{x'}=\vec{x})]} & & & \{\vec{x}.\phi\}}
\]
along the arrow
\[  
\xymatrix {
\{\vec{x'}.\phi\wedge \psi\} \ar[rrr]^{[(\phi\wedge \psi\wedge \vec{x'}=\vec{x})]} & & & \{\vec{x}.\phi\}}
\]   
But the sieve generated by the arrow 
\[  
\xymatrix {
\{\vec{x'}.\phi\wedge \psi\} \ar[rrr]^{[(\phi\wedge \psi\wedge \vec{x'}=\vec{x})]} & & & \{\vec{x}.\phi\}}
\]   
is $\overline{{\cal F}(U)}^{\cal T}$-covering and hence, since the arrow
\[  
\xymatrix {
\{\vec{x'}.\phi\wedge (\psi \wedge \chi)\} \ar[rrrr]^{[(\phi\wedge (\psi \wedge \chi)\wedge \vec{x'}=\vec{x})]} & & & & \{\vec{x}.\phi\}}
\]
is equal to the composite of 
\[  
\xymatrix {
\{\vec{x'}.\phi\wedge \psi\} \ar[rrr]^{[(\phi\wedge \psi\wedge \vec{x'}=\vec{x})]} & & & \{\vec{x}.\phi\}}
\] 
and
\[
\xymatrix {
\{\vec{x'''}.\phi\wedge \psi \wedge \chi\} \ar[rrr]^{[(\phi\wedge \psi \wedge \chi\wedge \vec{x'''}=\vec{x'})]} & & & \{\vec{x'}.\phi\wedge \psi\}}
\]
we deduce, by property (iv) in Definition \ref{def2}, that the sieve generated by the arrow  
\[  
\xymatrix {
\{\vec{x'}.\phi\wedge (\psi \wedge \chi)\} \ar[rrrr]^{[(\phi\wedge (\psi \wedge \chi)\wedge \vec{x'}=\vec{x})]} & & & & \{\vec{x}.\phi\}}
\] 
is $\overline{{\cal F}(U)}^{\cal T}$-covering, as required.\\

\emph{Rule for infinitary disjunction:}
\[
\begin{array}{c}
\{(\phi_{i} \vdash_{\vec{x}} \chi) \textrm{ | } i\in I \}\\
\hline\\
(\mathbin{\mathop{\textrm{\huge $\vee$}}\limits_{i\in I}}\phi_{i}  \vdash_{\vec{x}} \chi)  
\end{array}   
\] 
We have to prove that if each of the sieves in ${\cal C}_{\mathbb T}$ respectively generated by the single arrow 
\[  
\xymatrix {
\{\vec{x'}.\phi_{i}\wedge \chi\} \ar[rrr]^{[(\phi_{i}\wedge \chi\wedge \vec{x'}=\vec{x})]} & & & \{\vec{x}.\phi_{i}\}}
\]
as $i$ varies in $I$ is $\overline{{\cal F}(U)}^{\cal T}$-covering then the sieve generated by the single morphism
\[  
\xymatrix {
\{\vec{x'}.(\mathbin{\mathop{\textrm{\huge $\vee$}}\limits_{i\in I}}\phi_{i})\wedge \chi\} \ar[rrrr]^{[((\mathbin{\mathop{\textrm{\huge $\vee$}}\limits_{i\in I}}\phi_{i}) \wedge \chi)\wedge \vec{x'}=\vec{x})]} & & & & \{\vec{x}.\mathbin{\mathop{\textrm{\huge $\vee$}}\limits_{i\in I}}\phi_{i}\}}
\]
is also $\overline{{\cal F}(U)}^{\cal T}$-covering.\\
The sieve on $\{\vec{x}.\mathbin{\mathop{\textrm{\huge $\vee$}}\limits_{i\in I}}\phi_{i}\}$ generated by the arrows
\[  
\xymatrix {
j_{i}:=\{\vec{x'}.\phi_{i}\} \ar[rr]^{[\phi_{i} \wedge \vec{x'}=\vec{x}]} & & \{\vec{x}.\mathbin{\mathop{\textrm{\huge $\vee$}}\limits_{i\in I}}\phi_{i}\}}
\]
as $i$ varies in $I$ is $\overline{{\cal F}(U)}^{\cal T}$-covering by definition of $J_{\mathbb T}$, since ${\cal T}^{\mathbb T}_{U} \supseteq J_{\mathbb T}$.\\
Now, for each $i\in I$ the diagram
\[  
\xymatrix {
\{\vec{x'''}.\phi_{i}\wedge \chi\} \ar[rrr]^{[(\phi_{i} \wedge \chi \wedge \vec{x'''}=\vec{x''})]} \ar[dd]^{[(\phi_{i} \wedge \chi\wedge \vec{x'''}=\vec{x'})]} & & & \{\vec{x''}.(\mathbin{\mathop{\textrm{\huge $\vee$}}\limits_{i\in I}}\phi_{i})\wedge \chi\} \ar[dd]^{[(\mathbin{\mathop{\textrm{\huge $\vee$}}\limits_{i\in I}}\phi_{i})\wedge \chi \wedge \vec{x''}=\vec{x}\}]} \\
& & & \\
\{\vec{x'}.\phi_{i}\} \ar[rrr]^{j_{i}} & & & \{\vec{x}.\mathbin{\mathop{\textrm{\huge $\vee$}}\limits_{i\in I}}\phi_{i}\}\}}
\] 
is a pullback in ${\cal C}_{\mathbb T}$. Our thesis then follows from the transitivity axiom for Grothendieck topologies.\\

\emph{Rules for existential quantification:}
\[
\begin{array}{c}
(\phi \vdash_{\vec{x}, \vec{y}} \psi)\\
\hline\hline
((\exists \vec{y})\phi \vdash_{\vec{x}} \psi)  
\end{array}   
\] 
where $\vec{y}$ is not free in $\psi$.\\
We have to prove that the sieve in ${\cal C}_{\mathbb T}$ generated by the single morphism 
\[  
\xymatrix {
\{\vec{x'}, \vec{y'}.\phi\wedge \psi\} \ar[rrrr]^{[(\phi\wedge \psi\wedge \vec{x'}=\vec{x} \wedge \vec{y'}=\vec{y})]} & & & & \{\vec{x}, \vec{y}.\phi\}}
\]
is $\overline{{\cal F}(U)}^{\cal T}$-covering if and only if the sieve generated by the single morphism
\[  
\xymatrix {
\{\vec{x'}. ((\exists \vec{y})\phi) \wedge \psi \} \ar[rrrr]^{[((\exists \vec{y})\phi) \wedge \wedge \vec{x'}=\vec{x}]} & & & & \{\vec{x'}. (\exists \vec{y})\phi \}}
\]
is $\overline{{\cal F}(U)}^{\cal T}$-covering.\\
The diagram
\[  
\xymatrix {
\{\vec{x'''}, \vec{y'''}.(\phi\wedge \psi)[\vec{x'''}\slash \vec{x}, \vec{y'''}\slash \vec{y}]\} \ar[rrrr]^{[(\phi \wedge \psi \wedge \vec{x'''}=\vec{x''} \wedge \vec{y'''}=\vec{y''})]} \ar[dd]^{[(\phi_{i} \wedge \chi\wedge \vec{x'''}=\vec{x'})]} & & & & \{\vec{x''}, \vec{y''}.\phi[\vec{x'''}\slash \vec{x}, \vec{y'''}\slash \vec{y}]\} \ar[dd]^{[\phi \wedge \vec{x''}=\vec{x}\}]} \\
& & & & \\
\{\vec{x'}.((\exists \vec{y})\phi \wedge \psi)[\vec{x'}\slash \vec{x}]\} \ar[rrrr]^{[((\exists \vec{y})\phi)\wedge \psi \wedge \vec{x'}=\vec{x}]} & & & & \{\vec{x}.(\exists \vec{y})\phi\}}
\] 
is a pullback in ${\cal C}_{\mathbb T}$. Indeed, this easily follows from the construction of pullbacks in ${\cal C}_{\mathbb T}$ given in section \ref{section_dualita} by invoking the rules for existential quantification, as in the proof for the substitution rule.\\
Now, the `if' part of our thesis clearly follows from the stability axiom for Grothendieck topologies. It remains to prove the `only if' part. To this end, notice that the arrow
\[  
\xymatrix {
\{\vec{x''}, \vec{y''}.\phi[\vec{x''}\slash \vec{x}, \vec{y''}\slash \vec{y}]\} \ar[rrr]^{[\phi \wedge \vec{x''}=\vec{x}\}]} & & & \{\vec{x}.(\exists \vec{y})\phi\}}
\]
is a cover in ${\cal C}_{\mathbb T}$; so the sieve generated by it is $\overline{{\cal F}(U)}^{\cal T}$-covering by definition of $J_{\mathbb T}$, since ${\cal T}^{\mathbb T}_{U} \supseteq J_{\mathbb T}$. Hence, by the commutativity of the square above, the sieve generated by the arrow
\[  
\xymatrix {
\{\vec{x'}. ((\exists \vec{y})\phi) \wedge \psi \} \ar[rrrr]^{[((\exists \vec{y})\phi) \wedge \wedge \vec{x'}=\vec{x}]} & & & & \{\vec{x'}. (\exists \vec{y})\phi \}}
\]
is $\overline{{\cal F}(U)}^{\cal T}$-covering by properties (ii) and (iv) in Definition \ref{def2}.\\ 
This completes the proof of part (i) of the theorem.\\
(ii) We have to prove that, given $V \in {\mathscr{P}}(Sieves({\cal C}_{\mathbb T}))$, if a sieve $R$ belongs to $\overline{V}^{\cal T}$ then any sequent in ${\cal G}(R)$ is provable in ${\mathbb T}\cup {\cal G}(V)$ using geometric logic, that is $J_{V}\supseteq \overline{V}^{\cal T}$. In fact, we will prove that $J_{V}$ is a Grothendieck topology containing $J_{\mathbb T}$ and all the sieves in $V$ as its covering sieves; this will clearly imply our thesis.\\
Clearly, by definition of $J_{V}$, the sieves in $V$ belong to $J_{V}$, and if $R$ is a $J_{\mathbb T}$-covering sieve then, by definition of $J_{\mathbb T}$, any sequent in ${\cal G}(R)$ is provable in $\mathbb T$, so that $R$ belongs to $J_{V}$. To prove that $J_{V}$ is a Grothendieck topology, we use Definition \ref{def2}. Property (i) is obvious, and property (ii) easily follows from the cut rule in geometric logic. Property (iii) follows from the proof of Theorem \ref{dualita}. It remains to prove property (iv). Since ${\cal G}(R)={\cal G}(\overline{R}^{J_{\mathbb T}})$ for any sieve $R$ in ${\cal C}_{\mathbb T}$ then, by Proposition \ref{regulargeometric} and Remark \ref{pbstable}, it suffices to prove that for any sieve $S$ generated by a monomorphism $m:d\to c$ and any sieve $T$ in ${\cal C}_{\mathbb T}$ on $d$, if both $S$ and $T$ are $J_{V}$-covering then $S\ast T$ is $J_{V}$-covering. Now, in view of the equality ${\cal G}(R)={\cal G}(\overline{R}^{J_{\mathbb T}})$, our claim easily follows from the cut rule in geometric logic, by using Proposition \ref{regulargeometric}, Remark \ref{pbstable} and Remark \ref{def2closure}.\\
This concludes the proof of part (ii) of the theorem.\\

(iii) Let us begin by proving that $\overline{{\cal G}({\cal F}(U))}^{\mathbb T}=\overline{U}^{\mathbb T}$. Note that $\overline{{\cal G}({\cal F}(U))}^{\mathbb T}$ is the collection of sequents of the form ${\cal G}({\cal F}(\sigma))$ as $\sigma$ varies in $U$. If $\sigma$ is $\phi \vdash_{\vec{x}} \psi$ then ${\cal G}({\cal F}(\sigma))$ is the $\mathbb T$-equivalence class of the sequent $\phi \vdash_{\vec{x}} \phi \wedge \psi$; but this sequent is clearly $\mathbb T$-equivalent to $\sigma$, and hence $\overline{{\cal G}({\cal F}(U))}^{\mathbb T}=\overline{U}^{\mathbb T}$, as required.\\
We have $G(F(U))=\overline{{\cal G}(F(U))}^{\mathbb T}=\overline{{\cal G}(\overline{{\cal F}(U)}^{\cal T})}^{\mathbb T}=\overline{\overline{{\cal G}({\cal F}(U))}^{\mathbb T}}^{\mathbb T}=\overline{{\cal G}({\cal F}(U))}^{\mathbb T}=\overline{U}^{\mathbb T}$, where the central equality follows from part (i) of the theorem. This proves the first of the two equalities in part (iii) of the theorem; it remains to show that $\overline{U}^{\mathbb T}={\mathbb T}^{U}$ holds. The inclusion $\overline{U}^{\mathbb T} \subseteq {\mathbb T}^{U}$ follows from part (i) of the theorem, while the other one one follows as a consequence of the first equality in part (iii) and from part (i) of the theorem: if $\sigma\in {\mathbb T}^{U}$ then $\sigma \in \overline{\{\sigma \}}^{\mathbb T}=\overline{{\cal G}({\cal F}(\{\sigma\}))}^{\mathbb T} \subseteq \overline{{\cal G}({\cal F}(U))}^{\mathbb T}=\overline{U}^{\mathbb T}$.\\ 

(iv) Let us begin by proving that $\overline{{\cal F}({\cal G}(V))}^{\cal T}=\overline{V}^{\cal T}$. Now, ${\cal F}({\cal G}(V))^{\cal T}$ is the collection of sieves of the form ${\cal F}({\cal G}(R))$ as $R$ varies in $V$, and it is immediate to see that ${\cal F}({\cal G}(R))=\overline{R}^{J_{\mathbb T}}$; hence our claim follows from Proposition \ref{gengrot}.\\
Now, by using the fact that $\overline{{\cal F}({\cal G}(V))}^{\cal T}=\overline{V}^{\cal T}$, one can prove the required equalities as in the proof of part (iii) of the theorem, with the only difference that part (ii) play the role of part (i) here.\\     
\end{proofs}  

Let $\mathbb T$ be a geometric theory over a signature $\Sigma$. Given a quotient ${\mathbb T}'$ of $\mathbb T$, let $J_{{\mathbb T}'}^{\mathbb T}$ be the associated $\mathbb T$-topology of ${\mathbb T}'$. Then the equalities $\overline{U}^{\mathbb T}={\mathbb T}^{U}$ and $\overline{V}^{\cal T}=J_{V}$ in Theorem \ref{prooftheory} give the following equivalences:\\
(1) for any sieve $R \in Sieves({\cal C}_{\mathbb T})$, $R\in J_{{\mathbb T}'}^{\mathbb T}$ if and only if any sequent in ${\cal G}(R)$ is provable in ${\mathbb T}'$;\\
(2) for any geometric sequent $\sigma$ over $\Sigma$, $\sigma$ is provable in ${\mathbb T}'$ if and only if ${\cal F}(\sigma)$ is $J_{{\mathbb T}'}^{\mathbb T}$-covering.\\
In particular, we obtain the following characterization of the syntactic topology $J_{\mathbb T}$ on ${\cal C}_{\mathbb T}$: a sieve $R$ is $J_{\mathbb T}$-covering if and only if any sieve in ${\cal G}(R)$ is provable in $\mathbb T$.  

\begin{rmk}
\emph{Let us briefly consider how much of Theorem \ref{prooftheory} survives for smaller fragments of geometric logic, e.g. cartesian, regular, or coherent logic. If $\mathbb T$ is a cartesian (resp. regular, coherent) theory over $\Sigma$, one can define exactly as above an assignment ${\cal F}:{\cal S}\to Sieves({\cal C}_{\mathbb T})$, where ${\cal S}$ is the collection of cartesian (resp. regular, coherent) sequents over $\Sigma$ and ${\cal C}_{\mathbb T}$ is the cartesian (resp. regular, coherent) syntactic category of the theory $\mathbb T$. Accordingly, the closure operator on ${\mathscr{P}}({\cal S})$ sends a collection $U$ of sequents in ${\cal S}$ to the collection of cartesian (resp. regular, coherent) sequents over $\Sigma$ which are derivable from $U \cup {\mathbb T}$ by using cartesian (resp. regular, coherent) logic, and it is immediate to see that the proof of part (i) of the theorem continues to hold. On the contrary, no assignment $\cal G$ with values in the class of cartesian (resp. regular, coherent) sequents over $\Sigma$ can be defined, since one should restrict to sieves generated by a monomorphism (resp. a single arrow, a finite number of arrows); however, if we consider $\cal G$ to take values in the class of geometric sequents over $\Sigma$ as in the geometric case then we still still have that part (ii) of the theorem holds and that for any presieve $V$ in the relevant cartesian (resp. regular, coherent) category ${\cal C}_{\mathbb T}^{\textrm{cart}}$ (resp. ${\cal C}_{\mathbb T}^{\textrm{reg}}$, ${\cal C}_{\mathbb T}^{\textrm{coh}}$) the theory ${\cal G}(\overline{V}^{\cal T})$ is classified by the topos $\Sh({\cal C}_{\mathbb T}^{\textrm{cart}}, \overline{V}^{\cal T})$ (cfr. the proof of Theorem \ref{dualita}).}  
\end{rmk}

\begin{rmk}\label{generazione}
\emph{Given a closed geometric quotient ${\mathbb T}'$ of ${\mathbb T}$, it is natural to look for axiomatizations of ${\mathbb T}'$ over ${\mathbb T}$ which are as simple as possible; this translates, via the duality theorem, into the problem of finding a `simple as possible' set of generators for the associated Grothendieck topology $J_{{\mathbb T}'}^{\mathbb T}$ over $J_{\mathbb T}$; in fact, if a collection $V$ of presieves in ${\cal C}_{\mathbb T}$ generates a Grothendieck topology $J$, then, by Theorem \ref{prooftheory}(ii), ${\mathbb T}^{J}$ is axiomatized over $\mathbb T$ by the collection of sequents in  ${\cal G}(V)$ (note that, conversely, if a collection $U$ of geometric sequents axiomatizes a quotient ${\mathbb T}'$ then, by Theorem \ref{prooftheory}(i) the collection of presieves ${\cal F}(U)$ generates over $J_{\mathbb T}$ the Grothendieck topology $J_{{\mathbb T}'}^{\mathbb T}$).\\
For example, one may ask if $\mathbb T$ can be axiomatized over $\mathbb T$ by geometric sequents of the form $\top \vdash_{\vec{x}} \phi$; this correspond to requiring that $J_{{\mathbb T}'}^{\mathbb T}$ should be generated over $J_{{\mathbb T}}$ by a collection of principal sieves generated by subobjects of objects of the form $\{\vec{x}. \top\}$; two notable classes of theories with this property are the classes of Booleanizations and DeMorganizations of a given geometric theory (cfr. \cite{OC3}).\\
It is often the case that, by adopting the point of view of Grothendieck topologies, one gets interesting insights at the level of theories. To give an illustration of this, let us discuss the case of the Booleanization ${\mathbb T}'$ of a geometric theory $\mathbb T$. Given a Heyting category $\cal C$, let us denote by $\tilde{\cal C}$ its full subcategory on the non-zero objects. Since $J_{\mathbb T}$ is subcanonical then $\tilde{{\cal C}_{\mathbb T}}$ is $J_{\mathbb T}$-dense in ${\cal C}_{\mathbb T}$, and the induced Grothendieck topology $J_{\mathbb T}|\tilde{{\cal C}_{\mathbb T}}$ is dense (as a Grothendieck topology on $\tilde{{\cal C}_{\mathbb T}}$); hence $J_{{\mathbb T}'}^{\mathbb T}$ is generated over $J_{\mathbb T}$ by the sieves generated in ${\cal C}_{\mathbb T}$ by the principal stably non-empty sieves in $\tilde{{\cal C}_{\mathbb T}}$. Now, given a Heyting category $\cal C$ and monomorphisms $f:d\mono c'$ and $g:c'\mono c$ in $\tilde{{\cal C}}$, it is immediate to see that if $(f)$, regarded as a sieve in $\tilde{\cal C}$, is stably non-empty then $(f)=g^{\ast}((g \circ f)\cup \neg (g \circ f))$ where $\cup$ and $\neg$ respectively denote the union and pseudocomplementation in the Heyting algebra $\Sub_{\cal C}(c)$, and hence $(f)$ is the pullback of a stably non-empty sieve in $\tilde{\cal C}$ on $c$. Therefore, since every object in ${\cal C}_{\mathbb T}$ has a monomorphism to an object of the form $\{\vec{x}. \top\}$, we deduce that $J_{{\mathbb T}'}^{\mathbb T}$ is generated over $J_{\mathbb T}$ by a collection of principal sieves generated by subobjects of objects of the form $\{\vec{x}. \top\}$, as required.}
\end{rmk}

\section{The lattice structure}\label{latticestructure}

In this section we study the structure of the lattice of subtoposes of a given Grothendieck topos. It is well-known that this lattice, endowed with the obvious order relation given by the inclusion of subtoposes, is a coHeyting algebra (see for example section A4.5 \cite{El}). Our aim is to describe this structure in terms of Grothendieck topologies and later of theories, in view of Theorem \ref{dualita}. In fact, as we see below, it suffices to describe the lattice operations on the collection of subtoposes of a given presheaf toposes.\\   
Given an Heyting algebra $H$ and an element $a\in H$, the collection $\princfil{a}$ of all the elements $h\in H$ such that $h\geq a$ is closed under the operations of conjunction, disjunction and Heyting implication and it is (therefore) an Heyting algebra with respect to these operations. Indeed, the assertion about the conjunction and disjunction is obvious, while the fact that $b\imp c$ is in $\princfil{a}$ if $b$ and $c$ are follows from the inequality $c\leq (b\imp c)$.\\
This remark allows us to restrict our attention to the case of subtoposes of a presheaf topos in order to describe the effect of the operations of union, intersection and coHeyting implication on a pair of subtoposes of a given Grothendieck topos; indeed, the union (resp. intersection, coHeyting implication) of two subtoposes of $\Sh({\cal C}, J)$ is the same as the union (resp. intersection, coHeyting implication) of them in the coHeyting algebra of subtoposes of $[{\cal C}^{\textrm{op}}, \Set]$, since the order-relation in the former lattice is clearly the restriction of the order relation in the second (in both cases the order being the dual of the relation `to be a subtopos of').\\

\subsection{The lattice operations on Grothendieck topologies}

Let ${\cal E}=[{\cal C}^{\textrm{op}}, \Set]$ be a presheaf topos, with subobject classifier $\Omega$. Recall that $\Omega: {\cal C}^{\textrm{op}} \to \Set$ is defined by:\\
$\Omega(c)=\textrm{\{$R$ | $R$ is a sieve on $c$\}}$ (for any object $c\in \cal C$),\\
$\Omega(f)=f^{\ast}(-)$ (for any arrow $f$ in $\cal C$),\\
where $f^{\ast}(-)$ denotes the operation of pullback of sieves in $\cal C$ along $f$.\\
We know from Theorem 1 p. 233 \cite{MM} that, given a small category $\cal C$, the Grothendieck topologies $J$ on $\cal C$ correspond exactly to local operators on the topos $[{\cal C}^{\textrm{op}}, \Set]$; this correspondence, to which we refer as $(\ast)$, sends a local operator $j:\Omega\to \Omega$ to the subobject $J\mono \Omega$ which it classifies, that is to the Grothendieck topology $J$ on $\cal C$ defined by: $S\in J(c)$ if and only if $j(c)(S)=M_{c}$, and conversely a subobject $J\in \Omega$ to the map $j:\Omega\to \Omega$ which classifies it.\\
Let us recall from \cite{MM} (formula (7) p. 38) that, given a subobject $A\mono \Omega$, its characteristic map $\chi_{A}:\Omega\to \Omega$ is given by the formula:\\
\[
\chi_{A}(c)(S)=\{f:d\to c \textrm{ | } f^{\ast}(S)\in A(d)\}
\]   
Let us now give an explicit description of the internal Heyting operations $\wedge,\vee, \imp:\Omega\to \Omega$ on our presheaf topos $\cal E$ (defined for example in the proof of Lemma A1.6.3 \cite{El}); this will be convenient for our purposes.\\
The internal conjunction map $\wedge:\Omega\times \Omega\to \Omega$ is the classifying map of the subobject $(\top,\top):1\mono \Omega\times \Omega$ , so we immediately get the following expression:\\
\[
\wedge(c)(S,T)=S\cap T
\]
for any object $c\in {\cal C}$ and sieves $S$ and $T$ on $c$.\\
The internal disjunction map $\vee:\Omega \times \Omega\to \Omega$ is the classifying map of the union of subobjects $\pi_{1}^{\ast}(\top)$ and $\pi_{2}^{\ast}(\top)$, where $\pi_{1}$ and $\pi_{2}$ are the two product projections $\Omega\times \Omega\to \Omega$ so we get
\[
\vee(c)(S,T)=\{f:d\to c \textrm{ | } f^{\ast}(S)\cup f^{\ast}(T)=M_{d}\}
\]
for any object $c\in {\cal C}$ and sieves $S$ and $T$ on $c$.\\
   
The internal implication map $\imp:\Omega \times \Omega\to \Omega$ is the classifying map of the equalizer $\Omega_{1}\mono \Omega \times \Omega$ of $\wedge$ and $\pi_{1}$ so we obtain
\[
\imp(c)(S,T)=\{f:d\to c \textrm{ | } f^{\ast}(S)\subseteq f^{\ast}(T)\}
\] 
for any object $c\in {\cal C}$ and sieves $S$ and $T$ on $c$.\\
 
It is immediate to check that the order relation between local operators on $\cal E$ given by the opposite of the natural order between subtoposes transfers via $(\ast)$ to the following order between Grothendieck topologies on $\cal C$: $J\leq J'$ if and only if for every $c\in {\cal C}$, $J(c)\subseteq J'(c)$ i.e. every $J$-covering sieve is $J'$-covering. Hence, from $(\ast)$ we deduce that the relation $\leq$ defines an Heyting algebra structure on the collection of Grothendieck topologies on the category $\cal C$; in particular, for any two Grothendieck topologies $J$ and $J'$ on $\cal C$, there exists a meet $J\wedge J'$, a join $J\vee J'$ and a Heyting implication $J\imp J'$. We note that the bottom element of this lattice is the Grothendieck topology $\bot$ on $\cal C$ given by $\bot(c)=\{M_{c}\}$ for every $c\in {\cal C}$, while the top element is the topology $\top$ defined by: $\top(c)=\{S \textrm{ | $S$ sieve on } c \}$, for every $c\in {\cal C}$.\\
We can easily get an explicit expression for $J\wedge J'$: $S\in J\wedge J'(c)$ if and only if $S\in J(c)$ and $S\in J'(c)$; indeed, the class of Grothendieck topologies is clearly closed under intersection. The join $J\vee J'$ is the smallest Grothendieck topology $K$ such that $J\leq K$ and $J'\leq K$, so it is the Grothendieck topology generated by the collection of sieves which are either $J$-covering or $J'$-covering. In order to get a more explicit description of it, and also of the Heyting implication between Grothendieck topologies, we specialize A. Joyal's theory as it is described in A4.5 \cite{El} to the context of Grothendieck toposes; this will lead in particular to an explicit description of the Grothendieck topology generated by a family of sieves which is stable under pullbacks.\\
First, let us make explicit in terms of the category $\cal C$ the Galois connection from $\Sub_{\cal E}(\Omega)$ to itself given by the mappings $D\to D^{r}$ and $D\to D^{l}$ decribed p. 213 \cite{El}.\\
Given a subobject $D\mono \Omega$, $D^{r}\mono \Omega$ and $D^{l}\mono \Omega$ are defined to be respectively
\[
\forall_{\pi_{2}}((\pi^{\ast}_{1}(D)\imp \Theta)\mono \Omega
\]
and
\[
\forall_{\pi_{1}}((\pi^{\ast}_{2}(D)\imp \Theta)\mono \Omega
\]
where $\pi_{1}$ and $\pi_{2}$ are the two product projections $\Omega \times \Omega \to \Omega$, $\pi^{\ast}_{1}$ and $\pi^{\ast}_{2}$ are the pullback functors $\Sub(\Omega)\to \Sub(\Omega \times \Omega)$ respectively along $\pi_{1}$ and $\pi_{2}$, and $\Theta \mono \Omega \times \Omega$ is the equalizer of $\pi_{2}, \imp: \Omega \times \Omega \to \Omega$.\\
First, note that the subobjects of $\Omega$ can be identified with collections of sieves in $\cal C$ which are stable under pullback; in fact, from now on we will use this identification.\\
From the formulas above, we get the following expression for $\Theta$:
\[
\begin{array}{ccl}
\Theta(c) & = & \{(S, T)  \textrm{ | $S$ and $T$ are sieves on $c$ s.t. for all $f:d\to c$,}\\
& & \textrm{$f^{\ast}(S)\subseteq f^{\ast}(T)$ implies $f\in T$}\}
\end{array}
\]
for any object $c\in {\cal C}$.\\
Now, by using formula (7) p. 146 \cite{MM}, we obtain:
\[
\begin{array}{ccl}
\pi^{\ast}_{1}(D)\imp \Theta & = & \{(S, T) \textrm{ | $S$ and $T$ are sieves on $c$ s.t. for all $f:d\to c$, }\\
& & \textrm{($f^{\ast}(S)\in D(d)$ and $f^{\ast}(S)\subseteq f^{\ast}(T)$) implies $f\in T$} \}
\end{array}
\]
By using formula (15) p. 148 \cite{MM}, we get the following description of $\forall_{\pi_{2}}(A)$ for a subobject $A$ of $\Omega \times \Omega$:
\[
\forall_{\pi_{2}}(A)(c)=\{R \textrm{ sieve on $c$ | for all $f:d\to c$, } \Omega \times f^{\ast}(R)\subseteq A \}
\] 
for any object $c\in {\cal C}$. If we apply this expression to the subobject $\pi^{\ast}_{1}(D)\imp \Theta$ calculated above we thus obtain
\[
\begin{array}{ccl}
D^{r} & = & \{T \textrm{ sieve on $c$ | for all arrows $e\stackrel{h}{\to}d \stackrel{g}{\to} c$ and sieve $S$ on $d$ }\\
& & \textrm{[$h^{\ast}(S)\in D(e)$ and $h^{\ast}(S)\subseteq h^{\ast}(g^{\ast}(T))]$ implies $h\in g^{\ast}(T)$} \}
\end{array}
\]
Similarly, one can derive the following expression for $D^{l}$:
\[
\begin{array}{ccl}
D^{l} & = & \{S \textrm{ sieve on $c$ | for all arrows $e\stackrel{h}{\to}d \stackrel{g}{\to} c$ and sieve $T$ on $d$ }\\
& & \textrm{[$h^{\ast}(T)\in D(e)$ and $h^{\ast}(g^{\ast}(S))\subseteq h^{\ast}(T)]$ implies $h\in T$} \}
\end{array}
\]
Notice that the formulas above can alternatively be put in the following form:

\[
\begin{array}{ccl}
D^{r} & = & \{T\textrm{ sieve on $c$ | for any arrow $d \stackrel{f}{\to} c$ and sieve $S$ on $d$, }\\
& & \textrm{[$S\in D(d)$ and $S\subseteq f^{\ast}(T)]$ implies $f\in T$} \}
\end{array}
\]

\[
\begin{array}{ccl}
D^{l} & = & \{S \textrm{ sieve on $c$ | for any arrow $d \stackrel{f}{\to} c$ and sieve $Z$ on $d$, }\\
& & \textrm{[$Z\in D(d)$ and $f^{\ast}(S)\subseteq Z]$ implies $Z=M_{d}$} \}
\end{array}
\]
Let us for example verify the equivalence of the previous expression for $D^{l}$ with this latter formulation: take $g=f$, $h=1_{d}$ and $T=Z$ in one direction and $f=g\circ h$ and $Z=h^{\ast}(T)$ in the other direction.\\

From these expressions one immediately obtains the following formula:
\[
\begin{array}{ccl}
(D^{r})^{l} & = & \{S \textrm{ sieve on $c$ | for any arrow $d\stackrel{f}{\to} c$ and sieve $T$ on $d$, }\\
& &  [(\textrm{for any arrow $e\stackrel{g}{\to}d$ and sieve $Z$ on $e$ }\\
& & \textrm{($Z\in D(e)$ and $Z \subseteq g^{\ast}(T)$) implies $g\in T$) and ($f^{\ast}(S)\subseteq T)]$}\\
& & \textrm{implies } T=M_{d}\}
\end{array}
\]
 
We recall from the proof of Corollary A4.5.13(i) \cite{El} that the classifying map of $(D^{r})^{l}$ is the smallest local operator $j$ on $\cal E$ such that all the monomorphisms in $\cal E$ whose classifying map factors through $D\mono \Omega$ are $j$-dense. Let us now show that, via the identification $(\ast)$ local operators on ${\cal E}=[{\cal C}^{\textrm{op}}, \Set]$ with Grothendieck topologies on $\cal C$, this topology corresponds exactly to the Grothendieck topology generated by $D$, that is the smallest Grothendieck topology $J$ on $\cal C$ such that all the sieves in $D$ (regarded here as a collection of sieves in $\cal C$) are $J$-covering. To this end, it suffices to recall from \cite{El} that, given a local operator $j$ on a topos $\cal E$, the $j$-dense monomorphisms are exactly those whose classifying map factors through the subobject classified by $j$; notice that if ${\cal E}=[{\cal C}^{\textrm{op}}, \Set]$ and $j$ corresponds to a Grothendieck topology $J$ on $\cal C$, this subobject is exactly $J$ (regarded as a subobject of $\Omega_{[{\cal C}^{\textrm{op}}, \Set]}$). Now, clearly, all the sieves in $D$ are $J$-covering if and only if $D\leq J$ as subobjects of $\Omega$, so our claim immediately follows.\\
Thus, our formula for $(D^{r})^{l}$ gives an explicit description of the Grothendieck topology generated by $D$. Similarly, starting from Corollary A4.5.13(i) \cite{El}, one can prove that our formula for $D^{l}$ gives an explicit description of the largest Grothendieck topology $J$ on $\cal C$ via $(\ast)$ such that all the sieves in $D$ are $J$-closed (one replaces, in the discussion above, the subobject $J$ classifying dense monomorphisms by the subobject $\Omega_{J}$ classifying $J$-closed monomorphisms, i.e. the equalizer of the arrows $j, 1_{\Omega}:\Omega\to \Omega$).\\
As an application, let us derive an explicit formula for the Heyting operation on the collection of Grothendieck topologies on a given small category.\\
Example 4.5.14(f) \cite{El} provides a description of the Heyting operation on the collection of local operators on a topos: given local operators $j_{1}$ and $j_{2}$ on a topos $\cal E$, $j_{1}\imp j_{2}=(J_{1}\cap \Omega_{j_{2}})^{l}$. If ${\cal E}=[{\cal C}^{\textrm{op}}, \Set]$ and $j_{1}, j_{2}$ correspond to Grothendieck topologies $J_{1}, J_{2}$ on $\cal C$ via $(\ast)$ then our (second) formula for $D^{l}$ gives the following expression for $J_{1}\imp J_{2}$:

\[
\begin{array}{ccl}
J_{1}\imp J_{2}(c) & = & \{S \textrm{ sieve on $c$ | for any arrow $d \stackrel{f}{\to} c$ and sieve $Z$ on $d$ }\\
& & \textrm{[$Z$ is $J_{1}$-covering and $J_{2}$-closed and $f^{\ast}(S)\subseteq Z]$ implies $Z=M_{d}$} \}
\end{array}
\]

In particular the pseudocomplement $\neg J$ of a Grothendieck topology $J$ on $\cal C$ is given by the following formula:

\[
\begin{array}{ccl}
\neg J(c) & = & \{S \textrm{ sieve on $c$ | for any arrow $d \stackrel{f}{\to} c$ and sieve $Z$ on $d$ }\\
& & \textrm{[$Z$ is $J$-covering and $f^{\ast}(S)\subseteq Z]$ implies $Z=M_{d}$} \}
\end{array}
\]
Let us now prove directly that, given a category $\cal C$ and a collection $D$ of sieves in $\cal C$ which is closed under pullback, the above formula for $D^{l}$ always defines a Grothendieck topology on $\cal C$ and that $(D^{r})^{l}$ is the Grothendieck topology on $\cal C$ generated by $D$. This will ensure that our results hold also for a general, not necessarily small, category $\cal C$. 
In passing, note that the Grothendieck topology on $\cal C$ generated by a given family of sieves $\cal F$ in $\cal C$ can be obtained as $({{\cal F}_{\textrm{p.b.}}}^{r})^{l}$ where ${\cal F}_{\textrm{p.b.}}$ is the collection of all the sieves in $\cal C$ which are pullbacks in $\cal C$ of sieves in $\cal F$.\\
To prove that $D^{l}$ is a Grothendieck topology on $\cal C$, observe that $D^{l}$ clearly satisfies the maximality and stability axioms for Grothendieck topologies; it remains to verify that it satisfies the transitivity axiom. Let $R$ and $S$ be sieves on $c\in {\cal C}$ such that $S\in D^{l}(c)$ and for each $s:a\to c$ in $S$, $s^{\ast}(R)\in D^{l}(a)$; we want to prove that $R\in D^{l}(c)$, that is given any arrow $f:d\to c$ and sieve $Z$ on $d$, ($Z\in D(d)$ and $f^{\ast}(S)\subseteq Z$) implies $Z=M_{d}$. Now, for any $h\in f^{\ast}(S)$, $h^{\ast}(f^{\ast}(R))\subseteq h^{\ast}(Z)$ and hence $h\in Z$ since $(f\circ h)^{\ast}(R)\in D^{l}(dom(h))$. So $f^{\ast}(S)\subseteq Z$, which implies $Z=M_{d}$ since $S\in D^{l}(c)$.\\
Let us now show that $(D^{r})^{l}$ is the Grothendieck topology on $\cal C$ generated by $D$; since we already know that $(D^{r})^{l}$ is a Grothendieck topology, this amounts to verifying that for any Grothendieck topology $K$ on $\cal C$ which contains $D$, $(D^{r})^{l}\leq K$. Let $S$ be a sieve in $(D^{r})^{l}(c)$; then $S$ is $K$-covering if and only if $\overline{S}^{K}=M_{c}$. Now, if we take $f=1_{c}$ and $T=\overline{S}^{K}$ in the formula for $(D^{r})^{l}$, we have that for any arrow $e\stackrel{g}{\to}d$ and sieve $Z$ on $e$, [$Z\in D(e)$ and $Z \subseteq g^{\ast}(T)$] implies that $g^{\ast}(T)$ is $K$-covering and hence maximal (being $K$-closed), and $f^{\ast}(S)\subseteq T$; hence the formula gives that $T$ is maximal, as required.\\       
 
Also, we can verify directly that the formula for $J_{1}\imp J_{2}$ satisfies the property of the Heyting implication between $J_{1}$ and $J_{2}$, i.e. that for any Grothendieck topology $K$ on $\cal C$, $K\wedge J_{1}\leq J_{2}$ if and only if $K\leq J_{1}\imp J_{2}$. Indeed, $(J_{1}\imp J_{2})\wedge J_{1}\leq J_{2}$ since for every $S\in (J_{1}\imp J_{2})\wedge J_{1}(c)$, $S\subseteq \overline{S}^{J_{2}}$ and hence $\overline{S}^{J_{2}}$ is maximal i.e. $S$ is $J_{2}$-covering; in the other direction, if $K\wedge J_{1}\leq J_{2}$ then for any $K$-covering sieve $S$, $[Z$ is $J_{1}$-covering and $J_{2}$-closed and $f^{\ast}(S)\subseteq Z]$ implies that $Z$ is $K\wedge J_{1}$-covering and hence $J_{2}$-covering and $J_{2}$-closed i.e. maximal.\\

\subsection{The lattice operations on theories}        
By using the duality theorem, we can interpret the meaning of the lattice operations on the collection of Grothendieck topologies on the geometric syntactic category ${\cal C}_{\mathbb T}$ of a geometric theory $\mathbb T$ at the level of quotients of $\mathbb T$.\\
Let us denote by $\mathfrak{Th}_{\Sigma}^{\mathbb T}$ the collection of closed geometric theories over $\Sigma$ which are quotients of $\mathbb T$. By definition of the duality of Theorem \ref{dualita}, it is clear that the order on $\mathfrak{Th}_{\Sigma}^{\mathbb T}$ corresponding to the order $\leq$ between Grothendieck topologies on ${\cal C}_{\mathbb T}$ is the following: ${\mathbb T}'\leq {\mathbb T}''$ if and only if all the axioms of ${\mathbb T}'$ (equivalently, all the geometric sequents provable in ${\mathbb T}'$) are provable in ${\mathbb T}''$. So Theorem \ref{dualita} gives the following result
\begin{theorem}
Let $\mathbb T$ be a geometric theory over a signature $\Sigma$. Then the collection $\mathfrak{Th}_{\Sigma}^{\mathbb T}$ of closed geometric theories over $\Sigma$ which are quotients of $\mathbb T$, endowed with the order defined by `${\mathbb T}\leq {\mathbb T}'$ if and only if all the axioms of $\mathbb T$ are provable in ${\mathbb T}'$' is an Heyting algebra. 
\end{theorem}\qed
Note in particular that, by taking $\mathbb T$ to be the empty (geometric) theory over $\Sigma$, we obtain that the collection $\mathfrak{Th}_{\Sigma}^{\emptyset}$ of all the closed geometric theories over $\Sigma$ is an Heyting algebra.\\  
By definition of the order in $\mathfrak{Th}_{\Sigma}^{\mathbb T}$, we get the following description of the lattice operations in $\mathfrak{Th}_{\Sigma}^{\mathbb T}$:\\
(i) the bottom element is the closure of $\mathbb T$;\\
(ii) the top element is the contradictory theory (that is the collection of all the geometric sequents over $\Sigma$);\\
(iii) the wedge ${\mathbb T}'\wedge {\mathbb T}''$ is the largest geometric theory over $\Sigma$ which is contained in both ${\mathbb T}'$ and ${\mathbb T}''$, i.e. the collection of geometric sequents $\sigma$ over $\Sigma$ such that $\sigma$ is provable in both ${\mathbb T}$ and ${\mathbb T}'$;\\
(iv) the join ${\mathbb T}'\vee {\mathbb T}''$ is the smallest closed geometric theory over $\Sigma$ which contains both ${\mathbb T}'$ and ${\mathbb T}''$, i.e the closure of the union of the axioms of ${\mathbb T}'$ and of ${\mathbb T}''$;\\ 
(v) the implication ${\mathbb T}' \imp {\mathbb T}''$ is the largest closed geometric theory $\mathbb S$ over $\Sigma$ such that ${\mathbb S}\wedge {\mathbb T}'\leq {\mathbb T}''$, i.e. such that every geometric sequent $\sigma$ which is provable in both ${\mathbb S}$ and ${\mathbb T}$ is provable in ${\mathbb T}'$; in particular, the pseudocomplement $\neg{{\mathbb T}'}$ is the largest closed geometric theory over $\Sigma$ such that every geometric sequent $\sigma$ which is provable in both ${\mathbb S}$ and ${\mathbb T}'$ is provable in $\mathbb T$.\\
We note that these operations are quite natural from the logical perspective; however it is by no means obvious from the point of view of geometric logic that there should exist an Heyting operation on the lattice of closed geometric theories over a given signature, while this fact follows as a formal consequence of our duality theorem. Another consequence of the theorem is the fact that our lattices $\mathfrak{Th}_{\Sigma}^{\mathbb T}$ are complete (i.e. they are locales); indeed, any intersection of Grothendieck topologies is a Grothendieck topology.\\
Let us discuss, from the point of view of geometric logic, the fact that our lattice $\mathfrak{Th}_{\Sigma}^{\mathbb T}$ is distributive; this is a formal consequence of the fact that it is an Heyting algebra, so it is true by the duality theorem, but is seems instructive to justify this from the point of view of geometric logic. Explicitly, this means that for any closed geometric theories ${\mathbb T}'$ and $\{{\mathbb T}_{k} \textrm{ | } k\in K\}$, ${\mathbb T}'\wedge (\mathbin{\mathop{\textrm{\huge $\vee$}}\limits_{k\in K}} {\mathbb T}_{k})=\mathbin{\mathop{\textrm{\huge $\vee$}}\limits_{k\in K}} ({\mathbb T}'\wedge {\mathbb T}_{k})$; since the inequality $\geq$ is trivially satisfied, this amounts to verifying that for any geometric sequent $\sigma$ over $\Sigma$, if $\sigma$ is in ${\mathbb T}'$ and is derivable from axioms of the ${\mathbb T}_{k}$, then $\sigma$ is derivable from axioms of the ${\mathbb T}'\wedge {\mathbb T}_{k}$. To this end, we need the following lemma. 
\begin{lemma}
Let $\Sigma$ be a signature. If a geometric sequent $\sigma\equiv \phi \vdash_{\vec{x}} \psi$ over $\Sigma$ is provable in the theory ${\mathbb S}=\{\tau \equiv \phi_{\tau} \vdash \psi_{\tau} \textrm{ | } \tau\in {\mathbb S}\}$ using geometric logic then $\sigma$ is provable in the theory ${\mathbb S}_{\sigma}=\{ \phi_{\tau} \wedge \phi  \vdash \psi_{\tau} \vee \psi  \textrm{ | } \tau\in {\mathbb S}\}$ using geometric logic.
\end{lemma}          
\begin{proofs}
Given a geometric sequent $\tau \equiv \chi \vdash \xi$ over $\Sigma$, for a string of variables $\vec{x'}$ of the same kind as $\vec{x}$ denote by ${\cal W}_{\vec{x'}}(\tau)$ the sequent $\chi \wedge \phi[\vec{x'}\slash \vec{x}] \vdash \xi \vee \psi[\vec{x'}\slash \vec{x}]$. Then one can easily check that for any instance of an inference rule of geometric logic, if we choose a string $\vec{x'}$ of variables which are not free in any of the sequents involved in it then the image via ${\cal W}_{\vec{x'}}$ of the conclusion of the rule is derivable in geometric logic from the images via ${\cal W}_{\vec{x'}}$ of the premises of the rule. And this fact clearly implies our thesis.      
\end{proofs}
The lemma easily implies our claim. Indeed, if we have a derivation of $\sigma\in {\mathbb T}$ from axioms $\tau \equiv \phi_{\tau} \vdash \psi_{\tau}$ of any of the ${\mathbb T}_{k}$ then, by the lemma, we have a derivation of $\sigma$ from the sequents $\phi \wedge \phi_{\tau} \vdash \psi \vee \psi_{\tau}$, each of which belongs to ${\mathbb T}$, since it is derivable from $\sigma$, and from ${\mathbb T}_{k}$ whenever $\sigma_{\tau}$ lies in ${\mathbb T}_{k}$, since $\phi \wedge \phi_{\tau} \vdash \psi \vee \psi_{\tau}$ is derivable from $\tau$.\\

This is an illustration of the fact that it can be very useful to use the duality theorem to get insights into geometric logic; we will discuss other applications of this kind below.

\subsection{The Heyting implication in $\mathfrak{Th}_{\Sigma}^{\mathbb T}$}  
The purpose of this section is to give an explicit logical description of the Heyting operation between closed quotients of a given geometric theory $\mathbb T$. We will achieve this by interpreting the formula for the Heyting implication of Grothendieck topologies obtained above at the level of theories via the duality theorem.\\
The following fact about local operators will be useful for our purposes.
\begin{lemma}
Let $\cal E$ be an elementary topos and $j,j'$ two local operators on $\cal E$ with associated universal closure operators $c_{j}$ and $c_{j'}$. Then $j\leq j'$ if and only if for every subobject $m:A'\mono A$ in $\cal E$, $c_{j}(m)\leq c_{j'}(m)$; specifically, if $j\leq j'$ then for any subobject $m$ in $\cal E$, $c_{j'}(m)=c_{j'}(c_{j}(m))$.
\end{lemma}  

\begin{proofs}
Let $L_{j}$ and $L_{j'}$ the cartesian reflectors on $\cal E$ associated respectively to the local operators $j$ and $j'$. Recall that $c_{j}(m)$ is given by the pullback 
\[  
\xymatrix {
c_{j}(A') \ar[r] \ar[d]^{c_{j}(m)} & L_{j}A' \ar[d]^{L_{j}m}\\
A \ar[r]^{\eta_{a}} & L_{j}A }
\] 
If $j\leq j'$ then $L_{j'}(L_{j}(m))\cong L_{j'}(m)$ since $L_{j'}$ factors through $L_{j}$ and they are both cartesian reflectors, so if we apply the pullback-preserving functor $L_{j'}$ to the pullback above we get $L_{j'}(c_{j}(m))\cong L_{j'}(L_{j}(m))\cong L_{j'}(m)$; from this it immediately follows by definition of $c_{j'}$ in terms of $L_{j'}$ that $c_{j'}(m)=c_{j'}(c_{j}(m))$. In particular, $c_{j}(m)\leq c_{j'}(m)$.\\
The converse is clear, since $j$ is the classifying map of $c_{j}(\top)$ for each local operator $j$.   
\end{proofs} 

\begin{rmk}\label{passaggiochiusi}
\emph{We observe that it follows immediately from the lemma that if $j\leq j'$ then for any subobject $m$, if $m$ is $c_{j'}$-closed then $m$ is $c_{j}$-closed.}   
\end{rmk}

We shall also need the following results.

\begin{proposition}\label{reggeom}
Let $\cal C$ be a regular category, $J$ a Grothendieck topology on $\cal C$ such that $J\supseteq J^{\textrm{reg}}_{\cal C}$ and $r:d\to c$ be a cover in $\cal C$. Then\\
(i) for any sieve $R$ on $c$, $R\in J(c)$ if and only if $r^{\ast}(R)\in J(d)$;\\
(ii) for any sieve $R$ on $c$ generated by a monomorphism, $R$ is $J$-closed if and only if $r^{\ast}(R)$ is $J$-closed;\\
(iii) for any sieve $R$ on $c$, $R$ is $J$-closed if and only if for any monomorphism $f:d\to c$, $f^{\ast}(R)\in J(d)$ implies $f\in R$;\\
(iv) for any sieves $R$ and $T$ on $c$ such that $T$ is generated by a monomorphism, $r^{\ast}(R)\subseteq r^{\ast}(T)$ if and only if $R\subseteq T$.   
\end{proposition}

\begin{proofs}
(i) This immediately follows from the stability and transitivity axioms for Grothendieck topologies.\\
(ii) The `only if' part is obvious; let us prove the `if' part. Given an arrow $f:d\to c$ such that $f^{\ast}(R)\in J(d)$ we want to prove that $f\in R$. Consider the pullback in $\cal C$
\[  
\xymatrix {
a \ar[r]^{h} \ar[d]^{g}  & d \ar[d]^{r} \\
b \ar[r]^{f} & c }
\] 
By the commutativity of this square and the stability axiom for Grothendieck topologies, it follows that $h^{\ast}(r^{\ast}(R))\in J(a)$ and hence $h\in r^{\ast}(R)$ i.e. $r\circ h \in R$. But $f\circ g=r\circ h\in R$ so $f\circ g \in R$. But $g$ is a cover and $R$ is generated by a monomorphism so, since covers are orthogonal to monomorphisms (cfr. Lemma A1.3.2 \cite{El}), we conclude that $f\in R$, as required.\\
(iii) The `only if' part is obvious, so it remains to prove that if for any monomorphism $f:d\mono c$, $f^{\ast}(R)\in J(d)$ implies $f\in R$, then $R$ is $J$-closed. Let $g:e\to c$ be an arrow such that $g^{\ast}(R)\in J(e)$; we want to prove that $g\in R$. Denoted by $e\stackrel{g''}{\epi} u \stackrel{g'}{\mono} c$ the cover-mono factorization of $g$, we have by part (i) of the proposition that $g'^{\ast}(R)\in J(u)$; so $g'\in R$ by our hypothesis and hence $g\in R$, as required.\\
(iv) The `if' part is obvious, so it remains to prove that if $r^{\ast}(R)\subseteq r^{\ast}(T)$ then $R\subseteq T$. Given $f\in R$, consider the pullback in $\cal C$
\[  
\xymatrix {
a \ar[r]^{h} \ar[d]^{g}  & d \ar[d]^{r} \\
b \ar[r]^{f} & c }
\] 
Now, $h$ belongs to $r^{\ast}(R)$ and hence to $r^{\ast}(T)$, so $f\circ g=r\circ h\in T$. But $g$ is a cover and $T$ is generated by a monomorphism so, since covers are orthogonal to monomorphisms (cfr. Lemma A1.3.2 \cite{El}), we conclude that $f\in T$, as required.         
\end{proofs}

\begin{proposition}\label{coveringchiusi}
Let $\mathbb T$ be a geometric theory over a signature $\Sigma$, ${\mathbb T}'$ a quotient of $\mathbb T$ and $\{\{\vec{x_{i}}.\phi_{i}\} \stackrel{[\theta_{i}]}{\to} \{\vec{y}.\psi\} \textrm{ | } i\in I\}$ a set of generators for a sieve $S$ in the syntactic category ${\cal C}_{\mathbb T}$ of $\mathbb T$. Then\\
(i) $S$ is $J_{{\mathbb T}'}^{\mathbb T}$-covering if and only if $\psi \vdash_{\vec{y}} \mathbin{\mathop{\textrm{\huge $\vee$}}\limits_{i\in I}}(\exists \vec{x_{i}})\theta_{i}$ is provable in ${\mathbb T}'$;\\
(ii) $S$ is $J_{{\mathbb T}'}^{\mathbb T}$-closed if and only if it is generated by a single monomorphism and for any geometric formula $\psi'(\vec{y})$ such that $\psi' \vdash_{\vec{y}} \psi$ is provable in $\mathbb T$, the sequent $\psi' \vdash_{\vec{y}} \mathbin{\mathop{\textrm{\huge $\vee$}}\limits_{i\in I}}(\exists \vec{x_{i}})\theta_{i}$ is provable in ${\mathbb T}'$ (if and) only if it is provable in $\mathbb T$.  
\end{proposition}

\begin{proofs}
(i) This is precisely equivalence (1) after the proof of Theorem \ref{prooftheory}.\\
(ii) This follows at once from Remark \ref{passaggiochiusi}, Proposition \ref{regulargeometric}(ii), Proposition \ref{reggeom}(iii) and part (i) of this proposition, by recalling the well-known identification of subobjects of $\{\vec{y}. \psi\}$ in ${\cal C}_{\mathbb T}$ with $\mathbb T$-provable equivalence classes of geometric formulae $\psi'(\vec{y})$ over $\Sigma$ such that $\psi' \vdash_{\vec{y}} \psi$ is provable in $\mathbb T$. 
\end{proofs}

Having in mind Remark \ref{generazione}, let us look for a simple as possible set of generators of $J_{1}\imp J_{2}$.\\
We note that the collection $K$ given by
\[
\begin{array}{ccl}
K(c) & = & \{S \textrm{ sieve on $c$ | for any arrow $d\stackrel{f}{\to} c$ and sieve $T$ on $c$ }\\
& & \textrm{[$f^{\ast}(T)$ is $J_{1}$-covering and $J_{2}$-closed and $f^{\ast}(S)\subseteq f^{\ast}(T)]$}\\
& & \textrm{implies $f\in T$} \}
\end{array}
\]  
for each $c\in {\cal C}$, generates the Grothendieck topology $J_{1}\imp J_{2}$. Indeed, all the sieves in $K$ are clearly $(J_{1}\imp J_{2})$-covering and if $S\in J_{1}\imp J_{2}(c)$ then $g^{\ast}(S)\in K(d)$ for any arrow $g:d\to c$ so that our claim follows from the maximality and transitivity axioms for Grothendieck topologies.\\
Now, let us suppose that ${\cal C}$ is the syntactic category ${\cal C}_{\mathbb T}$ of a geometric theory $\mathbb T$ and that $J_{1}$ and $J_{2}$ are respectively the associated topologies $J_{{\mathbb T}_{1}}^{\mathbb T}$ and $J_{{\mathbb T}_{2}}^{\mathbb T}$ of two quotients ${\mathbb T}_{1}$ and ${\mathbb T}_{2}$ of $\mathbb T$. By Proposition \ref{gengrot}(ii), $K$ is generated over $J_{\mathbb T}$ by sieves generated by a single monic arrow. This remark enables us to arrive at a simplified axiomatization of the Heyting implication ${\mathbb T}_{1} \imp {\mathbb T}_{2}$, as follows.\\
Before applying the formula obtained above in our case, it is convenient to make a series of simplifications.\\
First, we observe that      
\[
\begin{array}{ccl}
K(c) & = & \{S \textrm{ sieve on $c$ | for any arrow $d\stackrel{f}{\to} c$ and sieve $T=(t)$ on $c$}\\
& & \textrm{with $t$ monic,} \\
& & \textrm{$[f^{\ast}(T)$ is $J_{1}$-covering and $J_{2}$-closed and $f^{\ast}(S)\subseteq f^{\ast}(T)]$}\\
& & \textrm{implies $f\in T$}\}.
\end{array}
\]   
Indeed, by Proposition \ref{regulargeometric}, $\overline{T}^{J_{\mathbb T}}_{{\cal C}_{\mathbb T}}$ is generated by a monic arrow, and if $f^{\ast}(T)$ is ($J_{1}$-covering and) $J_{2}$-closed then $f^{\ast}(\overline{T}^{J_{\mathbb T}})=\overline{f^{\ast}(T)}^{J_{\mathbb T}}=f^{\ast}(T)$, where the second equality follows from the fact that, since $J_{\mathbb T}\subseteq J_{2}$, $f^{\ast}(T)$ is $J_{\mathbb T}$-closed by Remark \ref{passaggiochiusi}.\\
Second, we note that the quantification over all the arrows $f$ in the preceding expression can be restricted to all the arrows $f$ which are monic, that is we have
\[
\begin{array}{ccl}
K(c) & = & \{S \textrm{ sieve on $c$ | for any monic arrow $d\stackrel{f}{\to} c$ and sieve $T=(t)$ on $c$}\\
& & \textrm{with $t$ monic,}\\
& & \textrm{[$f^{\ast}(T)$ is $J_{1}$-covering and $J_{2}$-closed and $f^{\ast}(S)\subseteq f^{\ast}(T)]$}\\
& & \textrm{implies $f\in T$} \}
\end{array}
\]   
Indeed, this immediately follows from Proposition \ref{reggeom} by considering the cover-mono factorization of the arrow $f$.\\
Now, we can make a futher rewriting of our formula: since, given a monic arrow $f:d\to c$ and a sieve $R$ on $d$, $R=f^{\ast}(R')$ where $R'$ is the sieve $\{f\circ g \textrm{ | } g\in R\}$, we obtain the following equality:
\[
\begin{array}{ccl}
K(c) & = & \{S \textrm{ sieve on $c$ | for any monic arrow $d\stackrel{f}{\to} c$ and sieve $T=(t)$ on $d$}\\
& & \textrm{with $t$ monic,} \\
& & \textrm{[$T$ is $J_{1}$-covering and $J_{2}$-closed and $f^{\ast}(S)\subseteq T]$}\\
& & \textrm{implies $1_{d}\in T$} \}
\end{array}
\]   
We are now ready to apply this formula to the syntactic category of our geometric theory $\mathbb T$. In view of Propositions \ref{reggeom}(iii) and \ref{coveringchiusi}, we get the following result.

\begin{theorem}
Let $\mathbb T$ be a geometric theory over a signature $\Sigma$ and ${\mathbb T}_{1}, {\mathbb T}_{2}$ two quotients of $\mathbb T$. Then ${\mathbb T}_{1}\imp {\mathbb T}_{2}$ is the theory obtained from $\mathbb T$ by adding all the axioms $\psi \vdash_{\vec{y}} \psi'$ with the property that $\psi' \vdash_{\vec{y}} \psi$ is provable in $\mathbb T$ and for any geometric formulae $\chi, \phi$ over $\Sigma$ in the context $\vec{y}$ such that $\chi \vdash_{y} \psi$ and $\phi \vdash_{\vec{y}} \chi$ are provable in $\mathbb T$, the conjunction of the facts\\
(i) $\chi \vdash_{\vec{y}} \phi$ provable in ${\mathbb T}_{1}$,\\
(ii) for any geometric formula $\xi(\vec{y})$ such that $\xi \vdash_{\vec{y}} \chi$ is provable in $\mathbb T$, the sequent $\xi \vdash_{\vec{y}} \phi$ is provable in ${\mathbb T}_{2}$ (if and) only if it is provable in $\mathbb T$,\\
(iii) $\psi' \wedge \chi \vdash_{\vec{y}} \phi$ provable in $\mathbb T$\\ 
implies that $\chi \vdash_{\vec{y}} \phi$ is provable in $\mathbb T$.
\end{theorem}\qed
In particular, we obtain that the pseudocomplement of a quotient ${\mathbb T}'$ in $\mathfrak{Th}_{\Sigma}^{\mathbb T}$ is the theory $\neg {\mathbb T}'$ obtained from $\mathbb T$ by adding all the axioms $\psi \vdash_{\vec{y}} \psi'$ with the property that $\psi \vdash_{\vec{y}} \psi'$ is provable in $\mathbb T$ and for any geometric formulae $\chi, \phi$ over $\Sigma$ in the context $\vec{y}$ such that $\chi \vdash_{y} \psi$ and $\phi \vdash_{\vec{y}} \psi'$ are provable in $\mathbb T$, the conjunction of the facts\\
(i) $\psi' \vdash_{\vec{y}} \phi$ provable in ${\mathbb T}'$,\\
(ii) $\psi' \wedge \chi \vdash_{\vec{y}} \phi$ provable in $\mathbb T$\\ 
implies that $\psi' \vdash_{\vec{y}} \phi$ is provable in $\mathbb T$.  

\newpage
 
\section{Relativization of local operators}\label{relat}
In this section we study the problem of relativizing a local operator with respect to another one, with applications to the calculations of open and quasi-closed local operators on a topos.\\
Let us recall from \cite{El} that for any topos $\cal E$ there is a bijection between universal closure operators on $\cal E$ and local operators on $\cal E$. This bijection sends a local operator $j:{\cal E}\to {\cal E}$ to the universal closure operator $c_{j}$ (also denoted $c_{L}$ where $L$ is the corresponding reflector on $\cal E$) defined, for each monomorphism $m:A'\mono A$ in $\cal E$, by the pullback square
\[  
\xymatrix {
c_{L}(A') \ar[r] \ar[d] & LA' \ar[d]^{Lm}\\
A \ar[r]^{\eta_{A}^{L}} & LA }
\] 
where $L$ is the cartesian reflector on $\cal E$ corresponding to $j$ and $\eta_{A}^{L}$ is the unit of the reflection, and a closure operator $c$ on $\cal E$ to the local operator $j_{c}:\Omega \to \Omega$ given by classifying map of the subobject $c(1\stackrel{\top}{\mono} \top)$. Let us also recall that given a local operator $j$ on $\cal E$, the domain $\Omega_{j}$ of the equalizer $e_{j}:\Omega_{j}\mono \Omega$ of the arrows $1_{\Omega}, j:\Omega \to \Omega$ is the subobject classifier of the topos $\sh_{j}({\cal E})$ and the classifying map $\chi_{m}:A\to \Omega$ of a monomorphism $m$ in $\cal E$ factors through $e_{j}$ if and only if $m$ is $c_{j}$-closed.\\
Given geometric inclusions $\palrl{{\cal F}'}{i'}{L'}{{\cal F}}$ and $\palrl{\cal F}{i}{L}{\cal E}$, let us denote by $j_{L'}$ and $j_{L}$ the corresponding local operators respectively on $\cal F$ and $\cal E$. Denoted by $\Omega$ the subobject classifier of $\cal E$, let us define $e_{L}:\Omega_{L}\mono \Omega$ to be the equalizer of $1_{\Omega_{L}}, j_{L}: \Omega \to \Omega$, $e_{L'}:{(\Omega_{L})}_{L'}\mono \Omega_{L}$ to be the equalizer of $1_{\Omega}, j_{L'}: \Omega_{L'}\to \Omega_{L'}$ and $e_{L'\circ L}:\Omega_{L' \circ L}\mono \Omega$ to be the equalizer of  $1_{\Omega}, j_{L'\circ L}: \Omega \to \Omega$.\\

\begin{lemma}\label{lemmarel}
With the above notation, the composite
\[  
\xymatrix {
{(\Omega_{L})}_{L'} \ar[r]^{e_{L'}} & \Omega_{L} \ar[r]^{e_{L}} & \Omega}
\]  
and the arrow
\[  
\xymatrix {
\Omega_{L'\circ L} \ar[rr]^{e_{L'\circ L}} & & \Omega}
\] 
are isomorphic (as objects of ${\cal E}\slash \Omega$).\\
\end{lemma}

\begin{proofs}
Let us prove that, given a subobject $m:A'\mono A$ in $\cal E$ with classifying map $\chi_{m}:A\to \Omega$,\\
(1) $\chi_{m}$ factors through $e_{L}\circ e_{L'}$ if and only if $c_{L}(m)=m$ and $c_{L'}(Lm)=Lm$;\\
(2) $\chi_{m}$ factors through $e_{L'\circ L}$ if and only if $c_{L'\circ L}(m)=m$;\\
(3) $c_{L}(m)=m$ and $c_{L'}(Lm)=Lm$ if and only if $c_{L'\circ L}(m)=m$.\\
(1) $\chi_{m}$ factors through $e_{L}\circ e_{L'}$ if and only if $\chi_{m}$ factors through $e_{L}$ and the factorization $\chi^{L}_{m}$ of $\chi_{m}$ through $e_{L}$ factors through $e_{L'}$; by definition of $e_{L}$, the first condition precisely means that $c_{L}(m)=m$, while the second, in view of the adjunction $Hom_{\sh_{j_{L}}({\cal E})}(LA, \Omega_{L})\cong Hom_{\cal E}(A, \Omega_{L})$, is equivalent to requiring that the subobject in $\sh_{j_{L}}(\cal E)$ classified by the factorization $\overline{\chi^{L}_{m}}:LA\to {\Omega}_{L}$ of $\chi^{L}_{m}$ through $\eta_{A}:A\to LA$ is $c_{L'}$-closed (by definition of $e_{L'}$). Now, consider the diagram
\[  
\xymatrix {
A' \ar[r]^{!} \ar[d]^{m} & 1 \ar[d]^{\top_{L}} \ar[r]^{!} & 1 \ar[d]^{\top} \\
A \ar[r]^{\chi^{L}_{m}} & \Omega_{L} \ar[r]^{e_{L}} & \Omega
}
\]               
where $\top_{L}$ is the factorization of $\top:1\to \Omega$ through $e_{L}$. The outer rectangle is the pullback witnessing that $\chi_{m}$ classifies $m$, while the right square is trivially a pullback (it being commutative and $e_{L}$ being monic); so we conclude from the pullback lemma that the left-hand square is a pullback. But $L$ preserves pullbacks so we obtain that the square
\[  
\xymatrix {
LA' \ar[r]^{!} \ar[d]^{Lm} & 1 \ar[d]^{\top_{L}} \\
LA \ar[r]^{\overline{\chi^{L}_{m}}} & \Omega_{L}
}
\]   
is a pullback, i.e. $\overline{\chi^{L}_{m}}$ classifies the subobject $Lm$ in $\sh_{j_{L}}(\cal E)$. This concludes the proof of (1).\\
(2) This is immediate by definition of $\Omega_{L'\circ L}$.\\
(3) By definition of $c_{L'}$ and $c_{L'\circ L}$, we have a rectangle
\[  
\xymatrix {
c_{L'\circ L}(A) \ar[d]^{c_{L'\circ L}(m)} \ar[r] & c_{L'}(LA') \ar[d]^{c_{L'}(Lm)} \ar[r] & L'(LA') \ar[d]^{L'(Lm)} \\
A \ar[r]^{\eta^{L}_{A}} & LA \ar[r]^{\eta^{L'}_{LA}} & L'(LA)
}
\] 
in which both squares are pullbacks; indeed, this follows as a consequence of the pullback lemma, since $\eta_{A}^{L'\circ L}=\eta^{L'}_{LA}\circ \eta_{A}^{L}$. In particular, notice that if $A$ is a $L$-sheaf then $c_{L'\circ L}(m)=c_{L'}(Lm)$.\\
Suppose $c_{L}(m)=m$ and $c_{L'}(Lm)=Lm$. The fact that $c_{L'}(Lm)=Lm$ implies, by definition of $c_{L}(m)$ and the fact that the left-hand square above is a pullback, that $c_{L'\circ L}(m)=c_{L}(m)$; hence, $c_{L}(m)=m$ implies $c_{L'\circ L}(m)=m$, as required. Conversely, suppose that $c_{L'\circ L}(m)=m$. Then, by applying the pullback-preserving functor $L$ to the left-hand square above, we obtain $Lm=c_{L'}(Lm)$; but then, by definition of $c_{L}(m)$, we have $c_{L'\circ L}(m)=c_{L}(m)$ and hence $c_{L}(m)=m$.\\
Now, from (1), (2) and (3) we deduce that for any subobject $m$ in $\cal E$, $\chi_{m}$ factors through $e_{L}\circ e_{L'}$ if and only if it factors through $e_{L'\circ L}$, so that the thesis of the lemma follows from the Yoneda Lemma.    
\end{proofs}   

The following definition will be central for the results in this section. 
\begin{definition}
Given a topos $\cal E$, local operators $j$ and $k$ on $\cal E$ and a local operator $k':\Omega_{j}\to \Omega_{j}$ in $\sh_{j}({\cal E})$, we say that $k$ relativizes to $k'$ at $j$ (or that $k'$ is the \emph{relativization of $k$ at $j$}) if the square
\[  
\xymatrix {
\Omega_{j} \ar[r]^{k'} \ar[d]^{e_{j}} & \Omega_{j} \ar[d]^{e_{j}} \\
\Omega \ar[r]^{k} & \Omega
}
\] 
in $\cal E$ commutes. 
\end{definition}

Notice that in the definition above, since $e_{j}$ is monic, there can be at most one relativization of $k$ at $j$.\\
The fundamental property of relativizations is given by the following result.
\begin{theorem}\label{relativization}
Let $k'$ be the relativization of $k$ at $j$ as above. Then\\
(i) $\sh_{k'}(\sh_{j}({\cal E}))=\sh_{k \vee j}({\cal E})$ (where $k \vee j$ is the join of $k$ and $j$ in the lattice of local operators on $\cal E$).\\
(ii) for any subobject $m$ in $\sh_{j}({\cal E})$, $c_{k'}(m)=c_{k}(m)$.\\
(iii) if $k\geq j$ then for any subobject $m$ in $\cal E$, $c_{k'}(L_{j}m)=c_{k}(m)$.
\end{theorem}
\begin{proofs}
(i) Let $s$ be the local operator on $\cal E$ corresponding to $\sh_{k'}(\sh_{j}({\cal E}))$, regarded as a subtopos of $\cal E$ via the composite geometric inclusion $\sh_{k'}(\sh_{j}({\cal E})) \hookrightarrow \sh_{j}({\cal E}) \hookrightarrow {\cal E}$. We have to prove that 
$e_{s}: \Omega_{s} \mono \Omega$ is isomorphic to $e_{k\vee j}:\Omega_{k\vee j}\mono \Omega$. By the Yoneda Lemma, it is equivalent to prove that for any subobject $m:A'\mono A$, $\chi_{m}$ factors through $e_{s}$ if and only if it factors through $e_{k\vee j}$. Now, by Lemma \ref{lemmarel}, $\chi_{m}$ factors through $e_{s}$ if and only if $m$ is $c_{j}$-closed and $Lm$ is $c_{k'}$-closed, where $L$ is the cartesian reflector corresponding to $j$, while, by Example A4.5.13 \cite{El}, $\chi_{m}$ factors through $e_{k\vee j}$ if and only if $m$ is both $c_{j}$-closed and $c_{k}$-closed. So we have to prove that, given a $c_{j}$-closed subobject $m:A'\mono A$, $Lm$ is $c_{k'}$-closed if and only if $m$ is $c_{k}$-closed. Consider the commutative diagram
\[  
\xymatrix {
A \ar[d]_{\eta_{A}^{L}} \ar[dr]^{\chi^{L}_{m}} \ar[r]^{\chi_{m}} & \Omega  \ar[r]^{k} & \Omega \\
LA \ar[r]_{\overline{\chi^{L}_{m}}} & \Omega_{j}  \ar[r]_{k'} \ar[u]_{e_{j}} & \Omega_{j} \ar[u]_{e_{j}}
}
\] 
where the notation is that of Lemma \ref{lemmarel}.\\
From the proof of Lemma \ref{lemmarel} we know that $\overline{\chi^{L}_{m}}$ is the characteristic map in $\sh_{j}({\cal E})$ of the subobject $Lm$. By definition of $\Omega_{k}$, $\chi_{m}$ factors through $e_{k}$ (i.e. $m$ is $c_{k}$-closed) if and only if $k\circ \chi_{m}=\chi_{m}$, while, by definition of $\Omega^{\sh_{j}({\cal E})}_{k'}$, $\overline{\chi^{j}_{m}}$ factors through $\Omega^{\sh_{j}({\cal E})}_{k'}\mono \Omega_{j}$ (i.e. $Lm$ is $c_{k'}$-closed) if and only if $k'\circ \overline{\chi^{j}_{m}}=\overline{\chi^{j}_{m}}$. Now, since $e_{j}$ is monic and $\eta_{A}$ is the unit of the reflection corresponding to $j$, $k'\circ \overline{\chi^{j}_{m}}=\overline{\chi^{j}_{m}}$ if and only if $e_{j} \circ k'\circ \overline{\chi^{j}_{m}} \circ \eta_{A}=e_{j} \circ \overline{\chi^{j}_{m}} \circ \eta_{A}$. But, by the commutativity of the diagram above, this is precisely equivalent to $k\circ \chi_{m}=\chi_{m}$.\\
(ii) The condition $k\circ e_{j}=e_{j}\circ k':\Omega_{j}\to \Omega$ is equivalent to the assertion that the subobjects classified by the maps $k\circ e_{j}$ and $e_{j}\circ k'$ are equal. Now, since $k$ classifies $c_{k}(\top)$ then $k\circ e_{j}$ classifies $e_{j}^{\ast}(c_{k}(\top))=c_{k}(e_{j}^{\ast}(\top))=c_{k}(\top_{j})$, where $\top_{j}$ is the factorization of $\top$ through $e_{j}$, while $e_{j}\circ k'$ is easily seen to classify $c_{k'}(\top_{j})$; so the condition amounts to requiring that $c_{k'}(\top_{j})=c_{k}(\top_{j})$. But every subobject in $\sh_{j}({\cal E})$ is a pullback (both in $\sh_{j}({\cal E})$ and in $\cal E$) of $\top_{j}$; thus for any subobject $m$ in $\sh_{j}({\cal E})$, $c_{k'}(m)=c_{k}(m)$, as required.\\
(iii) By (ii), it suffices to prove that if $m$ is a subobject in $\cal E$ then $c_{k}(L_{j}m)=c_{k}(m)$; this immediately follows from the definition of $c_{k}(-)$ as the pullback of $L_{k}(-)$ along the unit of the adjunction $i\vdash L_{k}$ and the fact that if $k\geq j$ then $L_{k}(m)\cong L_{k}(L_{j}(m))$.            
\end{proofs}

Now, let us consider some instances of relativizations.\\

\begin{proposition}\label{propfactor}
With the notation of Lemma \ref{lemmarel}, $j_{L'}:\Omega_{L}\to \Omega_{L}$ is the relativization of $j_{L'\circ L}$ at $j_{L}$, that is the square
\[  
\xymatrix {
\Omega_{j_{L}} \ar[r]^{j_{L'}} \ar[d]^{e_{j_{L}}} & \Omega_{j_{L}} \ar[d]^{e_{j_{L}}} \\
\Omega \ar[r]^{j_{L'\circ L}} & \Omega
}
\] 
commutes.
\end{proposition}

\begin{proofs}
We prove that the composites $e_{j_{L}}\circ j_{L'}$ and $j_{L'\circ L}\circ e_{j_{L}}$ classify the same subobject of $\sh_{j_{L}}({\cal E})$, namely $c_{L'}(\top_{L})$.\\
Consider the diagram
\[  
\xymatrix {
c_{L'}(1) \ar[r]^{!} \ar[d]^{c_{L'}(\top_{L})} & 1 \ar[d]^{\top_{L}} \ar[r]^{!} & 1 \ar[d]^{\top} \\
\Omega_{L} \ar[r]^{j_{L'}} & \Omega_{L} \ar[r]^{e_{j_{L}}} & \Omega
}
\]        
Since both squares in it are pullbacks we conlcude by the pullback lemma that $e_{j_{L}}\circ j_{L'}$ classifies $c_{L'}(\top_{L})$. 
On the other hand, if $j_{L'\circ L}$ classifies $c_{L'\circ L}(\top)$, then $j_{L'\circ L}\circ e_{j_{L}}$ classifies $e_{j_{L}}^{\ast}(c_{L'\circ L}(\top))=c_{L'\circ L}(e_{j_{L}}^{\ast}(\top))=c_{L'\circ L}(\top_{L})$; but $c_{L'\circ L}(\top_{L})=c_{L'}(\top_{L})$ (cfr. the proof of Lemma \ref{lemmarel}), so we are done. 
\end{proofs}

\begin{rmk}\label{unique}
\emph{We note that all the relativizations arising as in Proposition \ref{propfactor} have the property that $k\geq j$. We shall see below instances of relativization in which this condition does not hold. For the moment, let us note that if $k'$ is the relativization at $j$ of two local operators $k_{1}$ and $k_{2}$ then $k_{1}\vee j=k_{2}\vee j$. Indeed, this follows from Theorem \ref{relativization} by recalling the identification between subcategories of sheaves on a topos and local operators on it.}
\end{rmk}

\begin{rmk}\label{remk2}
\emph{Notice that, given $k$ and $j$ local operators on a topos $\cal E$, there exists a relativization of $k$ at $j$ if and only if $j\circ k \circ e_{j}=k\circ e_{j}$ (equivalently, $c_{k}(\top_{j})$ being classified by $k\circ e_{j}$, $c_{k}(\top_{j})$ is $j$-closed); in particular, if $k\geq j$  then $k$ relativizes at $j$.\\ Conversely, given $k'$ local operator on $\sh_{j}({\cal E})$, there always exists a local operator $k$ on $\cal E$ such that $k$ relativizes to $k'$ at $j$. Indeed, take $k$ to be the local operator on $\cal E$ corresponding to the composite of the geometric inclusions $\sh_{k'}(\sh_{j}({\cal E}))\hookrightarrow \sh_{j}({\cal E})$ and $\sh_{j}({\cal E})\hookrightarrow {\cal E}$; then, by Proposition \ref{propfactor} and Remark \ref{unique}, $k$ relativizes to $k'$ at $j$.}     
\end{rmk}

\begin{proposition}
Under the hypotheses of Theorem \ref{relativization}, if $k$ relativizes to $k'$ at $j$ then $k\vee j$ relativizes to $k'$ at $j$.
\end{proposition}

\begin{proofs}
The condition $(k\vee j)\circ e_{j}=e_{j}\circ k'$ is equivalent to the assertion that both maps classify the same subobject, equivalently that $c_{k\vee j}(\top_{j})=c_{k}(\top_{j})$. Now, since $k\leq k\vee j$, $c_{k\vee j}(\top_{j})\geq c_{k}(\top_{j})$. To show that $c_{k\vee j}(\top_{j})\leq c_{k}(\top_{j})$ it is enough to prove, by the characterization of the closure of a subobject as the smallest closed subobject containing it, that $c_{k}(\top_{j})$ is $(k\vee j)$-closed. Now, we observed in the proof of Theorem \ref{relativization} that the $(k\vee j)$-closed subobjects are exactly those which are both $j$-closed and $k$-closed so our thesis immediately follows from Remark \ref{remk2}.\\
Alternatively, our thesis follows as a consequence of Theorem \ref{relativization}(i) and Proposition \ref{propfactor}.        
\end{proofs}

Let us now show that the notions of open and quasi-closed subtopos - unlike the notion of closed subtopos - behave naturally with respect to relativizations.

\begin{proposition}\label{openqc}
Let $\cal E$ be a topos and $j$ a local operator on $\cal E$. Given a subterminal object $U$ in $\sh_{j}({\cal E})$, the open (resp. quasi-closed) local operator $o^{\sh_{j}({\cal E})}(U)$ (resp. $qc^{\sh_{j}({\cal E})}(U)$) in $\sh_{j}({\cal E})$ associated to $U$ is the relativization at $j$ of the open (resp. quasi-closed) local operator $o^{\cal E}(U)$ (resp. $qc^{\cal E}(U)$) in $\cal E$ associated to $U$ (regarded as a subterminal in $\cal E$).
\end{proposition} 

\begin{proofs}
Recall from \cite{El} that $o^{\cal E}(U)$ given by the composite
\[  
\xymatrix {
\Omega \cong 1\times \Omega \ar[rr]^{u\times 1} & & \Omega \times \Omega \ar[rr]^{\imp} & & \Omega}
\]
where $u:1\to \Omega$ is the classifying map of the subobject $U$, while $qc^{\cal E}(U)$ is the composite 
\[  
\xymatrix {
\Omega \cong \Omega \times 1 \ar[rr]^{1\times (u,u)} & & \Omega \times \Omega \times \Omega \ar[rr]^{\imp \times 1} & & \Omega \times \Omega \ar[r]^{\imp} & \Omega}
\]

From the description of the internal Heyting operations $\wedge_{\cal E},\vee_{\cal E}, \imp_{\cal E}:\Omega\to \Omega$ on $\cal E$ given in the proof of Lemma A1.6.3 \cite{El}, it easily follows that the diagrams
\[  
\xymatrix {
\Omega_{j}\times\Omega_{j} \ar[d]^{e_{j}\times e_{j}} \ar[rr]^{\wedge_{\sh_{j}({\cal E})}} & & \Omega_{j} \ar[d]^{e_{j}} & \Omega_{j}\times\Omega_{j} \ar[d]^{e_{j}\times e_{j}} \ar[rr]^{\imp_{\sh_{j}({\cal E})}} & & \Omega_{j} \ar[d]^{e_{j}} \\
\Omega\times\Omega \ar[rr]_{\wedge_{\cal E}} & & \Omega & \Omega\times\Omega \ar[rr]_{\imp_{\cal E}} & & \Omega}
\]
are commutative.\\
Let us begin by proving that the left-hand square commutes. The arrow ${\wedge_{\cal E}}: \Omega\times\Omega \to \Omega$ is the classifying map of $(\top,\top):1\mono \Omega\times\Omega$ and ${\wedge_{\sh_{j}({\cal E})}}: \Omega_{j}\times\Omega_{j} \to \Omega$ is the classifying map of $(\top_{j},\top_{j}):1\mono \Omega_{j}\times\Omega_{j}$, that is of the factorization of $(\top,\top)$ through $e_{j}\times e_{j}$.\\ 

Let us prove that the composites $e_{j}\circ \wedge_{\sh_{j}({\cal E})}$ and $\wedge_{\cal E}\circ e_{j}$ classify the same subobject of $\sh_{j}({\cal E})$, namely $(\top_{j}, \top_{j}):1\mono \Omega_{j}\times \Omega_{j}$.\\
Consider the diagram
\[  
\xymatrix {
1 \ar[rr]^{!} \ar[d]^{(\top_{j}, \top_{j})} & & 1 \ar[d]^{\top_{j}} \ar[r]^{!} & 1 \ar[d]^{\top} \\
\Omega_{j}\times \Omega_{j} \ar[rr]^{\wedge_{\sh_{j}({\cal E})}} & & \Omega_{j} \ar[r]^{e_{j}} & \Omega
}
\]        
Since both squares in it are pullbacks we conlcude by the pullback lemma that $e_{j}\circ \wedge_{\sh_{j}({\cal E})}$ classifies $(\top_{j}, \top_{j})$. On the other hand, if $\wedge_{\cal E}$ classifies $(\top,\top)$, then $\wedge_{\cal E}\circ e_{j}$ classifies $e_{j}^{\ast}((\top,\top))=(\top_{j},\top_{j})$. This proves that the square for $\wedge$ commutes.\\

Let us now prove that the square for $\imp$ commutes. $\Omega\times\Omega \stackrel{{\imp_{\cal E}}}{\to} \Omega$ is the classifying map of $r:E\mono \Omega\times\Omega$ and $\Omega_{j}\times\Omega_{j} \stackrel{{\imp_{\sh_{j}({\cal E})}}}{\to} \Omega$ is the classifying map of $r_{j}:E_{j}\mono \Omega_{j}\times\Omega_{j}$, where $r$ and $r_{j}$ are respectively the equalizer of $\wedge_{\cal E}, \pi^{\cal E}_{1}:\Omega\times\Omega\to \Omega$ and of $\wedge_{\sh_{j}({\cal E})}, \pi^{\sh_{j}({\cal E})}_{1}:\Omega_{j}\times\Omega_{j}\to \Omega_{j}$.\\   
It is easy to verify, by using the commutativity of the square for $\wedge$, that the pullback of $r$ along $e_{j}\times e_{j}$ is an equalizer for $\wedge_{\sh_{j}({\cal E})}, \pi^{\sh_{j}({\cal E})}_{1}:\Omega_{j}\times\Omega_{j}\to \Omega_{j}$, and hence isomorphic to $r_{j}$; from this our claim immediately follows.\\
Now, by definition of open and quasi-closed local operators, the commutativity of the diagrams for $\wedge$ and $\imp$ immediately implies our thesis, since if $U$ is a subterminal in $\sh_{j}({\cal E})$ then the classifying map of $U\mono 1$ in $\sh_{j}({\cal E})$ is the factorization of its classifying map in $\cal E$ through $e_{j}:\Omega_{j}\mono \Omega$.

\end{proofs}

As an application of Theorem \ref{relativization} and Proposition \ref{openqc}, we deduce the following well-known fact.

\begin{corollary}\label{cordense}
Let $\cal E$ be topos and $j$ be a dense (i.e. $j\leq \neg\neg$) local operator on $\cal E$. Then $\sh_{\neg\neg_{\cal E}}(\sh_{j}({\cal E}))=\sh_{\neg\neg}({\cal E})$.
\end{corollary}

\begin{proofs}
For any topos $\cal E$, $qc^{\cal E}(0_{\cal E})=\neg\neg_{\cal E}$ (cfr. \cite{El}). The corollary then follows from Theorem \ref{relativization} and Proposition \ref{openqc} by invoking the fact (remarked in \cite{El}) that for a dense local operator $j$ on $\cal E$, the inclusion $\sh_{j}({\cal E})\hookrightarrow {\cal E}$ preserves the initial object. 
\end{proofs}

To conclude this section, let us remark a useful fact. Given an elementary topos $\cal E$, we denote by $\textbf{Lop}({\cal E})$ the collection of local operators on $\cal E$, endowed with the Heyting algebra structure given by the canonical order between topologies (cfr. \cite{El}). Let us note that, given a local operator $j$ on $\cal E$, there is a bijection between the collection of local operators $k$ in $\cal E$ such that $k\geq j$ and the collection of local operators on $\sh_{j}({\cal E})$. Indeed, if $k\geq j$ then the geometric inclusion $\sh_{k}({\cal E})\hookrightarrow {\cal E}$ factors (uniquely up to isomorphism) through $\sh_{j}({\cal E})\hookrightarrow {\cal E}$ and hence it correponds to a unique local operator $k_{j}$ on $\sh_{j}({\cal E})$ such that $\sh_{k_{j}}(\sh_{j}({\cal E}))=\sh_{k}({\cal E})$, while conversely, given a local operator $s$ on $\sh_{j}({\cal E})$, the geometric inclusion given by the composite $\sh_{s}(\sh_{j}({\cal E})) \hookrightarrow \sh_{j}({\cal E}) \hookrightarrow {\cal E}$ corresponds to a unique local operator $s^{j}$ on $\cal E$ such that $\sh_{s}(\sh_{j}({\cal E}))=\sh_{s^{j}}({\cal E})$. It is clear that these two correspondences are inverse to each other. Moreover, since the order between local operators on a topos corresponds exactly to the reverse inclusion between the corresponding subcategories of sheaves, we see that these bijections are also order-preserving, where the order between local operators $k\geq j$ on $\cal E$ is the (restriction of the) order in $\textbf{Lop}({\cal E})$ and the order between local operators on $\sh_{j}({\cal E})$ is the order in $\textbf{Lop}(\sh_{j}({\cal E}))$. Now, recall that given a Heyting algebra $H$ and an element $a\in H$, $\princfil{a}$ is a Heyting algebra which is closed under the operations of conjunction, disjunction and Heyting implication in $H$ and hence the map $a\vee (-):H\to \princfil{a}$ is an Heyting algebra homomorphism. So the bijections $(-)_{j}$ and $(-)^{j}$ are isomorphisms of Heyting algebras between the subalgebra $\princfil{j}$ of $\textbf{Lop}({\cal E})$ and $\textbf{Lop}(\sh_{j}({\cal E}))$ and hence the map $(j\vee (-))_{j}:\textbf{Lop}({\cal E})\to \textbf{Lop}(\sh_{j}({\cal E}))$ is a Heyting algebra homomorphism.\\

\section{Open, closed, quasi-closed subtoposes}

\subsection{Open subtoposes}\label{opensub}
Let us recall from section A4.5 \cite{El} that an open subtopos of a topos $\cal E$ is a geometric inclusion of the form ${\cal E}\slash U \hookrightarrow {\cal E}$ for a subterminal object $U$ in $\cal E$. The relevant universal closure operation sends a subobject $A'\mono A$ to the implication $(A\times U)\imp A'$ in the Heyting algebra $\Sub(A)$; so, if $L_{U}:{\cal E}\to {\cal E}$ is the corresponding cartesian reflector, then a monomorphism $A'\mono A$ is $L_{U}$-dense if and only if $(A \times U)\leq A'$ in $\Sub(A)$. Thus $A\times U$ is the smallest $L_{U}$-dense subobject of $A$, from which it follows that $L_{U}$ is the smallest local operator on $\cal E$ such that the monomorphism $U\mono 1$ is dense (cfr. the discussion preceding Lemma A4.5.10 \cite{El}). From Proposition A4.3.11 \cite{El} we then deduce that a geometric morphism $f:{\cal F}\to {\cal E}$ factors through the inclusion ${\cal E}\slash U \hookrightarrow {\cal E}$ if and only if $f^{\ast}(U)=1$.\\

Let $\cal E$ be the classifying topos $\Set[{\mathbb T}]\simeq \Sh({\cal C}_{\mathbb T}, J_{\mathbb T})$ of a geometric theory $\mathbb T$ over a signature $\Sigma$; we now describe the quotient of $\mathbb T$ corresponding via Theorem \ref{dualita} to an open subtopos ${\cal E}\slash U \hookrightarrow {\cal E}$ of $\cal E$. Recall that the geometric syntactic category ${\cal C}_{\mathbb T}$ of $\mathbb T$ embeds into its $\infty$-pretopos completion $\Sh({\cal C}_{\mathbb T}, J_{\mathbb T})$ via the Yoneda embedding $y:{\cal C}_{\mathbb T}\hookrightarrow \Sh({\cal C}_{\mathbb T}, J_{\mathbb T})$, and under this identification all the subobjects in $\Set[{\mathbb T}]$ of an object in ${\cal C}_{\mathbb T}$ lie again in ${\cal C}_{\mathbb T}$. Since the terminal object of $\cal E$ can be identified with $\{[].\top\}$ and the subobjects of a given object $\{\vec{x}.\psi\}$ of ${\cal C}_{\mathbb T}$ can be identified with the geometric formulae $\phi(\vec{x})$ which $\mathbb T$-provably imply $\psi(\vec{x})$ (Lemma D1.4.4(iv) \cite{El2}), we conclude that the subterminal object $U$ of $1$ in $\cal E$ corresponds to a unique (up to $\mathbb T$-provable equivalence) geometric sentence $\phi$ over $\Sigma$. Alternatively, $\{[].\phi\}$ arises as the domain of the subobject of $\{[].\top\}$ which is the union of all the images of the morphisms from objects in $\{c\in {\cal C}_{\mathbb T} \textrm{ | } U(c)=\ast\}$ to the terminal object of ${\cal C}_{\mathbb T}$ (cfr. Proposition \ref{gengrot}).\\
Let us recall that the Diaconescu's equivalence ${\bf Geom}({\cal F},\Sh({\cal C}_{\mathbb{T}},J_{\mathbb T})) \simeq {\bf Flat}_{J_{{\mathbb T}}}({\cal C}_{\mathbb T}, {\cal F})$ sends a geometric morphism $f:{\cal F}\to {\cal E}$ to the functor $f^{\ast}\circ y:{\cal C}_{\mathbb T}\to {\cal F}$ (where $y:{\cal C}_{\mathbb T}\to \Sh({\cal C}_{\mathbb{T}},J_{\mathbb T}$) is the Yoneda embedding) while the equivalence ${\mathbb T}\textrm{-mod}({\cal F}) \simeq {\bf Flat}_{J_{\mathbb T}}({\cal C}_{\mathbb T}, {\cal F})$ sends each model $M\in {\mathbb T}\textrm{-mod}({\cal E})$ to the functor $F_{M}:{\cal C}_{\mathbb T}\rightarrow {\cal E}$ assigning to a formula $\{\vec{x}.\phi\}$ its interpretation $[[\phi(\vec{x})]]_{M}$ in $M$. Thus, via the composite equivalence ${\bf Geom}({\cal F},\Sh({\cal C}_{\mathbb{T}},J_{\mathbb T})) \simeq {\mathbb T}\textrm{-mod}({\cal F})$, the geometric morphisms ${\cal F}\to {\cal E}$ which factor through ${\cal E}\slash U \hookrightarrow {\cal E}$ correspond to the $\mathbb T$-models $M$ such that $[[\phi]]_{M}=1$, i.e. such that $\phi$ is satisfied in $M$. Hence we deduce that the quotient ${\mathbb T}_{\phi}$ of $\mathbb T$ obtained by adding to $\mathbb T$ the axiom $\top \vdash_{[]} \phi$ is classified by the topos ${\cal E}\slash U$ and corresponds to it the via the duality of Theorem \ref{dualita}.\\

Let us now describe the effect of taking slices on the site representation of a Grothendieck topos $\cal E$ as the category of sheaves $\Sh({\cal C}, J)$ on a category $\cal C$ with respect to a Grothendieck topology $J$ on $\cal C$. The subterminal $U$ can be identified, by Remark C2.3.21 \cite{El2}, with a $J$-ideal on $\cal C$; if we regard this ideal as a full subcategory ${\cal C}'$ of $\cal C$ (that is, ${\cal C}'$ is the full subcategory of $\cal C$ on the objects $c$ such that $U(c)\cong 1_{\Set}$) then we have $\Sh({\cal C}, J)\slash U\simeq \Sh({\cal C}', J|_{{\cal C}'})$. Indeed, we may define an equivalence between $\Sh({\cal C}, J)\slash U$ and  $\Sh({\cal C}', J|_{{\cal C}'})$ as follows. Given a object $G\to U$ in $\Sh({\cal C}, J)/U$, for every $c\in {\cal C}$ not belonging to ${\cal C}'$, $G(c)=\emptyset$, since we have an arrow $G(c)\to U(c)$ and $U(c)=\emptyset$; if we associate to it the restriction $G|_{{\cal C}'}$ then we obtain a $J|_{{\cal C}'}$-sheaf by definition of induced Grothendieck topology on ${\cal C}'$. It is now clear that this assigment defines a geometric equivalence between our two toposes; moreover, it is easy to see that the inclusion ${\cal E}\slash U \hookrightarrow {\cal E}$ corresponds, via the equivalence ${\cal E}\slash U \simeq \Sh({\cal C}', J|_{{\cal C}'})$ to the geometric inclusion $\Sh({\cal C}', J|_{{\cal C}'})\to \Sh({\cal C}, J)$ induced by the morphism of sites $({\cal C}', J|_{{\cal C}'})\to ({\cal C}, J)$ given by the inclusion ${\cal C}' \hookrightarrow {\cal C}$.\\ 
Given a topos $\Sh({\cal C}, J)$, and a subterminal object $U$ in it, the topos ${\Sh({\cal C}, J)}\slash U$ is a subtopos of $\Sh({\cal C}, J)$, so it corresponds to a unique Grothendieck topology $J_{U}^{\textrm{open}}$ on $\cal C$ such that $J_{U}^{\textrm{open}}\supseteq J$; let us now describe this topology explicitly. By Theorem \ref{relativization} and Proposition \ref{openqc}, this topology is $J\vee J_{o(U)}$, where $J_{o(U)}$ is the Grothendieck topology on $\cal C$ corresponding via $(\ast)$ to the open local operator $o(U)$ on $[{\cal C}^{\textrm{op}}, \Set]$ associated to $U$. Now, $o(U)$ is by definition given by the composite
\[  
\xymatrix {
\Omega \cong 1\times \Omega \ar[rr]^{u\times 1} & & \Omega \times \Omega \ar[rr]^{\imp} & & \Omega}
\]
where $u:1\to \Omega$ is the classifying map of the subobject $U$. If $\cal E$ is the topos $[{\cal C}^{\textrm{op}}, \Set]$ then $U$ can be identified with the full subcategory ${\cal C}_{U}$ of $\cal C$ on the objects $c$ such that $U(c)=\{\ast\}$. So $u(c)(\ast)=\{f:d\to c \textrm{ | } d\in {\cal C}_{U}\}$ for any object $c\in {\cal C}$. Let us put, for any $c\in {\cal C}$, $Z(c)=u(c)(\ast)$. Then an easy calculation shows that $o(U)$ sends a sieve $R$ on an object $c\in {\cal C}$ to $\{g:e\to c \textrm{ | } g^{\ast}(Z(c))\subseteq g^{\ast}(R)\}$. Hence $J_{o(U)}$ is given by:
\[
R\in J_{o(U)}(c) \textrm{ if and only if } R\supseteq Z(c)
\]for any $c\in {\cal C}$. In particular, by property (ii) in Definition \ref{def2}, $J_{o(U)}$ is generated by the sieves $Z(c)$, as $c$ varies in $\cal C$. In passing, notice that for any arrow $f:d\to c$ in $\cal C$, $f^{\ast}(Z(c))=Z(d)$.\\            
Finally, let us apply this discussion to the syntactic representation $\Sh({\cal C}_{\mathbb T}, J_{\mathbb T})$ of the classiying topos $\Set[{\mathbb T}]$ of a geometric theory $\mathbb T$ over a signature $\Sigma$. From our discussion above it is clear that the subterminal in $\Sh({\cal C}_{\mathbb T}, J_{\mathbb T})$ corresponding to a sentence $\phi$ is the representable $y(\{[].\phi\})$, so that the subcategory ${\cal C}_{\phi}$ corresponding to it is the full subcategory of ${\cal C}_{\mathbb T}$ on the objects $\{\vec{x}.\psi \}$ of ${\cal C}_{\mathbb T}$ such that there exists (exactly) one morphism $\{\vec{x}.\psi \} \to \{[].\phi\}$ in ${\cal C}_{\mathbb T}$. Thus, by recalling the definition of morphism in the syntactic category ${\cal C}_{\mathbb T}$, one immediately obtains the following characterization for the objects of ${\cal C}_{\phi}$: $\{\vec{x}.\psi \} \in {\cal C}_{\phi}$ if and only if the sequent $\psi \vdash_{\vec{x}} \phi$ is provable in $\mathbb T$.\\
By definition of ${\cal C}_{\phi}$, the sieve $Z(\{[].\top\}))$ is generated over $J$ by the morphism $\{[].\phi\}\mono \{[].\top\}$ so, since $J_{U}^{\textrm{open}}$ is generated by the sieves $Z(c)$, and for any $c\in {\cal C}_{\mathbb T}$ $Z(c)$ is the pullback of $Z(\{[].\top\})$ along the unique arrow $c\to \{[].\top\}$, Theorem \ref{dualita} implies that the theory over $\Sigma$ classified by $\Sh({\cal C}_{\mathbb T}, J_{\mathbb T})\slash U$ is axiomatized over $\mathbb T$ by the sequent $\top \vdash_{[]} \phi$ (cfr. Remark \ref{generazione}). We have thus recovered the result obtained at the beginning of this section.
 
\subsection{Closed subtoposes}\label{subclosed}
We recall from \cite{El} that, given an elementary topos $\cal E$ and a subterminal object $U$ in $\cal E$, the closed local operator $c(U)$ associated to $U$ is the composite 
\[  
\xymatrix {
\Omega \cong 1\times \Omega \ar[rr]^{u\times 1} & & \Omega \times \Omega \ar[rr]^{\vee} & & \Omega}
\]
where $u:1\to \Omega$ is the classifying map of the subobject $U$. Unlike open and quasi-closed local operators, a closed local operator on $\cal E$ associated to a subterminal $U$ in a subtopos $\sh_{j}({\cal E})$ does not relativize to the closed local operator on $\sh_{j}({\cal E})$ associated to $U$; however, if $\cal E$ is the topos $[{\cal C}^{\textrm{op}}, \Set]$ we may easily find a local operator on $\cal E$ which relativizes to $c^{\Sh({\cal C}, J)}(U)$. Indeed, $c^{\Sh({\cal C}, J)}(U)$ is easily seen to be the map which sends a $J$-closed sieve $R$ on $c\in {\cal C}$ to the ($J$-closed) sieve $\{f:d\to c \textrm{ | } f^{\ast}(Z(c))\cup f^{\ast}(R)\in J(d)\}$; this naturally leads us to consider the arrow $\Omega_{[{\cal C}^{\textrm{op}}, \Set]} \to \Omega_{[{\cal C}^{\textrm{op}}, \Set]}$ in $[{\cal C}^{\textrm{op}}, \Set]$ sending a sieve $R$ on $c\in {\cal C}$ to the sieve $\{f:d\to c \textrm{ | } f^{\ast}(Z(c))\cup f^{\ast}(R)\in J(d)\}$. It is easily checked that this arrow is a local operator on $[{\cal C}^{\textrm{op}}, \Set]$ (since it corresponds via $(\ast)$ to a Grothendieck topology, say $J_{U}^{\textrm{closed}}$, on $\cal C$) and that it relativizes to $c^{\Sh({\cal C}, J)}(U)$. Thus $J_{U}^{\textrm{closed}}$ is given by:
\[
R\in J_{U}^{\textrm{closed}}(c) \textrm{ if and only if } Z(c)\cup R \in J(c) 
\]
for any $c\in {\cal C}$. Since $J_{U}^{\textrm{closed}}\supseteq J$ then $J_{U}^{\textrm{closed}}$ is, by Theorem \ref{relativization}, the (unique) Grothendieck topology $J_{U}^{\textrm{closed}}$ on $\cal C$ which corresponds to $c^{\Sh({\cal C}, J)}(U)$ of \Sh({\cal C}, J) (here regarded as a subtopos of $[{\cal C}^{\textrm{op}}, \Set]$ via the canonical geometric inclusion $\Sh({\cal C}, J)\hookrightarrow [{\cal C}^{\textrm{op}}, \Set]$).\\ 
Now, let us give a description of the theory ${\mathbb T}^{\textrm{closed}}_{\phi}$ over $\Sigma$ corresponding via Theorem \ref{dualita} to the closed subtopos $c^{\Sh({\cal C}_{\mathbb T}, J_{\mathbb T})}(U)$ of the classifying topos $\Set[\mathbb{T}]\simeq \Sh({\cal C}_{\mathbb T}, J_{\mathbb T})$ of $\mathbb T$ where $\phi$ is the geometric sentence over $\Sigma$ corresponding to $U$ as above. Since for any $c\in {\cal C}_{\mathbb T}$ $Z(c)$ is the pullback of $Z(\{[].\top\})$ along the unique arrow $c\to \{[].\top\}$, Theorem \ref{dualita} and Proposition \ref{generazione} give the following axiomatization for ${\mathbb T}^{\textrm{closed}}_{\phi}$: ${\mathbb T}^{\textrm{closed}}_{\phi}$ is obtained from $\mathbb T$ by adding the axiom 
\[
\psi \vdash_{\vec{y}} \psi'
\]
for any sequents $\psi' \vdash_{\vec{y}} \psi$ and $\psi \vdash_{\vec{y}} \psi' \vee (\phi \wedge \psi)$ which are provable in $\mathbb T$.     

\subsection{Quasi-closed subtoposes}
We recall from \cite{El} that, given an elementary topos $\cal E$ and a subterminal object $U$ in $\cal E$, the quasi-closed local operator $qc^{\cal E}(U)$ associated to $U$ is the composite 
\[  
\xymatrix {
\Omega \cong \Omega \times 1 \ar[rr]^{1\times (u,u)} & & \Omega \times \Omega \times \Omega \ar[rr]^{\imp \times 1} & & \Omega \times \Omega \ar[r]^{\imp} & \Omega}
\]
where $u:1\to \Omega$ is the classifying map of the subobject $U$.\\
If $\cal E$ is the topos $\Sh({\cal C}, J)$ and $U$ is a subterminal object $U$ in $\Sh({\cal C}, J)$, then $qc^{\cal E}(U)$ corresponds to a unique Grothendieck topology $J_{U}^{\textrm{qc}}$ on $\cal C$ such that $J_{U}^{\textrm{qc}}\supseteq J$; let us describe this topology explicitly. By Theorem \ref{relativization} and Proposition \ref{openqc}, this topology is $J\vee J_{qc(U)}$, where $J_{qc(U)}$ is the Grothendieck topology on $\cal C$ corresponding via $(\ast)$ to the quasi-closed local operator $qc^{[{\cal C}^{\textrm{op}}, \Set]}(U)$ on $[{\cal C}^{\textrm{op}}, \Set]$ associated to $U$ (regarded here as a subterminal in $[{\cal C}^{\textrm{op}}, \Set]$).\\
As above, let us identify $U$ with the full subcategory ${\cal C}_{U}$ of $\cal C$ on the objects $c$ such that $U(c)=\{\ast\}$ and put, for any $c\in {\cal C}$, $Z(c)=u(c)=\{f:d\to c \textrm{ | } d\in {\cal C}_{U}\}$ for any object $c\in {\cal C}$.  
In the case ${\cal E}=[{\cal C}^{\textrm{op}}, \Set]$ the local operator $qc^{\cal E}(U)$ is easily seen to send a sieve $R$ on $c\in {\cal C}$ to the sieve $\imp(c)(\{f:d\to c \textrm{ | } f^{\ast}(R)\subseteq f^{\ast}(Z(c))\}, Z(c))$, and hence $J_{qc(U)}$ is given by:  
\[
R\in J_{qc(U)}(c) \textrm{ if and only if for any $f:d\to c$, } (f^{\ast}(R)\subseteq Z(d) \textrm{ implies } f\in Z(c)) 
\]
for any $c\in {\cal C}$.\\
In order to specialize the above expression to the syntactic site of a geometric theory, let us observe that, if $\cal C$ is a geometric category and $J$ contains the geometric topology $J^{\textrm{geom}}$ on $\cal C$ then the condition in the right-hand side of the equivalence is satisfied for $f:d\to c$ if and only if it is satisfied by the image $f':d'\mono c$ of $f$ in $\cal C$. Indeed, since ${\cal C}_{U}$ is a $J$-ideal and every cover generates a $J$-covering sieve then $f\in Z(c)$ if and only if $f'\in Z(c)$. Now, let us prove that for any $f:d\to c$, $f'^{\ast}(R)\subseteq Z(d')$ if and only if $f^{\ast}(R)\subseteq Z(d)$.\\
Since the $Z(c)$ are stable under pullback, $f'^{\ast}(R)\subseteq Z(d')$ clearly implies $f^{\ast}(R)\subseteq Z(d)$. Conversely, let $r:d\epi d'$ be the factorization of $f$ through $f'$; given $g'\in f'^{\ast}(R)$, consider the pullback
\[  
\xymatrix {
e \ar[d]^{r'} \ar[r]^{g'} & d \ar[d]^{r}\\
e' \ar[r]^{g} & d'
}
\]
in $\cal C$. Clearly, since $R$ is a sieve, $g'\in f^{\ast}(R)\subseteq Z(d)$ and hence $e\in {\cal C}_{U}$; but $r'$ is a cover, which implies that $e'\in {\cal C}_{U}$ and hence that $g'\in Z(d')$, as required.\\
This remark enables us to achieve a simplified description of the theory ${\mathbb T}^{\textrm{qc}}_{\phi}$ over $\Sigma$ corresponding via Theorem \ref{dualita} to the quasi-closed subtopos $qc^{\Set[\mathbb{T}]}(U)$ of the classifying topos $\Set[\mathbb{T}]\simeq \Sh({\cal C}_{\mathbb T}, J_{\mathbb T})$ of $\mathbb T$, where $\phi$ is the geometric sentence over $\Sigma$ corresponding to a subterminal $U$ of $\Set[{\mathbb T}]$. Indeed, by recalling the identification between $\mathbb T$-provable equivalence classes of geometric formulae $\psi'(\vec{y})$ such that $\psi' \vdash_{\vec{y}} \psi$ is provable in $\mathbb T$ and subobjects of $\{\vec{y}. \psi\}$ in ${\cal C}_{\mathbb T}$ given by Lemma D1.4.4 \cite{El2}, we get, by Theorem \ref{dualita}, Proposition \ref{gengrot}(ii) and the syntactic characterization of ${\cal C}_{\phi}$ given in section \ref{opensub}, the following axiomatization for ${\mathbb T}^{\textrm{qc}}_{\phi}$: ${\mathbb T}^{\textrm{qc}}_{\phi}$ is obtained from $\mathbb T$ by adding the axioms
\[
\psi \vdash_{\vec{y}} \psi'
\]
where $\psi' \vdash_{\vec{y}} \psi$ is provable in $\mathbb T$ and for any geometric formula $\chi(\vec{y})$ over $\Sigma$ such that $\chi \vdash_{\vec{y}} \psi$ is provable in $\mathbb T$, $\chi \wedge \psi' \vdash_{y} \phi$ implies $\chi \vdash_{\vec{y}} \phi$.\\
Notice that if $\phi$ is $\bot$ then, in view of Remark \ref{generazione}, we recover the Booleanization of $\mathbb T$ defined in \cite{OC3}, that is the geometric theory over $\Sigma$ obtained from $\mathbb T$ by adding the axiom\\
\[
\top \vdash_{\vec{y}} \psi
\]
for any stably consistent formula $\psi(\vec{y})$ with respect to $\mathbb T$ (i.e. a formula-in-context $\psi(\vec{y})$ such that for any geometric formula $\chi(\vec{y})$ in the same context such that $\chi \vdash_{\vec{y}} \bot$ is not provable in $\mathbb T$, $\chi \wedge \psi \vdash_{\vec{y}} \bot$ is not provable in $\mathbb T$).

\section{The dense-closed factorization of a geometric inclusion}\label{denseclosed}

We recall from \cite{El} that the dense-closed factorization of a geometric inclusion $\sh_{j}({\cal E})\hookrightarrow {\cal E}$ in elementary topos theory is defined to be $\sh_{j}({\cal E})\hookrightarrow \sh_{c(ext(j))}({\cal E}) \hookrightarrow{\cal E}$, where $ext(j)$ is the $c_{j}$-closure of $0\mono 1$; the local operator $c(ext(j))$ is said to be the the \emph{closure of $j$} and denoted by $\overline{j}$.\\
In this section we interpret the meaning of this construction at the level of Grothendieck toposes and later, via the duality theorem, in terms of theories.\\
Let $\Sh({\cal C}, J)$ be a Grothendieck topos, $a_{J}:[{\cal C}^{\textrm{op}}, \Set]\to \Sh({\cal C}, J)$ the associated sheaf functor and $J'$ a Grothendieck topology on $\cal C$ which contains $J$.\\
Let us calculate the dense-closed factorization of the obvious geometric inclusion $\Sh({\cal C}, J') \hookrightarrow \Sh({\cal C}, J)$. Let us denote by $\tau^{J}_{J'}$ the corresponding local operator on $\Sh({\cal C}, J)$.\\     
The monomorphism $0\mono 1$ in $\Sh({\cal C}, J)$ is the image of the morphism $0\mono 1$ in $[{\cal C}^{\textrm{op}}, \Set]$ via the associated sheaf functor $a_{J}:[{\cal C}^{\textrm{op}}, \Set] \to \Sh({\cal C}, J)$; from Proposition \ref{propfactor} and Theorem \ref{relativization}(iii) we then deduce that the closure of $0\mono 1$ in $\Sh({\cal C}, J)$ with respect to the local operator corresponding to the geometric inclusion $\Sh({\cal C}, J') \hookrightarrow \Sh({\cal C}, J)$ is equal to the $J'$-closure of $0\mono 1$ in $[{\cal C}^{\textrm{op}}, \Set]$. Now, recall from \cite{MM} (formula (6) p. 235) that, for any Grothendieck topology $K$ on $\cal C$, the $K$-closure $c_{K}(A')$ of a subobject $A'\mono E$ in $[{\cal C}^{\textrm{op}}, \Set]$ is given by:\\
\[
e\in c_{K}(A')(c) \textrm{ if and only if } \{f:d\to c \textrm{ | } E(f)(e)\in A'(d)\}\in K(c)
\]
Given a subterminal $U$ in $[{\cal C}^{\textrm{op}}, \Set]$, identified with the full subcategory ${\cal C}_{U}$ of $\cal C$ as above in this paper, it is immediate to check that the $K$-closure $c_{K}(U\mono 1)$ of $U\mono 1$ in $[{\cal C}^{\textrm{op}}, \Set]$ identifies with the full subcategory ${\cal C}_{U}^{K}$ on the objects $c\in {\cal C}$ such that $\{f:d\to c \textrm{ | } d\in {\cal C}_{U}\}\in K(c)$; in particular, if $U=0$ then the objects of ${\cal C}_{0}^{K}$ are exactly the objects $c\in {\cal C}$ such that $\emptyset \in K(c)$.\\
By applying this discussion to our topology $J'$ we obtain that $ext(\tau^{J}_{J'})$ identifies (as a subterminal object in $\Sh({\cal C}, J)$) with the $J$-ideal ${\cal C}_{0}^{J'}=\{c\in {\cal C} \textrm{ | } \emptyset \in J'(c)\}$. So, by recalling the description of closed local operators on Grothendieck toposes given in section \ref{subclosed}, we obtain that the dense-closed factorization of the inclusion $\Sh({\cal C}, J') \hookrightarrow \Sh({\cal C}, J)$ is given by $\Sh({\cal C}, J') \hookrightarrow \Sh({\cal C}, J_{{\cal C}_{0}^{J'}}^{\textrm{closed}}) \hookrightarrow \Sh({\cal C}, J)$ where the topology $J_{{\cal C}_{0}^{J'}}^{\textrm{closed}}$ is defined by:  
\[
R\in J_{{\cal C}_{0}^{J'}}^{\textrm{closed}}(c) \textrm{ if and only if } Z(c)\cup R \in J(c) 
\]
where, for any $c\in {\cal C}$, $Z(c)=\{f:d\to c \textrm{ | } \emptyset \in J'(d)\}$.\\              
Finally, let us study the effect of the dense-closed factorization on theories via the duality theorem.\\
Given a geometric theory $\mathbb T$ over a signature $\Sigma$ and a quotient ${\mathbb T}'$ of $\mathbb T$, let us describe the geometric theory ${\mathbb T}'^{\textrm{dc}}_{\mathbb T}$ over $\Sigma$ such that $\Sh({\cal C}_{\mathbb T}, J^{\mathbb T}_{{\mathbb T}'}) \hookrightarrow \Sh({\cal C}_{\mathbb T}, J^{\mathbb T}_{{\mathbb T}'^{\textrm{dc}}_{\mathbb T}}) \hookrightarrow \Sh({\cal C}_{\mathbb T}, J_{\mathbb T})$ is the dense-closed factorization of the inclusion $\Sh({\cal C}_{\mathbb T}, J^{\mathbb T}_{{\mathbb T}'}) \hookrightarrow \Sh({\cal C}_{\mathbb T}, J_{\mathbb T})$.\\ 
By equivalence (1) after the proof of Theorem \ref{prooftheory}, we have that $\emptyset \in J^{\mathbb T}_{{\mathbb T}'}(\{\vec{y}. \psi\})$ if and only if $\psi \vdash_{\vec{y}} \bot$ is provable in ${\mathbb T}'$. So, if $\phi$ is the geometric sentence corresponding to the subterminal identified with ${\cal C}_{0}^{J^{\mathbb T}_{{\mathbb T}'}}$ (equivalently, $\{[]. \phi\}\mono \{[]. \top\}$) is the union in ${\cal C}_{\mathbb T}$ of the images of all the arrows $\{\vec{y}. \psi\}\to \{[]. \top\}$ such that $\psi \vdash_{\vec{y}} \bot$ is provable in ${\mathbb T}'$, cfr. section \ref{opensub}) then, in view of the results in section \ref{subclosed}, we have that ${\mathbb T}'^{\textrm{dc}}_{\mathbb T}$ is obtained from $\mathbb T$ by adding the axiom 
\[
\psi \vdash_{\vec{y}} \psi'
\]
for any sequents $\psi' \vdash_{\vec{y}} \psi$ and $\psi \vdash_{\vec{y}} \psi' \vee (\phi \wedge \psi)$ which are provable in $\mathbb T$.   

\section{The surjection-inclusion factorization}\label{surincl}
We recall from \cite{El} (Theorem A4.2.10) that every geometric morphism can be factored, uniquely up to canonical equivalence, as a surjection followed by an inclusion. In this section we discuss the meaning of this factorization in terms of theories via the duality theorem.\\ 
Let us recall from the theory of classifying toposes that, given a geometric theory $\mathbb T$ over a signature $\Sigma$ with classifying topos $\cal E$, there exists a $\Sigma$-structure $M_{\mathbb T}$ in $\cal E$ which is `universal' among $\mathbb T$-models i.e. which satisfies the following property: $M$ is a $\mathbb T$-model and for any $\mathbb T$-model $N$ in a Grothendieck topos $\cal F$ there exists a unique (up to isomorphism) geometric morphism $f_{M}:{\cal F}\to {\cal E}$ such that $f_{M}^{\ast}(M_{\mathbb T})=N$. Thus, any geometric morphism $f$ into $\cal E$ is (up to isomorphism) of the form $f_{M}$ for a (unique up to isomorphism) $\mathbb T$-model $M$; indeed, $M\cong f^{\ast}(M_{\mathbb T})$.\\
Given a $\Sigma$-structure $M$ in a topos $\cal G$, let us define $Th(M)$ to be the theory over $\Sigma$ consisting of all the geometric sequents $\sigma$ over $\Sigma$ which hold in $M$; note that, by the soundess theorem for geometric logic, $Th(M)$ is a closed theory. 

\begin{theorem}\label{surinclthm}
Let $\mathbb T$ be a geometric theory over a signature $\Sigma$ and $f:{\cal F}\to {\cal E}$ be a geometric morphism into the classifying topos $\cal E$ for $\mathbb T$, corresponding to a $\mathbb T$-model $M$ in $\cal F$ as above. Then the topos ${\cal E}'$ in the surjection-inclusion factorization ${\cal F}\epi {\cal E}' \hookrightarrow {\cal E}$ of $f$ classifies the quotient $Th(M)$ of $\mathbb T$.
\end{theorem} 

\begin{proofs}
Let us denote by ${\cal F}\stackrel{f'}{\epi} {\cal E}' \stackrel{i}{\hookrightarrow} {\cal E}$ the surjection-inclusion factorization of $f$.
Since $i$ is a geometric inclusion to the classifying topos of $\mathbb T$, $i$ corresponds via the duality theorem to a unique closed quotient ${\mathbb T}'$ of $\mathbb T$ such that ${\cal E}'$ is a classifying topos of ${\mathbb T}'$. We want to prove that ${\mathbb T}'=Th(M)$. From the proof of Theorem \ref{dualita} we know that, if $M_{\mathbb T}$ is the universal model of $\mathbb T$ then $i^{\ast}(M_{\mathbb T})$ is a universal model $M_{{\mathbb T}'}$ for ${\mathbb T}'$. So $f'$ corresponds to the ${\mathbb T}'$-model $M$ via the universal property of the classifying topos of ${\mathbb T}'$, since $f'^{\ast}(M_{{\mathbb T}'})=f'^{\ast}(i^{\ast}(M_{\mathbb T}))\cong f^{\ast}(M_{\mathbb T})=M$. Now, since $f'$ is a surjection then, by Lemma D1.2.13 \cite{El2}, $M$ is a conservative ${\mathbb T}'$-model, from which it follows that ${\mathbb T}'=Th(M)$.        
\end{proofs}    

\begin{rmk}
\emph{The theorem implies that if $\mathbb T$ is a closed geometric theory over a (many-sorted) signature $\Sigma$ and $M$ is a conservative ${\mathbb T}$-model then $f_{M}$ is a surjection. Indeed, the subtopos of $\Set[{\mathbb T}]$ arising in the surjection-inclusion factorization of $f$ coincides with $\Set[{\mathbb T}]$, since it corresponds via Theorem \ref{dualita} to $Th(M)={\mathbb T}$. This result generalizes Corollary D3.2.6 \cite{El2}, which was proved under the assumption that $\Sigma$ be one-sorted.}  
\end{rmk}

\section{Atoms}\label{atoms}
In this section we describe the atoms of the lattice of subtoposes of a given elementary topos, that is the non-trivial toposes having no proper subtoposes (we recall that a topos $\cal E$ is said to be trivial if it is naturally equivalent to the category one having just one object and the idenity morphism on it, equivalently if it is degenerate i.e. $0_{\cal E}\cong 1_{\cal E}$). 
\begin{proposition}\label{atomic}
Let $\cal E$ be an elementary topos. Then the atoms of the lattice of subtoposes of $\cal E$ are exactly the two-valued Boolean subtoposes of $\cal E$.
\end{proposition}

\begin{proofs}
Our thesis follows as an immediate consequence of the following two facts. First, every non-trivial topos contains a non-trivial Boolean subtopos; second, a non-trivial Boolean topos does not contain any proper subtoposes if and only if it is two-valued. 
To prove the first assertion, we note that if $\cal E$ is non-trivial then $\sh_{\neg\neg}({\cal E})$ is again non-trivial; indeed, $1_{\cal E}$ clearly belongs to $\sh_{\neg\neg}({\cal E})$ while $0_{\cal E}$ belongs to $\sh_{\neg\neg}({\cal E})$ since $\neg\neg$ is a dense local operator on $\cal E$ (cfr. p. 219 \cite{El}), so if $\sh_{\neg\neg}({\cal E})$ is trivial then $0_{\cal E}\cong 1_{\cal E}$ i.e. $\cal E$ is trivial. The fact that $\sh_{\neg\neg}({\cal E})$ is Boolean is well-known (see for example Lemma A4.5.22 \cite{El}). This completes the proof of the first fact. It remains to prove the second assertion. Let us observe that, given two subterminal objects $U$ and $V$ in $\cal E$, the subtopos ${\cal E}\slash U \hookrightarrow {\cal E}$ is contained in the subtopos ${\cal E}\slash V \hookrightarrow {\cal E}$ if and only if $U\leq V$ in the lattice $\Sub_{\cal E}(1_{\cal E})$. Indeed, it follows from our discussion in section \ref{opensub} above that ${\cal E}\slash U \hookrightarrow {\cal E}$ factors through ${\cal E}\slash V \hookrightarrow {\cal E}$ if and only if the projection $U\times V \to U$ is isomorphic to the terminal object $1_{U}:U\to U$ in ${\cal E}\slash U$, and, since for any object there can be at most one morphism from it to a given subterminal object, this condition is equivalent to requiring that $U\leq V$ (equivalently, $U\leq U\times V$). Now, if $\cal E$ is Boolean then all the subtoposes of $\cal E$ are open (by Proposition A4.5.22 \cite{El}), so that we have a lattice isomorphism between $\Sub_{\cal E}(1_{\cal E})$ and the lattice of subtoposes of $\cal E$; therefore a non-trivial Boolen topos does not contain any proper subtoposes if and only if it is two-valued.
\end{proofs}  
\begin{rmk}\label{rmk5}
\emph{We note that if a Grothendieck topos $\cal E$ has enough points then $\cal E$ is Boolean and two-valued if and only if it is atomic and connected. Indeed, we know from Corollary D3.5.2 \cite{El2} that every Boolean Grothendieck topos with enough points is atomic, and an atomic topos is two-valued if and only if it is connected (cfr. the proof of Theorem 2.5 \cite{OC2})}   
\end{rmk}   

Now we want to understand, in view of Theorem \ref{dualita}, the meaning of Proposition \ref{atomic} in terms of theories. To this end, let us recall from \cite{OC5} some definitions.

\begin{definition}
Let $\mathbb{T}$ be a geometric theory over a signature $\Sigma$. $\mathbb T$ is said to be Boolean if it classifying topos is a Boolean topos. 
\end{definition}
Given two geometric formulae $\phi$ and $\psi$ over $\Sigma$ in the same context $\vec{x}$, we write $\phi \stackrel{\mathbb T}{\sim} \psi$ to mean that both the sequents $\phi \vdash_{\vec{x}} \psi$ and $\psi \vdash_{\vec{x}} \phi$ are provable in $\mathbb T$.\\  
\begin{rmk}\label{rmk1}
\emph{We recall from \cite{OC3} that a geometric theory $\mathbb T$ over a signature $\Sigma$ is a Boolean if and only if for every geometric formula $\phi(\vec{x})$ over $\Sigma$ there exists a geometric formula $\psi(\vec{x})$ over $\Sigma$ in the same context, denoted $\neg \phi(\vec{x})$, such that $\phi(\vec{x}) \wedge \psi(\vec{x})\stackrel{\mathbb T}{\sim}\bot$ and $\phi(\vec{x}) \vee \psi(\vec{x})\stackrel{\mathbb T}{\sim}\top$.\\
From this criterion, it easily follows that if $\mathbb T$ is a Boolean then every infinitary first-order formula over $\Sigma$ is $\mathbb T$-provably equivalent using classical logic to a geometric formula in the same context; indeed, this can be proved by an inductive argument as in the proof of Theorem D3.4.6 p. 921 \cite{El2} (in the case of an infinitary conjunction $\mathbin{\mathop{\textrm{\huge $\wedge$}}\limits_{i\in I}}\phi_{i}$, we observe that this formula is equivalent in classical logic to the formula $\neg(\mathbin{\mathop{\textrm{\huge $\vee$}}\limits_{i\in I}}\neg \phi_{i}$), where the symbol $\neg$ here denotes the first-order negation. Notice that from the fact that every infinitary first-order formula is classically equivalent in $\mathbb T$ to a geometric formula, it follows from the axioms of infinitary first-order logic for implication and infinitary conjunction that the first-order implication between geometric formulae is classically provably equivalent in $\mathbb T$ to the Heyting implication between them in the relevant subobject lattice of ${\cal C}_{\mathbb T}$, while the infinitary conjunction of a family of geometric formulae is classically provably equivalent in $\mathbb T$ to the infimum of the family in that lattice.} 
\end{rmk}

\begin{definition}
Let $\mathbb{T}$ be a geometric theory. $\mathbb{T}$ is said to be atomic if its classifying topos $\Set[\mathbb T]$ is an atomic topos.  
\end{definition}

\begin{definition}
Let $\mathbb{T}$ be a geometric theory over a signature $\Sigma$. $\mathbb T$ is said to have enough models if for every geometric sequent $\sigma$ over $\Sigma$, $M\vDash \sigma$ for all the $\mathbb T$-models $M$ in $\Set$ implies that $\sigma$ is provable in $\mathbb T$.  
\end{definition}
\begin{rmk}\label{rmk2}
\emph{It was observed in \cite{OC5} (Proposition 2.3) that a theory has enough models if and only if its classifying topos has enough points.}
\end{rmk} 

\begin{definition}
Let $\mathbb{T}$ be a geometric theory over a signature $\Sigma$. $\mathbb T$ is said to be complete if every geometric sentence $\phi$ over $\Sigma$ is $\mathbb T$-provably equivalent to $\top$ or $\bot$, but not both.  
\end{definition}
\begin{rmk}\label{rmk3}
\emph{A geometric theory $\mathbb T$ over a signature $\Sigma$ is complete if and only if its classifying topos is two-valued (i.e. it has exactly two subobjects of $1$); indeed, we observed in section \ref{opensub} that the subobjects of the classifying topos $\Set[{\mathbb T}]$ can be identified with the $\mathbb T$-provable equivalence classes of geometric sentences over $\Sigma$.}
\end{rmk}

\begin{definition}
Let $\mathbb{T}$ be a geometric theory over a signature $\Sigma$. $\mathbb{T}$ is said to be contradictory if $\top \vdash_{[]} \bot$ is provable in $\mathbb T$. 
\end{definition}
\begin{rmk}\label{rmk4}
\emph{A geometric theory is contradictory if and only if its classifying topos is trivial. Indeed, it is easy to verify that if $\mathbb T$ is contradictory then the trivial topos satisfies the universal property of the classifying topos of $\mathbb T$, and that, conversely, if the classifying topos of $\mathbb T$ is trivial then $\bot$ holds in it and hence $\top \vdash_{[]} \bot$ is provable in $\mathbb T$.}\\ 
\end{rmk}

The following proposition represents the translation of Proposition \ref{atomic} in terms of theories via Theorem \ref{dualita}.
\begin{proposition}
Let $\mathbb T$ be a geometric theory over a signature $\Sigma$. Then the non-contradictory quotients ${\mathbb T}'$ of $\mathbb T$ such that for every geometric sequent $\sigma$ over $\Sigma$ either $\sigma$ is provable in $\mathbb T$ or the theory ${\mathbb T}\cup \{\sigma\}$ is contradictory are exactly the Boolean and complete theories. 
\end{proposition}
\qed
\begin{rmk}
\emph{We note that the `if' direction in the proposition above can be easily proved without appealing to the duality theorem as follows. If $\mathbb T$ is Boolean then given a geometric sequent $\phi \: \vdash_{\vec{x}} \psi$ over $\Sigma$, it is clear that $\phi \: \vdash_{\vec{x}} \psi$ is provable in $\mathbb T$ if and only if the infinitary first-order sentence $\forall \vec{x}(\phi \to \psi)$ is. Now, by Remark \ref{rmk1}, this formula is $\mathbb T$-provably equivalent using classical logic to a geometric sentence, and this sentence is $\mathbb T$-provably equivalent to $\top$ or $\bot$ since $\mathbb T$ is complete.}
\end{rmk}  

\newpage
\section{Toposes with enough points}
A point of a Grothendieck topos $\cal E$ is a geometric morphism $p:\Set\to {\cal E}$; if ${\cal E}$ is the classifying topos $\Set[{\mathbb T}]$ of a geometric theory $\mathbb T$ then the points of ${\cal E}$ correspond precisely to the models of $\mathbb T$ in $\Set$. Let us recall from \cite{El2} that a Grothendieck topos $\cal E$ is said to have enough points if the inverse image functors $f^{\ast}$ of the geometric morphisms $f:\Set \rightarrow {\cal E}$ are jointly conservative. If $\cal E$ is the classifying topos $\Set[{\mathbb T}]$ of a geometric theory $\mathbb T$ over a signature $\Sigma$ then $\cal E$ has enough points if and only if $\mathbb T$ has enough models (cfr. Proposition 2.3 \cite{OC5}).\\
Recall that a model $M$ of a geometric theory $\mathbb T$ is said to be conservative if for any geometric sequent $\sigma$ over $\Sigma$, $M\vDash \sigma$ implies $\sigma$ provable in $\mathbb T$. Thus a geometric theory has enough models if and only if its $\Set$-models are jointly conservative.\\
Given a point $p$ of a topos $\cal E$, let us denote by ${\cal E}_{p}\hookrightarrow {\cal E}$ the inclusion part of the surjection-inclusion factorization of $p$. By Theorem \ref{surinclthm}, if ${\cal E}=\Set[{\mathbb T}]$ then ${\cal E}_{p}$ classifies $Th(M)$ where $M$ is the $\mathbb T$-model corresponding to $p$.\\    
Given a Grothendieck topos $\cal E$, let us define the \emph{subtopos ${\cal E}^{\textrm{points}}$ of points} of $\cal E$ to be the union of all the subtoposes ${\cal E}_{p}$ of $\cal E$ as $p$ varies among the points of $\cal E$ (such union exists because, dually, any intersection of Grothendieck topologies is a Grothendieck topology).\\
From Theorem \ref{dualita} and the description of the (infinitary) wedge in $\mathfrak{Th}_{\Sigma}^{\mathbb T}$, the topos $\Set[{\mathbb T}]^{\textrm{points}}$ classifies the intersection of all the theories $Th(M)$ as $M$ varies among the $\mathbb T$-models $M$ in $\Set$; in particular $\Set[{\mathbb T}]$ coincides with $\Set[{\mathbb T}]^{\textrm{points}}$ if and only if it has enough points. Notice that, obviously, any intersection in $\mathfrak{Th}_{\Sigma}^{\mathbb T}$ of theories of the form $Th(M)$ (for a $\mathbb T$-model $M$ in $\Set$) has enough models; in particular, all the toposes of the form $\Set[{\mathbb T}]^{\textrm{points}}$ have enough points. So we conclude that, given a geometric theory $\mathbb T$, the quotients of $\mathbb T$ having enough models are exactly the intersections in $\mathfrak{Th}_{\Sigma}^{\mathbb T}$ of theories of the form $Th(M)$ (where $M$ is a $\mathbb T$-model in $\Set$). Hence, since every Grothendieck topos is (equivalent to) the classifying topos of a geometric theory, we obtain the following equivalent topos-theoretic statement: the subtoposes of a Grothendieck topos $\cal E$ which have enough points are exactly the unions of subtoposes of the form ${\cal E}_{p}$ where $p$ is a point of $\cal E$.\\
Finally, we note that, given an atom ${\cal F}$ in the lattice of subtoposes of a Grothendieck topos ${\cal E}$ i.e. a Boolean and two-valued subtopos $\cal F$ of $\cal E$ (cfr. section \ref{atoms} above), if $\cal F$ has enough points then $\cal F$ is of the form ${\cal E}_{p}$ for a point $p$ of $\cal E$. Indeed, it is clear that a topos with enough points has a point if and only if it is non-trivial.      

\section{Skeletal inclusions}
Recall from \cite{El2} that a geometric morphism $f:{\cal F}\to {\cal E}$ is said to be skeletal if it restricts to a geometric morphism $\sh_{\neg\neg}({\cal F})\to \sh_{\neg\neg}({\cal E})$. By Lemma D4.6.10 \cite{El2}, a geometric inclusion $f:{\cal F}\to {\cal E}$ corresponding to a local operator $j$ on $\cal E$ is skeletal if and only if $ext(j)$ is a $\neg\neg$-closed subterminal object of $\cal E$.\\
Let us use the notation of section \ref{denseclosed} above. Given the canonical geometric inclusion $\Sh({\cal C}, J')\hookrightarrow \Sh({\cal C}, J)$ corresponding to an inclusion $J\subseteq J'$, $ext(\tau^{J}_{J'})$ identifies (as a subterminal object in $\Sh({\cal C}, J)$) with the $J$-ideal ${\cal C}_{0}^{J'}=\{c\in {\cal C} \textrm{ | } \emptyset \in J'(c)\}$. Now, consider the full subcategory $\tilde{\cal C}$ of $\cal C$ on the objects which are not $J$-covered by the empty sieve; $\tilde{\cal C}$ is $J$-dense in $\cal C$, and hence, by the Comparison Lemma, $\Sh({\cal C}, J)\simeq \Sh(\tilde{\cal C}, J|\tilde{{\cal C}})$, where $J|\tilde{{\cal C}}$ is the induced Grothendieck topology on $\tilde{\cal C}$. Moreover, $J|\tilde{{\cal C}}$ is dense i.e. $J|\tilde{{\cal C}}\leq \neg\neg_{[{\tilde{\cal C}}^{\textrm{op}}, \Set]}$. Thus, by Corollary \ref{cordense} and Theorem \ref{relativization}(ii), $ext(\tau^{J}_{J'})$ is $\neg\neg_{\Sh(\tilde{\cal C}, J|\tilde{{\cal C}})}$-closed (as a subterminal in $\Sh(\tilde{\cal C}, J|\tilde{{\cal C}})$) if and only if $ext(\tau^{J}_{J'})$ is $\neg\neg_{[{\tilde{\cal C}}^{\textrm{op}}, \Set]}$-closed (as a subterminal in $[{\tilde{\cal C}}^{\textrm{op}}, \Set]$). But $\neg\neg_{[{\tilde{\cal C}}^{\textrm{op}}, \Set]}$ is well-known to correspond to the dense topology on $\tilde{\cal C}$ i.e. to the Grothendieck topology on $\tilde{\cal C}$ whose covering sieves are exactly the stably non-empty ones; so, by formula (6) p. 235 \cite{MM}, we obtain that $ext(\tau^{J}_{J'})$ is $\neg\neg_{[{\tilde{\cal C}}^{\textrm{op}}, \Set]}$-closed if and only if for any $c\in \tilde{{\cal C}}$, `$\{f:d\to c \textrm{ in $\tilde{\cal C}$ | } d\in {\cal C}_{0}^{J'}\}$ stably non-empty in $\tilde{\cal C}$' implies `$c \in {\cal C}_{0}^{J'}$'.\\
Hence the geometric inclusion $\Sh({\cal C}, J')\hookrightarrow \Sh({\cal C}, J)$ is skeletal if and only if for any $c\in \tilde{{\cal C}}$, `$Z(c)=\{f:d\to c \textrm{ in $\tilde{\cal C}$ | } \emptyset \in J'(d)\}$ stably non-empty in $\tilde{\cal C}$' implies `$\emptyset \in J'(c)$'.\\
Now, let us interpret the meaning of the notion of skeletal inclusion at the level of theories, via the duality theorem. Specifically, given a geometric theory $\mathbb T$ over a signature $\Sigma$, let us describe the quotients ${\mathbb T}'$ of ${\mathbb T}$ such that the geometric inclusion $\Sh({\cal C}_{\mathbb T}, J_{{\mathbb T}'}^{\mathbb T})\hookrightarrow \Sh({\cal C}_{\mathbb T}, J_{\mathbb T})$ is skeletal.\\            
By the equivalence (1) after the proof of Theorem \ref{prooftheory}, we have that $\emptyset \in J^{\mathbb T}_{{\mathbb T}'}(\{\vec{y}. \psi\})$ if and only if $\psi \vdash_{\vec{y}} \bot$ is provable in ${\mathbb T}'$. Given an object $\{\vec{y}. \psi\}\in {\cal C}_{\mathbb T}$, let us denote by $\{\vec{y}. \psi_{{\mathbb T}'}\}\mono \{\vec{y}. \psi\}$ the subobject in ${\cal C}_{\mathbb T}$ given by the union in ${\cal C}_{\mathbb T}$ of all the subobjects $\{\vec{y}. \psi'\}\to \{\vec{y}. \psi\}$ such that $\psi' \vdash_{\vec{y}} \bot$ is provable in ${\mathbb T}'$. Then, recalling the results in \cite{OC3}, we obtain the following condition for $\Sh({\cal C}_{\mathbb T}, J_{{\mathbb T}'}^{\mathbb T})\hookrightarrow \Sh({\cal C}_{\mathbb T}, J_{\mathbb T})$ to be skeletal (below by a $\mathbb T$-consistent geometric formula we mean a geometric formula $\phi(\vec{x})$ such that $\phi \vdash_{\vec{x}} \bot$ is not provable in $\mathbb T$):\\
`for any geometric formula $\psi(\vec{y})$ over $\Sigma$, if $\psi_{{\mathbb T}'}(\vec{y})$ is $\mathbb T$-consistent and for any $\mathbb T$-consistent geometric formula $\chi(\vec{y})$ over $\Sigma$ such that $\chi \vdash_{\vec{y}} \psi$ is provable in $\mathbb T$, $(\chi \wedge \psi_{{\mathbb T}'})(\vec{y})$ is $\mathbb T$-consistent then $\psi \vdash_{\vec{y}} \bot$ is provable in ${\mathbb T}'$'.

\section{Some applications}

\subsection{Open and closed quotients}
Let $\mathbb T$ be a geometric theory over a signature $\Sigma$. Given an elementary topos $\cal E$, it is well-known that the open and closed subtoposes associated to a given subterminal object are complementary to each other in $\textbf{Lop}({\cal E})$. From this we deduce, by the duality theorem, that the open and closed quotients ${\mathbb T}_{\phi}$ and ${\mathbb T}^{\textrm{closed}}_{\phi}$ of $\mathbb T$ corresponding to a given geometric sentence $\phi$ are complementary to each other in $\mathfrak{Th}_{\Sigma}^{\mathbb T}$; note that this can also be proved directly by logical arguments. Also, we know from the theory of elementary toposes that if $U$ and $V$ are complemented subterminals in a topos $\cal E$ then $o(U)=c(V)$; this implies, by the duality theorem, that if $\phi$ and $\psi$ are two geometric sentences such that $\top \vdash_{[]} \phi \vee \psi$ and $\phi \wedge \psi \vdash_{[]} \bot$ then ${\mathbb T}_{\phi}={\mathbb T}^{\textrm{closed}}_{\psi}$; again, this can be easily proved directly by logical arguments.\\
Now, let us recall the following fact about elementary toposes (cfr. Proposition A4.5.22 \cite{El}): an elementary topos is Boolean if and only if every subtopos of it is open. It is interesting to interpret the `only if' part of this statement at the level of theories via the duality theorem.\\
If $\mathbb T$ is a Boolean geometric theory over a signature $\Sigma$ and ${\mathbb T}'$ is a quotient of $\mathbb T$, we want to show that there exists a geometric sentence $\phi$ over $\Sigma$ such that ${\mathbb T}'$ is syntactically equivalent to ${\mathbb T}_{\phi}$. For any axiom $\sigma=\phi \vdash_{\vec{x}} \psi$ of ${\mathbb T}'$, consider the geometric formula $U(\sigma)$ over $\Sigma$ classically equivalent in $\mathbb T$ (as in Remark \ref{rmk1}) to the infinitary first-order formula $\forall \vec{x}(\phi \to \psi)$. Now, there is only a set of such formulae $U(\sigma)$ over $\Sigma$ up to $\mathbb T$-provable equivalence, the geometric syntactic category ${\cal C}_{\mathbb T}$ being well-powered, so we can take $\phi$ to be a geometric sentence which is classically equivalent in $\mathbb T$ to their infinitary conjunction (as in Remark \ref{rmk1}); it is now immediate to see that $\phi$ has the required property.

\subsection{A deduction theorem for geometric logic}
The following result is the analogue for geometric logic of the deduction theorem in classical first-order logic; we will derive it by using our duality theorem.

\begin{theorem}
Let $\mathbb T$ be a geometric theory over a signature $\Sigma$ and $\phi, \psi$ two geometric sentences over $\Sigma$ such that the sequent $\top \vdash_{[]} \psi$ is provable in the theory ${\mathbb T}\cup \{\top \vdash_{[]} \phi\}$. Then the sequent $\phi \vdash_{[]} \psi$ is provable in the theory $\mathbb T$. 
\end{theorem}     
 
\begin{proofs}
By the duality theorem and Lemma D1.4.4 \cite{El2}, we can rephrase our thesis as follows: if $\{[].\psi\}\stackrel{[\psi]}{\mono} \{[].\top\}$ belongs to the Grothendieck topology generated by the $J_{\mathbb T}$-covering sieves and the principal sieve generated by $\{[].\phi\}\stackrel{[\phi]}{\mono} \{[].\top\}$, then $[\phi]\leq [\psi]$ in $\Sub_{{\cal C}_{\mathbb T}}(\{[].\top\})$.\\    
Now, by recalling that the syntactic topology $J_{\mathbb T}$ is the geometric topology on the category ${\cal C}_{\mathbb T}$ and Proposition \ref{loccat}, we can further rewrite our thesis as follows: if $\cal C$ is a geometric category and $J^{\textrm{geom}}_{\cal C}$ is the geometric topology on it then, given subobjects $m:a\mono 1$ and $n:b\mono 1$ of the terminal object $1$ in $\cal C$ such that $(n)$ belongs to the Grothendieck topology generated by the $J^{\textrm{geom}}_{\cal C}$-covering sieves and the sieve $(m)$, $m\leq n$ in $\Sub_{\cal C}(1)$.\\
Let us use the formula for the Grothendieck topology $(D^{r})^{l}$ generated by a family of sieves $D$ that is stable under pullback, which we obtained in section \ref{latticestructure}. Here we take $D$ to be the collection of all the sieves which are either $J^{\textrm{geom}}_{\cal C}$-covering or of the form $f^{\ast}((m))$ for a arrow $f$ with codomain $1$; so, starting from the assumption that $(n)\in (D^{r})^{l}(b)$, we want to deduce that $m\leq n$ in $\Sub_{\cal C}(1)$.\\
We note that $m\leq n$ if and only if $m\leq (m\imp n)$, if and only if $m^{\ast}(m\imp n)\cong 1_{a}$ (where $\imp$ denotes the Heyting implication in $\Sub_{\cal C}(1)$).
Now, from the simplified formula for $D^{l}$ we see that, since $n\leq (m\imp n)$, in order to prove that $m^{\ast}(m\imp n)\cong 1_{a}$ it suffices to show that $m^{\ast}((m\imp n))\in D^{r}(a)$ (in the formula one takes $Z$ to be $m^{\ast}((m\imp n))$, $S$ to be $(n)$ and $f$ to be $m$); in fact, we will prove that $(m\imp n)\in D^{r}(1)$, which implies that $m^{\ast}((m\imp n))\in D^{r}(a)$ since $D^{r}$ is stable under pullback.\\
By the simplified formula for $D^{r}$, we are reduced to prove that for any arrow $f:d\to 1$ with codomain $1$ and any sieve $S$ on $d$ such that $S\in D(d)$, $S\subset f^{\ast}((m\imp n))$ implies $1_{d}\in f^{\ast}((m\imp n))$. Now, if $S\in D(d)$ then there are two options: either $S$ is  $J^{\textrm{geom}}_{\cal C}$-covering or (since $1$ is a terminal object) $S$ is equal to $f^{\ast}((m))$. In the first case, we have that $f^{\ast}((m\imp n))$ is therefore $J^{\textrm{geom}}_{\cal C}$-covering and hence, being generated by a monomorphism, maximal, as required. In the second case, we have that $f^{\ast}(m)\leq f^{\ast}(m\imp n)$. But $f^{\ast}(m\imp n)=f^{\ast}(m)\imp f^{\ast}(n)$ (cfr. p. 41 \cite{El}) and hence $f^{\ast}(m)\leq f^{\ast}(m\imp n)$ implies $f^{\ast}(m)\leq f^{\ast}(n)$ i.e. $1_{d}\in f^{\ast}((m\imp n))$.        
\end{proofs}

\newpage
\section{The quotients of a theory of presheaf type}
In the first part of this paper, we have described the classifying topos of the quotient of a geometric theory in a syntactic way. Often, it is natural to present theories as quotients of a theory of presheaf type; as we shall see below, this approach has the advantage that, under appropriate hypotheses, it is possible to obtain a `semantic' representation for the classifying topos of the given quotient. The purpose of this section is in fact to discuss the relationship between these syntactic and semantic representations of a given classifying topos.\\

The notation in this section is borrowed from \cite{El2}.\\

Let us recall that an object $c$ of a finitely accessible category is said to be finitely presentable if the representable functor $Hom_{\cal C}(c,-):{\cal C}\to \Set$ preserves filtered colimits.\\

\begin{definition}
A geometric theory $\mathbb T$ is said to be of presheaf type if it is classified by a presheaf topos.
\end{definition}
\begin{rmk}
\emph{Note that a theory $\mathbb T$ is of presheaf type if and only if it is classified by the topos $[{\cal C},\Set]$, where ${\cal C}:=\textrm{f.p.} {\mathbb T}\textrm{-mod}(\Set)$ is the category of finitely presentable $\mathbb T$-models in $\Set$ i.e. the full subcategory of ${\mathbb T}\textrm{-mod}(\Set)$ on the finitely presentable objects. To prove this recall that, by Diaconescu's theorem, we have an equivalence of categories ${\mathbb T}\textrm{-mod}(\Set)\simeq {\bf Flat}({\cal C}^{\textrm{op}},\Set)=\Ind{\cal C}$. Hence the category ${\mathbb T}\textrm{-mod}(\Set)$ is finitely accessible and the Cauchy completion $\cal{\check{C}}$ of the category $\cal C$ is recoverable (up to equivalence) from $\Ind{\cal C}$ as the full subcategory $\cal{\check{C}} \simeq \textrm{f.p.} {\mathbb T}\textrm{-mod}(\Set)$ of finitely presentable objects (cfr. Proposition C4.2.2 \cite{El2}); but $[{\cal C}, \Set]$ and $[{\cal{\check{C}}}, \Set]$ are naturally equivalent (cfr. Corollary A1.1.9 \cite{El}), from which our claim follows. Thus, by Diaconescu's theorem, any theory of presheaf type $\mathbb T$ is Morita-equivalent to the theory of flat functors on $\textrm{f.p.} {\mathbb T}\textrm{-mod}(\Set)^{\textrm{op}}$, that is we have an equivalence of categories ${\mathbb T}\textrm{-mod}({\cal E}) \simeq {\bf Flat}(\textrm{f.p.} {\mathbb T}\textrm{-mod}(\Set)^{\textrm{op}}, {\cal E})$ natural in ${\cal E}\in \mathfrak{BTop}$.}  
\end{rmk}

\subsection{The axiomatization of homogeneous models with respect to a Grothendieck topology}
Let $\mathbb T$ be a theory of presheaf type, together with an equivalence $\xi^{\cal E}: {\bf Flat}(\textrm{f.p.} {\mathbb T}\textrm{-mod}(\Set)^{\textrm{op}}, {\cal E}) \to {\mathbb T}\textrm{-mod}({\cal E})$ natural in ${\cal E}\in \mathfrak{Btop}$. If $y:\textrm{f.p.} {\mathbb T}\textrm{-mod}(\Set) \to [\textrm{f.p.} {\mathbb T}\textrm{-mod}(\Set)^{\textrm{op}}, \Set]$ is the Yoneda embedding then the factorization of the composite $\xi^{\Set}\circ y:\textrm{f.p.} {\mathbb T}\textrm{-mod}(\Set) \to {\mathbb T}\textrm{-mod}({\Set})$ through the inclusion $i:\textrm{f.p.} {\mathbb T}\textrm{-mod}(\Set)\hookrightarrow {\mathbb T}\textrm{-mod}(\Set)$ is an equivalence $\tau^{\xi}:\textrm{f.p.} {\mathbb T}\textrm{-mod}(\Set) \to \textrm{f.p.} {\mathbb T}\textrm{-mod}(\Set)$.\\  

Let us recall from \cite{OC1} that, given a flat functor $F:\textrm{f.p.} {\mathbb T}\textrm{-mod}(\Set)^{\textrm{op}}\to {\cal E}$, we have the `Yoneda representation'
\[
F \circ \tau^{\xi} \cong  Hom_{{\mathbb T}\textrm{-mod}(\cal E)}^{\cal E}(\gamma^{\ast}_{\cal E}(i(-)),M_{F}),
\]
where $\gamma_{\cal E}:{\cal E}\to \Set$ is the unique geometric morphism from $\cal E$ to $\Set$ and $M_{F}$ is the $\mathbb T$-model in $\cal E$ corresponding to $F\in {\bf Flat}(\textrm{f.p.} {\mathbb T}\textrm{-mod}(\Set)^{\textrm{op}}, {\cal E})$ via the equivalence $\xi^{\cal E}: {\bf Flat}(\textrm{f.p.} {\mathbb T}\textrm{-mod}(\Set)^{\textrm{op}}, {\cal E}) \to {\mathbb T}\textrm{-mod}({\cal E})$.\\ 
We note that, given an equivalence $\xi$ for a theory of presheaf type $\mathbb T$ as above, we can modify $\xi$ so that $\tau^{\xi}$ becomes the identity on $\textrm{f.p.} {\mathbb T}\textrm{-mod}(\Set)$. Indeed, composing with $(\tau^{\xi})^{-1}$ gives rise to an equivalence $((-)\circ (\tau^{\xi})^{-1})^{\cal E}:{\bf Flat}(\textrm{f.p.} {\mathbb T}\textrm{-mod}(\Set)^{\textrm{op}}, {\cal E}) \to {\bf Flat}(\textrm{f.p.} {\mathbb T}\textrm{-mod}(\Set)^{\textrm{op}}, {\cal E}) $ natural in ${\cal E}\in \mathfrak{Btop}$, and it easily follows from the Yoneda representation and the Yoneda Lemma that the composite equivalence $\xi':=\xi \circ ((-)\circ (\tau^{\xi})^{-1})$ is such that $\tau^{\xi'}\cong 1_{\textrm{f.p.} {\mathbb T}\textrm{-mod}(\Set)}$. In fact, given a theory of presheaf type $\mathbb T$, we will assume below that $\mathbb T$ comes equipped with an equivalence $\xi$ satisfying the condition $\tau^{\xi'}\cong 1_{\textrm{f.p.} {\mathbb T}\textrm{-mod}(\Set)}$; we will call such an equivalence \emph{canonical}, and, accordingly, we will say that an equivalence $\chi^{\cal E}:{\mathbb T}\textrm{-mod}({\cal E}) \simeq {\bf Geom}({\cal E},[\textrm{f.p.} {\mathbb T}\textrm{-mod}(\Set), \Set])$ natural in ${\cal E}\in \mathfrak{Btop}$ is canonical if it is induced by a canonical equivalence $\xi^{\cal E}: {\bf Flat}(\textrm{f.p.} {\mathbb T}\textrm{-mod}(\Set)^{\textrm{op}}, {\cal E}) \to {\mathbb T}\textrm{-mod}({\cal E})$ by composition with Diaconescu's equivalence.\\     

Let us also recall from \cite{OC1} the following definition.
\begin{definition}
Let $\mathbb T$ be a theory of presheaf type, $\cal E$ a Grothendieck topos and $S$ a sieve in $\textrm{f.p.} {\mathbb T}\textrm{-mod}(\Set)^{\textrm{op}}$ on an object $c\in \textrm{f.p.} {\mathbb T}\textrm{-mod}(\Set)$. A model $M\in {\mathbb T}\textrm{-mod}(\cal E)$ is said to be $S\textrm{-homogeneous}$ if and only if for each object $E\in \cal E$ and arrow $y:E^{\ast}(\gamma^{\ast}_{\cal E}(i(c)))\rightarrow E^{\ast}(M)$ in ${\mathbb T}\textrm{-mod}({\cal E}\slash E)$ there exists an epimorphic family $(p_{f}:E_{f}\rightarrow E, f\in S)$ and for each arrow $f:c\rightarrow d$ in $S$ an arrow $u_{f}:E_{f}^{\ast}(\gamma^{\ast}_{\cal E}(i(d)))\rightarrow E_{f}^{\ast}(M)$ in ${\mathbb T}\textrm{-mod}({\cal E}\slash E)$ such that $p_{f}^{\ast}(y)=u_{f}\circ E_{f}^{\ast}(\gamma^{\ast}_{\cal E}(i(f)))$.\\
If $J$ is a Grothendieck topology on $\textrm{f.p.} {\mathbb T}\textrm{-mod}(\Set)^{\textrm{op}}$ then $M$ is said to be $J\textrm{-homogeneous}$ if it is $S\textrm{-homogeneous}$ for every $J$-covering sieve $S$.  
\end{definition}

Thus, from the Yoneda representation above, it follows that $F$ is $J$-continuous if and only $M_{F}$ is $J$-homogeneous. Specifically, we have the following result (Theorem 4.6 \cite{OC1}): given a theory of presheaf type $\mathbb T$, together with a canonical equivalence $\chi^{\cal E}:{\mathbb T}\textrm{-mod}({\cal E}) \simeq {\bf Geom}({\cal E},[\textrm{f.p.} {\mathbb T}\textrm{-mod}(\Set), \Set])$ natural in ${\cal E}\in \mathfrak{Btop}$, a Grothendieck topology $J$ on $\textrm{f.p.} {\mathbb T}\textrm{-mod}(\Set)^{\textrm{op}}$, and a quotient ${\mathbb T}'$ of $\mathbb T$ with the corresponding inclusions $i^{\cal E}_{{\mathbb T}'}:{{\mathbb T}'}\textrm{-mod}({\cal E})\hookrightarrow {\mathbb T}\textrm{-mod}({\cal E})$ as in Remark \ref{diagrammafond}, the diagram in $\mathfrak{Cat}$ 

\[  
\xymatrix {
{{\mathbb T}'}\textrm{-mod}({\cal E}) \ar[rrr]^{\simeq} \ar[d]^{i^{\cal E}_{{\mathbb T}'}} & & & {\bf Geom}({\cal E},\Sh(\textrm{f.p.} {\mathbb T}\textrm{-mod}(\Set)^{\textrm{op}},J)) \ar[d]^{i\circ -} \\
{\mathbb T}\textrm{-mod}({\cal E}) \ar[rrr]^{\simeq}_{\chi^{\cal E}} & & & {\bf Geom}({\cal E},[\textrm{f.p.} {\mathbb T}\textrm{-mod}(\Set), \Set])}
\]
commutes (up to invertible natural equivalence) naturally in ${\cal E}\in \mathfrak{BTop}$ if and only if the ${\mathbb T}'$-models are exactly the $J$-homogeneous $\mathbb T$-models in every ${\cal E}\in \mathfrak{Btop}$.\\

The following theorem implies that $J$-homogeneous models are always axiomatizable by geometric sequents in the signature of $\mathbb T$.

\begin{theorem}\label{teoaxioms}
Let $\mathbb T$ be a theory of presheaf type and $J$ a Grothendieck topology on $\textrm{f.p.} {\mathbb T}\textrm{-mod}(\Set)^{\textrm{op}}$. Then there exists a (unique up to syntactic equivalence) geometric quotient ${\mathbb T}'$ of $\mathbb T$ such that the ${\mathbb T}'$-models are exactly the $J$-homogeneous $\mathbb T$-models in every Grothendieck topos.  
\end{theorem}

\begin{proofs}
Via the equivalence $[\textrm{f.p.} {\mathbb T}\textrm{-mod}(\Set), \Set] \simeq \Sh({\cal C}_{\mathbb T}, J_{\mathbb T})$, given by the uniqueness (up to equivalence) of the classifying topos of $\mathbb T$, the geometric inclusion $\Sh({\textrm{f.p.} {\mathbb T}\textrm{-mod}(\Set)}^{\textrm{op}}, J) \hookrightarrow [\textrm{f.p.} {\mathbb T}\textrm{-mod}(\Set), \Set]$ corresponds to a subtopos of $\Sh({\cal C}_{\mathbb T}, J_{\mathbb T})$, and hence the (closed) quotient of $\mathbb T$ corresponding to this inclusion via the duality theorem axiomatizes the $J$-homogeneous $\mathbb T$-models, by Remark \ref{diagrammafond} and the discussion preceding Theorem \ref{teoaxioms}.
\end{proofs}
In some cases of interest one can easily obtain an explicit axiomatization of the quotient ${\mathbb T}'$ in the theorem.
For example, if the category ${\textrm{f.p.} {\mathbb T}\textrm{-mod}(\Set)}^{\textrm{op}}$ satisfies the right Ore condition and $J_{at}$ is the atomic topology on it, then the geometric inclusion $\Sh({\textrm{f.p.} {\mathbb T}\textrm{-mod}(\Set)}^{\textrm{op}}, J_{at}) \hookrightarrow \textrm{f.p.} {\mathbb T}\textrm{-mod}(\Set)$ corresponds to the subtopos $\sh_{\neg\neg}({\Sh({\cal C}_{\mathbb T}, J_{\mathbb T})})$ of $\Sh({\cal C}_{\mathbb T}, J_{\mathbb T})$, and hence the $J_{at}$-homogeneous models are axiomatized by the Booleanization of $\mathbb T$ (cfr. \cite{OC3}).\\
Analogously, one can achieve a syntactic description of the geometric quotient of $\mathbb T$ corresponding to the De Morgan topology on the category ${\textrm{f.p.} {\mathbb T}\textrm{-mod}(\Set)}^{\textrm{op}}$; this is the DeMorganization of $\mathbb T$, as it is defined in \cite{OC3}.  

\newpage
\subsection{Finitely presented models of a theory of presheaf type}

The following definition will be central in this section.

\begin{definition}
Let $\mathbb T$ be a geometric theory over a signature $\Sigma$ and $\phi(x_{1}^{A_{1}}, \ldots, x_{n}^{A_{n}})$ be a geometric formula over $\Sigma$. We say that a $\mathbb T$-model $M$ in \Set is finitely presented by $\phi$ (or that $\phi$ presents $M$) if there exists a string of elements $(\xi_{1}, \ldots, \xi_{n})\in MA_{1}\times \ldots MA_{n}$, called the generators of $M$, such that for any $\mathbb T$-model $N$ in \Set and string of elements $(b_{1}, \ldots, b_{n})\in MA_{1}\times \ldots MA_{n}$ such that $(b_{1}, \ldots, b_{n})\in [[\phi]]_{N}$, there exists a unique arrow $f_{(b_{1}, \ldots, b_{n})}:M\to N$ in ${\mathbb T}\textrm{-mod}(\Set)$ such that $(fA_{1}\times \ldots fA_{n})((\xi_{1}, \ldots, \xi_{n}))=(b_{1}, \ldots, b_{n})$.      
\end{definition}
Of course, there can be at most one (up to isomorphism) $\mathbb T$-model finitely presented by a given formula $\phi$; we will denote such model by $M_{\phi}$.\\
Given a geometric theory $\mathbb T$ over a signature $\Sigma$ and a geometric formula $\phi(x_{1}^{A_{1}}, \ldots, x_{n}^{A_{n}})$ over $\Sigma$, let us consider the functor $F_{\phi}: {\mathbb T}\textrm{-mod}(\Set) \to \Set$ which sends to each model $N\in {\mathbb T}\textrm{-mod}(\Set)$ (the domain of) the interpretation $[[\phi]]_{N}$ of $\phi$ in $N$ and acts on arrows in the obvious way. The functor $F_{\phi}$ preserves filtered colimits (cfr. the proof of Lemma D2.4.9 \cite{El2}) so if it is representable then the representing object is a finitely presentable model. Notice that, by the Yoneda Lemma, $F_{\phi}$ is representable if and only if there exists a $\mathbb T$-model finitely presented by $\phi$. From this it follows that every finitely presented model of a geometric theory $\mathbb T$ is finitely presentable; the converse is always true if $\mathbb T$ is cartesian (cfr. pp. 882-883 \cite{El2}), but not in general (cfr. the coherent theory of fields in \cite{Diers}).\\
Suppose that $\mathbb T$ is a theory of presheaf type and ${\mathbb T}'$ is a quotient of $\mathbb T$ obtained from $\mathbb T$ by adding axioms $\sigma$ of the form $\phi \vdash_{\vec{x}} \mathbin{\mathop{\textrm{\huge $\vee$}}\limits_{i\in I}}(\exists \vec{y_{i}})\theta_{i}$, where, for any $i\in I$, $[\theta_{i}]:\{\vec{y_{i}}. \psi\}\to \{\vec{x}. \phi\}$ is an arrow in ${\cal C}_{\mathbb T}$ and $\phi(\vec{x})$, $\psi(\vec{y_{i}})$ are formulae presenting respectively ${\mathbb T}$-models $M_{\phi}$ and $M_{\psi_{i}}$.\\  
For each such axiom $\phi \vdash_{\vec{x}} \mathbin{\mathop{\textrm{\huge $\vee$}}\limits_{i\in I}}(\exists \vec{y_{i}})\theta_{i}$, consider the cosieve $S_{\sigma}$ on $M_{\phi}$ in $\textrm{f.p.} {\mathbb T}\textrm{-mod}(\Set)$ defined as follows. For each $i\in I$, $[[\theta_{i}]]_{M_{\psi_{i}}}$ is the graph of a morphism $[[\vec{y_{i}}.\psi_{i}]]_{M_{\psi_{i}}}\to [[\vec{x}.\phi]]_{M_{\psi_{i}}}$; then the image of the generators of $M_{\psi_{i}}$ via this morphism is an element of $[[\vec{x}.\phi]]_{M_{\psi_{i}}}$ and this in turn determines, by definition of $M_{\phi}$, a unique arrow $s_{i}:M_{\phi} \to M_{\psi_{i}}$ in ${\mathbb T}\textrm{-mod}(\Set)$.
We define $S_{\sigma}$ as the sieve in ${\textrm{f.p.} {\mathbb T}\textrm{-mod}(\Set)}^{\textrm{op}}$ on $M_{\phi}$ generated by the arrows $s_{i}$ as $i$ varies in $I$.\\
Let $F:{\textrm{f.p.} {\mathbb T}\textrm{-mod}(\Set)}^{\textrm{op}}\to {\cal E}$ be a flat functor; if $M_{\phi}\in {\textrm{f.p.} {\mathbb T}\textrm{-mod}(\Set)}$ is a finitely presented $\mathbb T$-model then $F(M_{\phi})=[[\phi]]_{M_{F}}$ where $M_{F}$ is the $\mathbb T$-model in $\cal E$ corresponding to $F$ via the Morita-equivalence. Indeed, denoted by $g:{\cal E}\to \Set[{\mathbb T}]\simeq [\textrm{f.p.} {\mathbb T}\textrm{-mod}(\Set), \Set]$ the geometric morphism corresponding to $F$ via the universal property of the classifying topos, we have that $F=g^{\ast}\circ y$ where $y:\textrm{f.p.} {\mathbb T}\textrm{-mod}(\Set)^{\textrm{op}}\to [\textrm{f.p.} {\mathbb T}\textrm{-mod}(\Set), \Set]$ is the Yoneda embedding; but $M_{F}=g^{\ast}(M_{\mathbb T})$ where $M_{\mathbb T}$ is the universal model of $\mathbb T$ lying in the classifying topos $\Set[{\mathbb T}]\simeq [\textrm{f.p.} {\mathbb T}\textrm{-mod}(\Set), \Set]$, and the representable $Hom(M_{\phi}, -)\in [\textrm{f.p.} {\mathbb T}\textrm{-mod}(\Set), \Set]$ is clearly (isomorphic to) $[[\phi]]_{M_{\mathbb T}}$. So, since inverse image functors of geometric morphisms preserve the interpretations of all geometric formulae, it follows that $F(M_{\phi})=[[\phi]]_{M_{F}}$, as required. It is also immediate to see that if $s_{i}:M_{\phi} \to M_{\psi_{i}}$ is an arrow in ${\textrm{f.p.} {\mathbb T}\textrm{-mod}(\Set)}$ induced as above by an arrow $[\theta_{i}]:\{\vec{y_{i}}. \psi\}\to \{\vec{x}. \phi\}$ in ${\cal C}_{\mathbb T}$ then $F(s_{i})=[[\theta]]_{M_{F}}$.\\
Given a geometric theory $\mathbb T$ over a signature $\Sigma$ and a geometric formula $\phi(x_{1}^{A_{1}}, \ldots, x_{n}^{A_{n}})$ over $\Sigma$, let us consider the functor $F^{\cal E}_{\phi}: {\mathbb T}\textrm{-mod}({\cal E}) \to {\cal E}$ which sends to each model $N\in {\mathbb T}\textrm{-mod}({\cal E})$ (the domain of) the interpretation $[[\phi]]_{N}$ of $\phi$ in $N$ and acts on arrows in the obvious way. If $\mathbb T$ is of presheaf type and $M_{\phi}$ is a $\mathbb T$-model model finitely presented by $\phi$ then $F^{\cal E}_{\phi}$ is $\cal E$-representable with representing object $\gamma^{\ast}_{\cal E}(i(M_{\phi}))$. Indeed, if $N\in {\mathbb T}\textrm{-mod}({\cal E})$ then from the Yoneda representation of the corresponding flat functor $F_{N}$ and the discussion above it follows that 
\[
F^{\cal E}_{\phi}(N)=[[\phi]]_{N}=F_{N}(M_{\phi}) \cong Hom_{{\mathbb T}\textrm{-mod}(\cal E)}^{\cal E}(\gamma^{\ast}_{\cal E}(i(M_{\phi})),N),
\]             
so that for any $E\in {\cal E}$ arrows $E\to [[\phi]]_{N}$ in $\cal E$ are in bijection with arrows $E^{\ast}(\gamma^{\ast}_{\cal E}(M_{\phi}))\to E^{\ast}(N)$ in ${\mathbb T}\textrm{-mod}({\cal E}\slash E)$.\\
Now, coming back to our sieve $S_{\sigma}$, it is clear that a model $N\in {\mathbb T}\textrm{-mod}({\cal E})$ is $S_{\sigma}$-homogeneous if and only if the sequent $\sigma$ holds in $N$; indeed, this follows directly from the discussion above by using Kripke-Joyal semantics, or alternatively by using that $N$ is $S_{\sigma}$-homogeneous if and only if $F_{N}$ sends $S_{\sigma}$ to an epimorphic family, if and only $\sigma$ holds in $N$. These remarks lead us to the following result.

\begin{theorem}\label{constru}
Let $\mathbb T$ be a theory of presheaf type such that all the finitely presentable $\mathbb T$-models in $\Set$ are finitely presented, and ${\mathbb T}'$ a quotient of $\mathbb T$ obtained from $\mathbb T$ by adding axioms $\sigma$ of the form $\phi \vdash_{\vec{x}} \mathbin{\mathop{\textrm{\huge $\vee$}}\limits_{i\in I}}(\exists \vec{y_{i}})\theta_{i}$, where, for each $i\in I$, $[\theta_{i}]:\{\vec{y_{i}}. \psi\}\to \{\vec{x}. \phi\}$ is an arrow in ${\cal C}_{\mathbb T}$ and $\phi(\vec{x})$, $\psi(\vec{y_{i}})$ are geometric formulae over the signature of $\mathbb T$ presenting respectively ${\mathbb T}$-models $M_{\phi}$ and $M_{\psi_{i}}$. With the notation above, if the collection of sieves $S_{\sigma}$ where $\sigma$ varies among the axioms of ${\mathbb T}'$ over $\mathbb T$ is stable under pullback then ${\mathbb T}'$ is classified by the topos $\Sh(\textrm{f.p.} {\mathbb T}\textrm{-mod}(\Set)^{\textrm{op}}, J)$ where $J$ is the Grothendieck topology on $\textrm{f.p.} {\mathbb T}\textrm{-mod}(\Set)^{\textrm{op}}$ generated by the sieves $S_{\sigma}$. 
\end{theorem}  

\begin{proofs}
This follows immediately from our discussions above (in particular, that preceding Theorem \ref{teoaxioms} together with Remark \ref{diagrammafond}) by using Lemma 3 \cite{blasce}.
\end{proofs}

\begin{rmk}
\emph{Our theorem generalizes the method of construction of the classifying topos of a quotient of a cartesian theory given by Propositions D3.1.7 and D3.1.10 \cite{El2}; indeed, it is well known (cfr. \cite{El2}) that the opposite of the category of finitely presentable models of a cartesian theory is equivalent (in the obvious way) to the cartesian syntactic category of the theory.}
\end{rmk}

Concerning the applicability of the theorem, we have seen above that, given geometric formulae $\phi(\vec{x})$ and $\psi(\vec{y})$ with finitely presented $\mathbb T$-models $M_{\phi(\vec{x})}$ and $N_{\psi(\vec{y})}$, any arrow $[\theta]:\{\vec{x}.\phi\}\to \{\vec{y}. \psi\}$ in the syntactic category ${\cal C}_{\mathbb T}$ gives rise to an arrow $N_{\psi}\to M_{\phi}$ in ${\mathbb T}\textrm{-mod}({\Set})$. If all the finitely presentable $\mathbb T$-models in $\Set$ are finitely presented and moreover all the homomorphisms of finitely presented $\mathbb T$-models arise in this way, then we say that the category $\textrm{f.p.} {\mathbb T}\textrm{-mod}(\Set)$ is syntactically presented; note that every cartesian theory satisfies this condition, by the results in \cite{El2}, and also the theory of undirected graphs p. 907 \cite{El2} and the theory of decidably linearly ordered objects p. 926 \cite{El2} enjoy it. If this condition is satisfied then we know from the proof of Theorem \ref{dualita} that it is superfluous to require the condition that the collection of sieves $S_{\sigma}$ ${\mathbb T}'$ should be stable under pullback, since we can always achieve it without modifying the syntactic-equivalence class of the theory ${\mathbb T}'$.\\
We remark that for theories $\mathbb T$ of presheaf type such that the category $\textrm{f.p.} {\mathbb T}\textrm{-mod}(\Set)$ is syntactically presented, every small presieve on $\textrm{f.p.} {\mathbb T}\textrm{-mod}(\Set)^{\textrm{op}}$ is of the form $S_{\sigma}$ for some geometric sequent in the signature of $\mathbb T$, so that, by the arguments in the proof of Theorem \ref{prooftheory}(ii), we can obtain axiomatizations of the quotient of $\mathbb T$ given by Theorem \ref{teoaxioms} starting from a collection of presieves on $\textrm{f.p.} {\mathbb T}\textrm{-mod}(\Set)^{\textrm{op}}$ which generates a given Grothendieck topology on $\textrm{f.p.} {\mathbb T}\textrm{-mod}(\Set)^{\textrm{op}}$, as in the following result.

\begin{theorem}\label{contruinv}
Let $\mathbb T$ be a theory of presheaf type such that the category $\textrm{f.p.} {\mathbb T}\textrm{-mod}(\Set)$ is syntactically presented and $J$ be a Grothendieck topology on $\textrm{f.p.} {\mathbb T}\textrm{-mod}(\Set)^{\textrm{op}}$. If a collection of presieves of the form $S_{\sigma}$ generates $J$ then the quotient ${\mathbb T}'$ of $\mathbb T$ corresponding to $J$ via Theorem \ref{teoaxioms} is axiomatized over $\mathbb T$ by the collection of the sequents $\sigma$; in particular, ${\mathbb T}'$ axiomatizes the $J$-homogeneous $\mathbb T$-models in every Grothendieck topos.    
\end{theorem}\qed
Note that this theorem also formally follows from Theorem \ref{constru} by using the duality theorem.\\

By applying the Theorem \ref{contruinv} in the case of the atomic topology we get the following result.

\begin{corollary}
Let $\mathbb T$ be a theory of presheaf type such that the category $\textrm{f.p.} {\mathbb T}\textrm{-mod}(\Set)$ is syntactically presented and $\textrm{f.p.} {\mathbb T}\textrm{-mod}(\Set)^{\textrm{op}}$ satisfies the right Ore condition. Then the theory ${\mathbb T}'$ corresponding to the atomic topology on $\textrm{f.p.} {\mathbb T}\textrm{-mod}(\Set)^{\textrm{op}}$ via Theorem \ref{teoaxioms} is obtained from $\mathbb T$ by adding
all the axioms of the form $\phi \vdash_{\vec{x}}(\exists \vec{y})\theta$, where $[\theta]:\{\vec{y}. \psi\}\to \{\vec{x}. \phi\}$ is any arrow in ${\cal C}_{\mathbb T}$ and $\phi(\vec{x})$, $\psi(\vec{y})$ are geometric formulae over the signature of $\mathbb T$ presenting respectively ${\mathbb T}$-models $M_{\phi}$ and $M_{\psi}$.
\end{corollary}  

\begin{proofs}
This follows from the theorem and Theorem \ref{teoaxioms} by observing that the collection of presieves on $\textrm{f.p.} {\mathbb T}\textrm{-mod}(\Set)^{\textrm{op}}$ formed by a single morphism generates the atomic topology on $\textrm{f.p.} {\mathbb T}\textrm{-mod}(\Set)^{\textrm{op}}$.  
\end{proofs}

In particular, we note that in the above Corollary all the axioms of the form $\phi \vdash_{\vec{x}} (\exists \vec{y})\psi$, where $\phi(\vec{x})$ and $\psi(\vec{x}, \vec{y})$ are geometric formulae over the signature of $\mathbb T$ presenting $\mathbb T$-models $M_{\phi}$ and $M_{\psi}$ and such that $\psi \vdash_{\vec{x}, \vec{y}} \phi$ is provable in $\mathbb T$, are provable in ${\mathbb T}'$.\\

We remark that if $\mathbb T$ is cartesian then the hypotheses of the Corollary are always satisfied. In this case, by recalling that the finitely presentable $\mathbb T$-models in $\Set$ are exactly those of the form $M_{\phi}$ for a cartesian formula $\phi$ and that the association of $M_{\phi}$ to $\phi$ defines an equivalence of categories ${\cal C}^{\textrm{cart}}_{\mathbb T} \simeq \textrm{f.p.} {\mathbb T}\textrm{-mod}(\Set)^{\textrm{op}}$, we obtain that the quotient ${\mathbb T}'$ over $\mathbb T$ in the Corollary is obtained from $\mathbb T$ by adding all the axioms of the form $\phi \vdash_{\vec{x}} (\exists \vec{y})\theta$, where $\phi(\vec{x})$ and $\theta(\vec{y}, \vec{x})$ are cartesian formulae over the signature of $\mathbb T$ such that the sequents $(\psi \vdash_{\vec{y}, \vec{x}} \phi)$ and $((\theta \wedge \theta[\vec{x'}\slash \vec{x}]) \vdash_{\vec{y}, \vec{x}, \vec{x'}} (\vec{x}=\vec{x'}))$ are provable in $\mathbb T$.

\section{Classifying toposes for theories with enough models}

In this section we extend some ideas and results from sections 3.2 and 3.4 of \cite{prest}, by rewriting them into a general topos-theoretic context; among other things, this will lead, under appropriate hypotheses, to a model-theoretic representation for the classifying topos of a quotient of a theory of presheaf type having enough models.\\
First, let us recall the definition of finitely accessible category and some basic facts which will be useful for our analysis; we refer the reader to section C4.2 \cite{El2} for the background.\\
A finitely accessible category $\cal L$ is a category which is equivalent to the Ind-completion $\Ind{\cal C}$ of a small category $\cal C$; $\Ind{\cal C}$ is defined to be the full subcategory of $[{\cal C}^{\textrm{op}}, \Set]$ on the flat functors $F:{\cal C}^{\textrm{op}} \to \Set$; recall that a functor $F:{\cal C}^{\textrm{op}}\to \Set$ is flat if it is a filtered colimit of representables, that is if the category of elements $\int_{\cal C} F$ of $F$ is filtered (recall that any presheaf $F$ is the colimit in $[{\cal C}^{\textrm{op}}, \Set]$ of the functor given by the composite $\int_{\cal C} F \stackrel{\pi}{\to} {\cal C} \stackrel{y}{\to} [{\cal C}^{\textrm{op}}, \Set]$ where $\pi$ is the obvious projection map and $y$ is the Yoneda embedding). Every representable functor is flat, so the Yoneda embedding $y:{\cal C}\to [{\cal C}^{\textrm{op}}, \Set]$ factors through the embedding $\Ind{\cal C} \hookrightarrow [{\cal C}^{\textrm{op}}, \Set]$; we will denote this factorization by $\overline{y}_{\cal C}:{\cal C}\to \Ind{\cal C}$. Moreover, the inclusion $\Ind{\cal C}\hookrightarrow [{\cal C}^{\textrm{op}}, \Set]$ creates filtered colimits.\\   
Given a finitely accessible category $\cal L$, we define $\textrm{f.p.}{\cal L}$ as the full subcategory of $\cal L$ on the finitely presentable objects; then the embedding $\textrm{f.p.}{\cal L}\hookrightarrow {\cal L}$ is (up to equivalence) of the form $\overline{y}_{\textrm{f.p.}{\cal L}}$ (cfr. Proposition C4.2.2 \cite{El2} and Corollary A1.1.9 \cite{El}).\\
We recall from \cite{El2} (Corollary C4.2.6) that the Ind-completion $\Ind{\cal C}$ of $\cal C$ is the free filtered-colimit completion of $\cal C$, that is, for any category $\cal D$ with filtered colimits, any functor $F:{\cal C}\to {\cal D}$ extends, via $\overline{y}_{\cal C}:{\cal C}\to \Ind{\cal C}$, uniquely up to canonical isomorphism, to a filtered-colimit preserving functor $\overline{F}:\Ind{\cal C}\to {\cal D}$.\\  
Now, generalizing \cite{prest}, given a small category $\cal C$, we construct correspondences between the collection ${\cal S}_{\Ind{\cal C}}$ of full subcategories of $\Ind{\cal C}$ and the collection ${\cal G}({\cal C})$ of Grothendieck cotopologies on $\cal C$. Given a cosieve $S$ in $\cal C$ on an object $c\in {\cal C}$, we denote by $\overline{S}$ the extension of $S:{\cal C}\to \Set$ (regarded here as a subfunctor of the representable ${\cal C}(c, -)$) along $\overline{y}_{\cal C}$ as above.\\
Let us define correspondences ${\cal K}:{\cal G}({\cal C}) \to {\cal S}_{\Ind{\cal C}}$ and ${\cal H}:{\cal S}_{\Ind{\cal C}}\to {\cal G}({\cal C})$ as follows.\\
Given a Grothendieck cotopology $J$ on $\cal C$, ${\cal K}(J)$ is the full subcategory of $\Ind{\cal C}$ defined by
\[
\begin{array}{ccl}
d\in {\cal K}(J) & \textrm{iff} & \overline{S}(d)=Hom_{\Ind{\cal C}}(\overline{y}_{\cal C}(c), d) \textrm{ for all } S\in J(c),  
\end{array}
\]
for any $d\in \Ind{\cal C}$.
Conversely, given a full subcategory $\cal D$ of $\Ind{\cal C}$, we define ${\cal H}({\cal D})$ by
\[
\begin{array}{ccl}
S\in {\cal H}({\cal D})(c) & \textrm{iff} & \overline{S}(d)=Hom_{\Ind{\cal C}}(\overline{y}_{\cal C}(c), d) \textrm{ for all } d\in {\cal D},  
\end{array}
\]
for any cosieve $S$ in $\cal C$ on an object $c\in {\cal C}$. Here by the equality $\overline{S}(d)=Hom_{\Ind{\cal C}}(\overline{y}_{\cal C}(c), d)$ we mean that the values at $d$ of the functors $\overline{S}$ and $Hom_{\Ind{\cal C}}(\overline{y}_{\cal C}(c), -)$ are canonically isomorphic i.e. (by the description of filtered colimits in $\Set$ p. 77 \cite{borceux}) if $d=colim(\overline{y}_{\cal C}\circ G)$ in $\Ind{\cal C}$ where $I$ is a filtered category and $G:I\to \cal C$ is a functor then for any arrow $r:c\to G(i)$ there exist objects $j, k\in I$ and arrows $s:c\to G(j)$ in $S$ and $\chi:i\to k$, $\xi:j\to k$ in $I$ such that $G(\chi)\circ r=G(\xi)\circ s$.\\
It is easy to verify that for any full subcategory $\cal D$ of $\Ind{\cal C}$, ${\cal H}({\cal D})$ is indeed a Grothendieck cotopology on ${\cal C}$; we provide the details for the reader's convenience.\\
It is clear that the maximality axiom holds. Let us verify the stability axiom. Given an arrow $f:c\to c'$ in $\cal C$ and a cosieve $S\in {\cal H}({\cal D})(c)$, we want to prove that $f^{\ast}(S)\in {\cal H}({\cal D})(c')$, that is for any arrow $r:c'\to G(i)$ there exists $j, k\in I$ and arrows $s:c'\to G(j)$ in $f^{\ast}(S)$ and $\chi:i\to k$, $\xi:j\to k$ in $I$ such that $G(\chi)\circ r=G(\xi)\circ s$. Consider the arrow $r\circ f$; since $S\in {\cal H}({\cal D})(c)$ then there exist $j', k'\in I$ and arrows $s':c\to G(j')$ in $S$ and $\chi':i\to k'$, $\xi:j'\to k'$ in $I$ such that $G(\chi')\circ r\circ f=G(\xi')\circ s'$. Then if we take $i=i'$, $j=k$, $\chi=\chi'$, $\xi=1_{k}$, $s=G(\chi)\circ r$, we have that $s\in f^{\ast}(S)$ and hence our thesis is satisfied. It remains to verify the transitivity axiom.
Given a cosieve $R$ on $c\in {\cal C}$ and a cosieve $S\in {\cal H}({\cal D})(c)$ such that $f^{\ast}(R)\in {\cal H}({\cal D})(cod(f))$ for any $f\in S$, we want to prove that $R\in {\cal H}({\cal D})(c)$. Since $S\in {\cal H}({\cal D})(c)$, given an arrow $r:c\to G(i)$ there exist $j, k\in I$ and arrows $f:c\to G(j)$ in $S$ and $\chi:i\to k$, $\xi:j\to k$ in $I$ such that $G(\chi)\circ r=G(\xi)\circ f$; now, since $f^{\ast}(R)\in {\cal H}({\cal D})(cod(f))$, there there exist $j', k'\in I$ and arrows $g:G(j)\to G(j')$ in $f^{\ast}(R)$ and $\chi':k\to k'$, $\xi':j'\to k'$ in $I$ such that $G(\chi')\circ G(\xi)=G(\xi')\circ g$; hence $G(\chi'\circ \chi)\circ r=G(\xi')\circ g\circ f$ and our thesis is satisfied.\\
Next, we note that ${\cal S}_{\Ind{\cal C}}$ and ${\cal G}({\cal C})$ are naturally equipped with partial orders (respectively the obvious inclusion between full subcategories of ${\cal S}_{\Ind{\cal C}}$ and the inclusion between Grothendieck cotopologies on ${\cal C}$) and if we regard them as poset categories then the correspondences $\cal H$ and $\cal K$ become contravariant functors; moreover, it is immediate to see that they form a Galois connection between ${\cal S}_{\Ind{\cal C}}$ and ${\cal G}({\cal C})$ i.e. they are adjoint to each other on the right. From the formal theory of Galois connections, it then follows that ${\cal H}({\cal K}({\cal H}({\cal D})))={\cal H}({\cal D})$ for any full subcategory $\cal D$ of $\Ind{\cal C}$ and ${\cal K}({\cal H}({\cal K}(J)))={\cal K}(J)$ for any Grothendieck cotopology on ${\cal C}$; we shall exploit this fact below.\\
The following lemma represents the extension of Lemma 3.11 \cite{prest} to the context of finitely accessible categories.

\begin{lemma}\label{lego}
Let $J$ be a Grothendieck cotopology on a small category $\cal C$. Then, with the notation above, $Hom_{\Ind{\cal C}}(-,d)$ is $J$-continuous if and only if $d\in {\cal K}(J)$, for any $d\in \Ind{\cal C}$. 
\end{lemma}
\begin{proofs}
We recall from \cite{MM} that via the equivalence ${\bf Geom}(\Set, [{\cal C}, \Set])\simeq {\bf Flat}({\cal C}^{\textrm{op}}, \Set)$ a flat functor $F:{\cal C}^{\textrm{op}}\to \Set$ is sent to the geometric morphism having as inverse image $F \otimes_{\cal C} - \cong - \otimes_{{\cal C}^\textrm{op}} F :[{\cal C}, \Set]\to \Set$; also, $F$ is $J$-continuous (for a Grothendieck topology $J$ on ${\cal C}^{\textrm{op}}$) if and only if for any $S\in J(c)$, $F \otimes_{\cal C} -$ sends the monomorphism $S\mono {\cal C}(c,-)$ to an isomorphism. Now, $Hom_{\Ind{\cal C}}(-,d):{\cal C}^{\textrm{op}}\to \Set$ is a flat functor (for any $d\in \Ind{\cal C}$), by the Yoneda representation of flat functors (cfr. \cite{El}), so it is $J$-continuous if and only if for any $S\in J(c)$, $Hom_{\Ind{\cal C}}(-,d) \otimes_{\cal C} -$ sends $S\mono {\cal C}(c,-)$ to an isomorphism. Now, given a functor $F:{\cal C}\to \Set$, $(Hom_{\Ind{\cal C}}(-,d) \otimes_{\cal C} - ) (F) = (- \otimes_{{\cal C}^{\textrm{op}}} F)(Hom_{\Ind{\cal C}}(-,d))$. If $d=colim(\overline{y}_{\cal C}\circ G)$ in $\Ind{\cal C}$ where $I$ is a filtered category and $G:I\to \cal C$ is a functor then $Hom_{\Ind{\cal C}}(-,d) \cong colim_{[{\cal C}^{\textrm{op}}, \Set]}Hom_{\Ind{\cal C}}(-,G(-))$ since all the objects in $\cal C$ are finitely presentable in $\Ind{\cal C}$ and colimits in functor categories are computed pointwise; so, since $(- \otimes_{\cal C} F)$ preserves filtered colimits (having a right adjoint) and for any $c\in {\cal C}$ $Hom_{\cal C}(-,c) \otimes_{\cal C} F \cong F(c)$ (by formula (4) p. 379 \cite{MM}), we deduce that $Hom_{\Ind{\cal C}}(-,d) \otimes_{\cal C} F\cong colim (F\circ G)=\overline{F}(d)$. Hence $Hom_{\Ind{\cal C}}(-,d) \otimes_{\cal C} -$ sends $S\mono {\cal C}(c,-)$ to the monomorphism $\overline{S}(d)\mono Hom_{\Ind{\cal C}}(\overline{y}_{\cal C}(c), d)$, from which our thesis follows.          
\end{proofs}

\begin{proposition}\label{inverso1}
Let $\mathbb T$ be a theory of presheaf type and ${\mathbb T}'$ be a geometric quotient of $\mathbb T$. Then, denoted by ${{\mathbb T}'}\textrm{-mod}(\Set)$ the full subcategory of ${{\mathbb T}}\textrm{-mod}(\Set)$ on the ${\mathbb T}'$-models, we have that ${\cal K}(\cal H({{\mathbb T}'}\textrm{-mod}(\Set)))={{\mathbb T}'}\textrm{-mod}(\Set)$. 
\end{proposition}

\begin{proofs}
By the duality theorem and Theorem 4.6 \cite{OC1}, we have that there exists a Grothendieck topology $J$ on $\textrm{f.p.} {\mathbb T}\textrm{-mod}(\Set)^{\textrm{op}}$ such that the ${\mathbb T}'$-models are exactly the $J$-homogeneous ones in any Grothendieck topos (cfr. also the proof of Theorem \ref{teoaxioms}); but by Lemma \ref{lego} a ${\mathbb T}$-model $M$ in $\Set$ is $J$-homogeneous if and only if $M\in {\cal K}(J)$, so that ${{\mathbb T}'}\textrm{-mod}(\Set)={\cal K}(J)$. Thus the thesis follows from the discussion preceding Lemma \ref{lego}.     
\end{proofs}
\begin{rmk} \emph{Conversely, we note that, by Theorem \ref{teoaxioms}, every full subcategory of ${{\mathbb T}}\textrm{-mod}(\Set)$ of the form ${\cal K}(J)$ for a Grothendieck topology $J$ on $\textrm{f.p.} {\mathbb T}\textrm{-mod}(\Set)^{\textrm{op}}$ is of the form ${{\mathbb T}'}\textrm{-mod}(\Set)$ for a geometric quotient ${\mathbb T}'$ of ${\mathbb T}$. So we conclude from Proposition \ref{inverso1} and the discussion preceding Lemma \ref{lego} that in the case of the category of models in $\Set$ of a theory of presheaf type $\mathbb T$, the `closed' full subcategories of our Galois correspondence are precisely the categories of models in $\Set$ of geometric quotients of $\mathbb T$.}
\end{rmk}

\newpage
We are now ready to prove the main result of this section.

\begin{theorem}\label{rappres}
Let $\mathbb T$ be a theory of presheaf type and ${\mathbb T}'$ be a geometric quotient of $\mathbb T$ having enough models. Then the topos $\Sh(\textrm{f.p.} {\mathbb T}\textrm{-mod}(\Set)^{\textrm{op}}, {\cal H}({{\mathbb T}'}\textrm{-mod}(\Set)))$ classifies ${\mathbb T}'$, provided that it has enough points.
\end{theorem}

\begin{proofs}
From Theorem \ref{teoaxioms} we know that there exists a geometric quotient ${\mathbb T}''$ of $\mathbb T$ such that the ${\mathbb T}''$-models are exactly the $J$-homogeneous ${\mathbb T}$-models in any Grothendieck topos. Now, ${\mathbb T}''$ has enough models, being classified by $\Sh(\textrm{f.p.} {\mathbb T}\textrm{-mod}(\Set)^{\textrm{op}}, {\cal H}({{\mathbb T}'}\textrm{-mod}(\Set)))$ which has enough points, and has the same models in $\Set$ as the theory ${\mathbb T}'$, by Lemma \ref{lego} and Proposition \ref{inverso1}. So, since both ${\mathbb T}'$ and ${\mathbb T}''$ have enough models and the same models in $\Set$, we conclude that they are syntactically equivalent and hence that they have equivalent classifying toposes; in particular $\Sh(\textrm{f.p.} {\mathbb T}\textrm{-mod}(\Set)^{\textrm{op}}, {\cal H}({{\mathbb T}'}\textrm{-mod}(\Set)))$ classifies ${\mathbb T}'$, as required. 
\end{proofs}

From the proof of the theorem, we can extract the following result.

\begin{proposition}\label{inverso2}
Let $\mathbb T$ be a theory of presheaf type and $J$ a Grothendieck topology on $\textrm{f.p.} {\mathbb T}\textrm{-mod}(\Set)^{\textrm{op}}$ such that both the toposes $\Sh(\textrm{f.p.} {\mathbb T}\textrm{-mod}(\Set)^{\textrm{op}}, J)$ and $\Sh(\textrm{f.p.} {\mathbb T}\textrm{-mod}(\Set)^{\textrm{op}}, {\cal H}({\cal K}(J)))$ have enough points. Then $J={\cal H}({\cal K}(J))$.
\end{proposition}

\begin{proofs}
By the theory of elementary toposes, $J={\cal H}({\cal K}(J))$ if and only if there exists a geometric equivalence between the toposes $\Sh(\textrm{f.p.} {\mathbb T}\textrm{-mod}(\Set)^{\textrm{op}}, J)$ and $\Sh(\textrm{f.p.} {\mathbb T}\textrm{-mod}(\Set)^{\textrm{op}}, {\cal H}({\cal K}(J)))$ which commute (in the obvious sense) with the canonical geometric inclusions $\Sh(\textrm{f.p.} {\mathbb T}\textrm{-mod}(\Set)^{\textrm{op}}, J)\hookrightarrow [\textrm{f.p.} {\mathbb T}\textrm{-mod}(\Set), \Set]$ and $\Sh(\textrm{f.p.} {\mathbb T}\textrm{-mod}(\Set)^{\textrm{op}}, {\cal H}({\cal K}(J))) \hookrightarrow [\textrm{f.p.} {\mathbb T}\textrm{-mod}(\Set), \Set]$; but this is equivalent, by the 2-dimensional Yoneda Lemma and the universal property of classifying toposes, to saying that the quotients ${\mathbb T}'$ and ${\mathbb T}''$ of $\mathbb T$ classified respectively by $\Sh(\textrm{f.p.} {\mathbb T}\textrm{-mod}(\Set)^{\textrm{op}}, J)$ and $\Sh(\textrm{f.p.} {\mathbb T}\textrm{-mod}(\Set)^{\textrm{op}}, {\cal H}({\cal K}(J)))$ via Theorem \ref{teoaxioms} have exactly the same models in any Grothendieck topos. Now, if both ${\mathbb T}'$ and ${\mathbb T}''$ have enough models then this happens precisely when they have exactly the same models in $\Set$, equivalently (by Lemma \ref{lego}) when ${\cal K}(J)={\cal K}({\cal H}({\cal K}(J)))$; but this always holds, by the discussion preceding Lemma \ref{lego}.     
\end{proofs}

\begin{rmk}\label{coherence}
\emph{Concerning the applicability of Theorem \ref{rappres}, let us mention the following fact. In \cite{flatcoh} the authors characterized the small categories $\cal C$ such that contravariant flat functors on them are coherently axiomatized in the language of presheaves on them; their condition amounts to requiring the existence of a certain kind of colimits in $\cal C$, and it is always satisfied if ${\cal C}^{\textrm{op}}$ is cartesian. Further, we note that if ${\cal C}$ satisfies this condition and $J$ is a finite type Grothendieck topology on ${\cal C}^{\textrm{op}}$ then the flat functors on ${\cal C}^{\textrm{op}}$ which are $J$-continuous can be coherently axiomatized in the language of presheaves on $\cal C$; thus the topos $\Sh({\cal C}^{\textrm{op}}, J)$ is coherent and hence has enough points by Deligne's theorem.}  
\end{rmk}

Finally, let us discuss how our results relate to those in sections 3.2 and 3.4 of \cite{prest}. There the authors only dealt with the case of embeddings ${\cal C}\hookrightarrow \Ind{\cal C}$ of the form $\textrm{f.p.}{\cal L}\hookrightarrow {\cal L}$ for a locally finitely presentable category $\cal L$. It is well-known that these categories $\cal L$ are precisely the categories of models in $\Set$ of cartesian theories, so that the category $\textrm{f.p.}{\cal L}$ always admits a syntactic description as the opposite of the syntactic category of the relevant cartesian theory; this fact is exploited in an essential way to derive some results in \cite{prest}, for example Proposition 3.5. Instead, we have arrived at Proposition \ref{inverso1}, which generalizes Proposition 3.5, by using the theory of classifying toposes and our duality theorem. Also, our Lemma \ref{lego} generalizes Lemma 3.11 \cite{prest}, whose proof relied on the locally finite presentability of the category $\cal L$, and our Proposition \ref{inverso2} implies Proposition 3.12 \cite{prest} (by Deligne's theorem, Proposition 3.4(a) \cite{prest} and Remark \ref{coherence}).

\section{A syntactic description of the finitely presented models of a cartesian theory}

In this section we give an explicit syntactic description of the finitely presented models of a given cartesian theory. We will derive this result from the well-known characterization of models of a cartesian theory as cartesian functors defined on the cartesian syntactic category of the theory. Specifically, recall from \cite{El2} (Theorem D1.4.7) that for any cartesian theory $\mathbb T$ over a signature $\Sigma$, there is an equivalence of categories $\mathfrak{Cart}({\cal C}_{\mathbb T}^{\textrm{cart}}, \Set)\simeq {\mathbb T}\textrm{-mod}(\Set)$. This equivalence is defined as follows. A cartesian functor $F:{\cal C}_{\mathbb T}^{\textrm{cart}}\to \Set$ is sent to the ${\mathbb T}$-model $F(M_{\mathbb T})$, where $M_{\mathbb T}$ is the `universal' model of $\mathbb T$ in ${\cal C}_{\mathbb T}^{\textrm{cart}}$, while a ${\mathbb T}$-model $M$ in $\Set$ is sent to the cartesian functor $F_{M}$ which sends an object $\{\vec{x}. \phi\}\in {\cal C}_{\mathbb T}^{\textrm{cart}}$ to the (domain of) its interpretation $[[\vec{x}. \phi]]_{M}$ in $M$ and an arrow $[\theta]:\{\vec{x}. \phi\}\to \{\vec{y}. \psi\}$ in ${\cal C}_{\mathbb T}^{\textrm{cart}}$ to the morphism $[[\theta]]_{M}:[[\vec{x}. \phi]]_{M}\to [[\vec{y}. \psi]]_{M}$ whose graph is the interpretation $[[\vec{x}, \vec{y}. \theta]]_{M}$ (cfr. p. 845 \cite{El2} for more details). The model $M_{\mathbb T}$ assigns to a sort $A$ the object $\{x^{A}.\top\}$ where $x^{A}$ is a variable of sort $A$, to a function symbol $f:A_{1}\cdots A_{n}\to B$ the morphism
\[  
\xymatrix {
\{x_{1}^{A_{1}}, \ldots , x_{n}^{A_{n}}. \top\} \ar[rrrr]^{[f(x_{1}^{A_{1}}, \ldots , x_{n}^{A_{n}})=y^{B}]} & & & & \{y^{B}.\top\} }
\]
and to a relation symbol $R\mono A_{1}\cdots A_{n}$ the subobject 
\[  
\xymatrix {
\{x_{1}^{A_{1}}, \ldots , x_{n}^{A_{n}}. R(x_{1}^{A_{1}}, \ldots , x_{n}^{A_{n}})\} \ar[rrrr]^{[R(x_{1}^{A_{1}}, \ldots , x_{n}^{A_{n}})]} & & & & \{x_{1}^{A_{1}}, \ldots , x_{n}^{A_{n}}. \top\} }
\]
as in Lemma D1.4.4(iv).\\
As it is remarked in \cite{El2} (Lemma D2.4.1), the $\mathbb T$-model $M_{\phi}$ corresponding via the equivalence $\mathfrak{Cart}({\cal C}_{\mathbb T}^{\textrm{cart}}, \Set)\simeq {\mathbb T}\textrm{-mod}(\Set)$ to the representable $Hom_{{\cal C}_{\mathbb T}^{\textrm{cart}}}(\{\vec{x}. \phi\}, -)\in \mathfrak{Cart}({\cal C}_{\mathbb T}^{\textrm{cart}}, \Set)$ is finitely presented by the formula $\phi(\vec{x})$. Indeed, we have the following equivalences natural in $N\in {\mathbb T}\textrm{-mod}(\Set)$:

\[
\begin{array}{ccl}
Hom_{{\mathbb T}\textrm{-mod}(\Set)}(M_{\phi}, N) & \simeq & Nat(Hom_{{\cal C}_{\mathbb T}^{\textrm{cart}}}(\{\vec{x}. \phi\}, -), F_{N})\\
& \simeq & F_{N}(\{\vec{x}.\phi\})=[[\vec{x}. \phi]]_{N},
\end{array}
\]

the second one being given by the Yoneda Lemma.\\
By recalling the definition of syntactic category ${\cal C}_{\mathbb T}^{\textrm{cart}}$, we thus obtain the following explicit description of $M_{\phi}A$.\\

$M_{\phi}$ assigns to to a sort $A$ the collection $M_{\phi}A$ of $\mathbb T$-provable equivalence classes $[\theta]$ of cartesian formulae $\theta(\vec{x}, x^{A})$ over $\Sigma$ such that the sequents $(\phi {\dashv \vdash}_{\vec{x}} (\exists x^{A})\theta)$ and $((\theta \wedge \theta[x'^{A}\slash x^{A}]) \vdash_{\vec{x}, x^{A}, x'^{A}} (x^{A}=x'^{A}))$ are provable in $\mathbb T$, where $x^{A}$ and $x'^{A}$ are distinct variables of sort $A$ not appearing in $\vec{x}$.\\  
Given a function symbol $f:A_{1}\cdots A_{n}\to B$, $M_{\phi}f:M_{\phi}A_{1}\times \cdots M_{\phi}A_{n}\to M_{\phi}B$ is the function assigning to a $n$-tuple $([\theta_{1}], \ldots, [\theta_{n}])\in M_{\phi}A_{1}\times \cdots M_{\phi}A_{n}$ the $\mathbb T$-provable equivalence class $[\exists x^{A_{1}}\ldots \exists x^{A_{n}}(\theta_{1}(\vec{x}, x^{A_{1}})\wedge \ldots \wedge \theta_{n}(\vec{x}, x^{A_{n}})\wedge y^{B}=f(x^{A_{1}}, \ldots, x^{A_{n}}))]$, where $y^{B}$ is a variable of sort $B$ not appearing in $\vec{x}$.\\
Given a relation symbol $R\mono A_{1}\cdots A_{n}$, $M_{\phi}R$ is the subset of $M_{\phi}A_{1}\times \cdots M_{\phi}A_{n}$ given by the $n$-tuples $([\theta_{1}], \ldots, [\theta_{n}])\in M_{\phi}A_{1}\times \cdots M_{\phi}A_{n}$ such that the sequent $(\theta_{1}(\vec{x}, x^{A_{1}})\wedge \ldots \wedge \theta_{n}(\vec{x}, x^{A_{n}}) \vdash_{\vec{x}, x^{A_{1}}, \ldots, x^{A_{n}}} R(x^{A_{1}}, \ldots, x^{A_{n}}))$ is provable in $\mathbb T$.\\   

Let us now verify directly that this model is finitely presented by $\phi$, by exhibiting its generators.\\
If $\vec{x}=(x_{1}^{A_{1}}, \ldots , x_{n}^{A_{n}})$ then the generators of $M_{\phi}$ are the $\mathbb T$-provable equivalence classes $\chi_{i}:=[\phi(\vec{x})\wedge x^{A_{i}}=x'^{A_{i}}]$ for $1\leq i\leq n$, where $x'^{A_{i}}$ is any variable of sort $A_{i}$ not appearing in $\vec{x}$.\\

We are now ready to describe the bijective correspondence, natural in $N\in {\mathbb T}\textrm{-mod}(\Set)$, between string of elements $(b_{1}, \ldots, b_{n})\in NA_{1}\times \ldots NA_{n}$ such that $(b_{1}, \ldots, b_{n})\in [[\phi]]_{N}$ and arrows $f:M_{\phi}\to N$ in ${\mathbb T}\textrm{-mod}(\Set)$, which witnesses the fact that $M_{\phi}$ is finitely presented by $\phi$.\\ 
To a given $\vec{a}=(a_{1}, \ldots a_{n})\in [[\phi]]_{N}$, we associate the ${\mathbb T}$-model homomorphism $f_{\vec{a}}:M_{\phi}\to N$ which assigns to each sort $A$ the function $f_{\vec{a}}A:M_{\phi}A=Hom_{{\cal C}_{\mathbb T}^{\textrm{cart}}}(\{\vec{x}. \phi\},\{x^{A}.\top\})\to NA=[[x^{A}. \top]]_{N}$ defined by $f_{\vec{a}}A([\theta])=[[\theta]]_{N}(\vec{a})$. Conversely, given a $\mathbb T$-model homomorphism $g:M_{\phi}\to N$, we associate to it the string $e_{g}:=(gA_{1}\times \ldots gA_{n})((\chi_{1}, \ldots, \chi_{n}))=(b_{1}, \ldots, b_{n})$. It is immediate to see that for any string $\vec{a}\in [[\phi]]_{N}$, $e_{f_{\vec{a}}}=\vec{a}$; to prove that $g=f_{(gA_{1}\times \ldots gA_{n})((\chi_{1}, \ldots, \chi_{n}))}$, it suffices to observe that $[[\vec{x}. \phi]]_{M_{\phi}}=Hom_{{\cal C}_{\mathbb T}^{\textrm{cart}}}(\{\vec{x}. \phi\},\{\vec{x}. \phi\})$ and then invoke the naturality of $F_{g}:F_{M_{\phi}}\imp F_{N}$ (cfr. the proof of Theorem D1.4.7 \cite{El2}).\\
It is natural to wonder how much of the preceding discussion can be adapted to regular, coherent or geometric theories. It is clear that the essential point is the fact that $Hom_{{\cal C}_{\mathbb T}^{\textrm{cart}}}(\{\vec{x}. \phi\}, -)$ is cartesian and hence corresponds via the equivalence $\mathfrak{Cart}({\cal C}_{\mathbb T}^{\textrm{cart}}, \Set)\simeq {\mathbb T}\textrm{-mod}(\Set)$ to a $\mathbb T$-model. For the above-mentioned fragments of logic, we instead have equivalences $\mathfrak{Reg}({\cal C}_{\mathbb T}^{\textrm{reg}}, \Set)\simeq {\mathbb T}\textrm{-mod}(\Set)$, $\mathfrak{Coh}({\cal C}_{\mathbb T}^{\textrm{coh}}, \Set)\simeq {\mathbb T}\textrm{-mod}(\Set)$ and $\mathfrak{Geom}({\cal C}_{\mathbb T}^{\textrm{geom}}, \Set)\simeq {\mathbb T}\textrm{-mod}(\Set)$; so, since the representables on the relevant syntactic categories are in general not regular (resp. coherent, geometric) functors we cannot conclude as above that for any regular (resp. coherent, geometric) theory there exist models which are finitely presented by given formulae in the appropriate fragments. Anyway, it is clear from the our discussion that, for any geometric theory over a signature $\Sigma$ and any geometric formula $\phi$ over $\Sigma$, there is a $\Sigma$-structure $M_{\phi}$ (in fact, a model of the cartesianization of $\mathbb T$, that is of the collection of all the cartesian formulae which are provable in $\mathbb T$) such that the $\Sigma$-structure homomorphisms $M_{\phi}\to N$ are in bijective correspondence with $[[\vec{x}. \phi]]_{N}$, naturally in $N\in {\mathbb T}\textrm{-mod}(\Set)$.
   
\section{Coherent theories and topologies of finite type}
Let us start this section with two remarks which will be important in what follows.\\

\begin{rmk}\label{provability}
\emph{Concerning the notion of provability in different fragments of logic, it is useful to remark this fact: the notion of provability in a cartesian (resp. regular, coherent) theory $\mathbb T$ with respect to cartesian (resp. regular, coherent) logic coincides with the notion of provability in $\mathbb T$ with respect to geometric logic. This can be deduced from the theory of classifying toposes as follows. As in Proposition D3.3.13 \cite{El2}, one can prove, by using the representation $[({\cal C}_{\mathbb T}^{\textrm{cart}})^{\textrm{op}}, \Set]$ (resp. $\Sh({\cal C}_{\mathbb T}^{\textrm{reg}}, J^{\textrm{reg}}_{\mathbb T})$, $\Sh({\cal C}_{\mathbb T}^{\textrm{coh}}, J^{\textrm{coh}}_{\mathbb T})$) of the classifying topos of $\mathbb T$, that the classical completeness theorem for cartesian (resp. regular, coherent) logic translates into the fact that the classifying topos $\Set[{\mathbb T}]$ of $\mathbb T$ has enough points; but this property is equivalent to the fact that $\mathbb T$, regarded as a geometric theory, has enough models (cfr. Proposition 2.3 \cite{OC5}). So all the notions of provability in question are equivalent to each other, and also to the notion of provability in classical first-order logic, being all equivalent to the notion of validity in all $\mathbb T$-models in $\Set$.}
\end{rmk}

\begin{rmk}\label{transfer}
\emph{Given two representations $\Sh({\cal C}, J)\simeq \Sh({\cal C}', J')$ of the same Grothendieck topos, we may construct a bijection between the class $\mathfrak{Groth}^{\cal C}_{J}$ of Grothendieck topologies on $\cal C$ which contain $J$ and the class $\mathfrak{Groth}^{{\cal C}'}_{J'}$ of Grothendieck topologies on ${\cal C}'$ which contain $J'$. Indeed, it is well-known that Grothendieck topologies on $\cal C$ (resp. ${\cal C}'$) which contain $J$ (resp. $J'$) are in bijection with the geometric inclusions into the topos $\Sh({\cal C}, J)$ (resp. $\Sh({\cal C}', J')$), so that we can pass from one class to the other by composing the corresponding geometric inclusions with the geometric equivalence $\Sh({\cal C}, J)\simeq \Sh({\cal C}', J')$. Moreover, via the bijections above, the natural order between geometric inclusions (i.e. one inclusion is less than another if and only if it factors through it) corresponds to the canonical order between Grothendieck topologies; thus our bijection between $\mathfrak{Groth}^{\cal C}_{J}$ and $\mathfrak{Groth}^{{\cal C}'}_{J'}$ is order-preserving and hence an Heyting algebra isomorphism. This fact will be exploited in the next section in order to obtain explicit descriptions of lattice operations between theories.}\\
\end{rmk}

Another notable application of this remark arises in the context of theories of presheaf type. Specifically, if $\mathbb T$ is a theory of presheaf type then its classifying topos can be represented either as $\Sh({\cal C}_{\mathbb T}, J_{\mathbb T})$ or as the presheaf topos $[\textrm{f.p.} {\mathbb T}\textrm{-mod}(\Set), \Set]$; thus, by the duality theorem, there is an order-preserving bijection between closed quotients ${\mathbb T}'$ of $\mathbb T$ and Grothendieck topologies $J$ on $\textrm{f.p.} {\mathbb T}\textrm{-mod}(\Set)^{\textrm{op}}$, with the property that for any ${\mathbb T}'$, the topos $\Sh(\textrm{f.p.} {\mathbb T}\textrm{-mod}(\Set)^{\textrm{op}}, J)$ of sheaves for the corresponding topology $J$ classifies ${\mathbb T}'$ (cfr. Theorem \ref{teoaxioms}).

\begin{definition}
Let $J$ be a Grothendieck topology on a category $\cal C$. Then $J$ is said to be of finite type if it is generated by a collection of finite presieves on $\cal C$.
\end{definition}
Recall that a Grothendieck topology on $\cal C$ is said to be generated by a given collection $\cal F$ of presieves on $\cal C$ if it is the smallest Grothendieck topology $J$ on $\cal C$ such that all the sieves generated by presieves in $\cal F$ are $J$-covering.

\begin{proposition}\label{finitetype}
Let $\cal C$ be a category and $J$ a Grothendieck topology on $\cal C$. Then $J$ is of finite type if and only if there exists an assignment $K$ sending to each object $c\in {\cal C}$ a collection $K(c)$ finite presieves in $\cal C$ on $c$ which satisfies the properties\\
(i) if $R\in K(c)$ then for any arrow $g:d\rightarrow c$ there exists a presieve $S\in K(c)$ such that for each arrow $f$ in $S$, $g\circ f\in R$;\\
(ii) if $\{f_{i}:c_{i}\rightarrow c \textrm{ | } i\in I \}\in K(c)$ and for each $i\in I$ we have a presieve $\{g_{ij}:d_{ij}\rightarrow c_{i} \textrm{ | } j\in I_{i} \}\in K(c_{i})$ then there exists a presieve $S\in K(c)$ such that $S\subseteq \{f_{i}\circ g_{ij}:d_{ij}\rightarrow c \textrm{ | } i\in I, j\in I_{i} \}$\\
and is such that for any sieve $S$ on $c\in {\cal C}$, $S\in J(c)$ if and only if $S\supseteq T$ for some $T\in K(c)$. 
\end{proposition}

\begin{proofs}
The `if' part of the proposition immediately follows from Definition \ref{def2}. Let us prove the `only if' part. We define $K$ as follows: for any presieve $V$ on $c\in {\cal C}$, $V\in K(c)$ if and only if $V$ is finite and the sieve generated by it is $J$-covering. By Definition \ref{def2}, $K$ satisfies properties (i) and (ii) of our proposition. Let us now define $K'$ by setting, for any sieve $R$ on $c\in {\cal C}$, $R\in K'(c)$ if and only if $R\supseteq T$ for some $T\in K(c)$. We want to prove that $J=K'$. Again, by Definition \ref{def2}, $K'$ is a Grothendieck topology and, clearly, $K'$ is contained in $J$. But the fact that $J$ is of finite type implies that $J\subseteq K'$, so that $J=K'$, as required.      
\end{proofs}

\begin{proposition}\label{finitemj}
Let $\cal C$ be a category and $J_{1}$, $J_{2}$ Grothendieck topologies on $\cal C$. Then\\
(i) If $J_{1}$, $J_{2}$ are of finite type then $J_{1}\wedge J_{2}$ is of finite type;\\
(ii) If $J_{1}$, $J_{2}$ are of finite type then $J_{1}\vee J_{2}$ is of finite type.
\end{proposition} 

\begin{proofs}
(i) Recall that for any sieve $S$ on $c\in {\cal C}$, $S\in (J_{1}\wedge J_{2})(c)$ if and only if $S\in J_{1}(c)$ and $S\in J_{2}(c)$. Let us denote by $K_{1}$ (resp. $K_{2}$) the collection of finite presieves for $J_{1}$ (resp. $J_{2}$) satisfying the conditions of Proposition \ref{finitetype}, and define $K$ as follows: for any presieve $V$ on $c\in {\cal C}$, $V\in K(c)$ if and only if there exist $V_{1}\in K_{1}(c)$ and $V_{2}\in K_{2}(c)$  such that $V=V_{1}\cup V_{2}$. Now, it is immediate to see that $K$ satisfies the conditions of Proposition \ref{finitetype}. But, clearly, for any sieve $S$ on $c\in {\cal C}$, $S\in J_{1}\wedge J_{2}(c)$ if and only if $S\supseteq T$ for some $T\in K(c)$, and hence $J_{1}\wedge J_{2}$ is of finite type by Proposition \ref{finitetype}.\\
(ii) Since $J_{1}\vee J_{2}$ is the smallest Grothendieck topology on $\cal C$ which contains both $J_{1}$ and $J_{2}$, the thesis immediately follows from the definition of Grothendieck topology of finite type; indeed, we can get a collection of finite presieves generating $J_{1}\vee J_{2}$ by taking the union of any two collections of finite presieves generating $J_{1}$ and $J_{2}$.     
\end{proofs}
Below, by a coherent theory over a signature $\Sigma$ we mean a geometric theory $\mathbb T$ over $\Sigma$ which can be axiomatized by coherent sequents over $\Sigma$.

\begin{theorem}\label{coh}
Let $\mathbb T$ be a cartesian theory over a signature $\Sigma$ and ${\cal C}^{\textrm{cart}}_{\mathbb T}$ the cartesian syntactic category of $\mathbb T$. Then the bijection between closed geometric quotients of $\mathbb T$ and Grothendieck topologies on ${\cal C}^{\textrm{cart}}_{\mathbb T}$ induced by the duality theorem via Remark \ref{transfer} restricts to a bijection between closed coherent quotients of $\mathbb T$ and finite type Grothendieck topologies on ${\cal C}^{\textrm{cart}}_{\mathbb T}$.  
\end{theorem}

\begin{proofs}
We can describe the bjection between closed geometric quotients of $\mathbb T$ and Grothendieck topologies on ${\cal C}^{\textrm{cart}}_{\mathbb T}$ induced by the duality theorem via the equivalence of classifying toposes $\Sh({\cal C}_{\mathbb T}, J_{\mathbb T})\simeq [{\cal C}^{\textrm{cart}}_{\mathbb T}, \Set]$ explicitly as follows. Given a Grothendieck topology $J$ on ${\cal C}^{\textrm{cart}}_{\mathbb T}$, the corresponding theory is axiomatized by all the sequents over $\Sigma$ of the form $\psi \vdash_{\vec{y}} \mathbin{\mathop{\textrm{\huge $\vee$}}\limits_{i\in I}}(\exists \vec{x_{i}})\theta_{i}$, where $\{[\theta_{i}] \textrm{ | } i\in I\}$ is any family of morphisms
\[  
\xymatrix {
\{\vec{x_{i}}.\phi_{i}\} \ar[r]^{[\theta_{i}]} & \{\vec{y}.\psi\}}
\]
in ${\cal C}^{\textrm{cart}}_{\mathbb T}$ forming a $J$-covering sieve. Conversely, by Proposition D1.3.10 \cite{El2}, any geometric (resp. coherent) theory over $\Sigma$ can be axiomatized by axioms of the form $\psi \vdash_{\vec{y}},  \mathbin{\mathop{\textrm{\huge $\vee$}}\limits_{i\in I}}(\exists \vec{x_{i}})\theta_{i}$ where $\psi$ and the $\theta_{i}$ are cartesian formulae over $\Sigma$ such that for any $i\in I$ $\theta_{i} \vdash_{\vec{x_{i}}, y} \psi$ is provable in geometric logic (where $I$ may be taken finite if $\mathbb T$ is coherent), so that the corresponding Grothendieck topology on ${\cal C}^{\textrm{cart}}_{\mathbb T}$ is generated by the sieves 
\[  
\xymatrix {
\{\vec{x_{i}}, \vec{y'}.\theta_{i}\} \ar[rr]^{[\theta_{i}\wedge \vec{y}=\vec{y'}]} & & \{\vec{y}.\psi\}}
\]
as $i$ varies in $I$.\\
Our thesis now follows from Remark \ref{generazione}.
\end{proofs}

\begin{corollary}
Let $\mathbb T$ be a cartesian theory over a signature $\Sigma$. Then the collection of closed coherent quotients of $\mathbb T$ form a sublattice of the collection $\mathfrak{Th}_{\Sigma}^{\mathbb T}$ of closed geometric quotients of $\mathbb T$.  
\end{corollary}

\begin{proofs}
This immediately follows from Theorem \ref{coh}, Theorem \ref{dualita} and Proposition \ref{finitemj}. 
\end{proofs}

Notice that the corollary implies that, more generally, the class of coherent theories in $\mathfrak{Th}_{\Sigma}^{\mathbb T}$ for a geometric theory $\mathbb T$ is closed under meets and joins in $\mathfrak{Th}_{\Sigma}^{\mathbb T}$; indeed, by the remark at the beginning of section \ref{latticestructure} and Remark \ref{transfer}, the meet and join of subtoposes of $\Sh({\cal C}_{\mathbb T}, J_{\mathbb T})\simeq \Sh({\cal C}_{\emptyset}, J^{\emptyset}_{\mathbb T})$ (where $\emptyset$ is the empty (cartesian) theory over $\Sigma$) are the same as those calculated in the lattice of subtoposes of $\Sh({\cal C}_{\emptyset}, J_{\emptyset})$.    

\begin{rmk}
\emph{Note that, by Remark \ref{provability}, the order-relation between coherent theories in $\mathfrak{Th}_{\Sigma}^{\mathbb T}$ is equivalent to the natural notion of order between coherent theories i.e. ${\mathbb T}_{1}\leq {\mathbb T}_{2}$ if and only if every (coherent) axiom of ${\mathbb T}_{1}$ is provable in ${\mathbb T}_{2}$ using coherent logic. Moreover, by the classical completeness theorem for coherent logic (Corollary D1.5.10 \cite{El2}), this order-relation also coincides with the well-known notion of order between first-order theories, ${\mathbb T}_{1}\leq {\mathbb T}_{2}$ being equivalent to the condition `for any $\Sigma$-structure $M$ in $\Set$, $M$ is a ${\mathbb T}_{2}$-model implies $M$ is a ${\mathbb T}_{1}$-model'.}
\end{rmk}

\section{An example}
As an application of the theory developed in the present paper, we calculate the meet of the theory of local rings and the theory of integral domains in the lattice of (coherent) theories over the signature of commutative rings with unit.\\
Let $\Sigma$ be the one-sorted signature consisting of two binary function symbols $+$ and $\cdot$, one unary function symbol $-$ and two constants $0$ and $1$, and $\mathbb T$ be the algebraic theory of commutative rings with unit over $\Sigma$; notice that the category $\textrm{f.p.} {\mathbb T}\textrm{-mod}(\Set)$ coincides with the category ${\bf Rng}_{f.g.}$ of finitely generated commutative rings with unit.\\
The theory ${\mathbb T}_{1}$ of local rings is obtained from $\mathbb T$ by adding the sequents
\[
((0=1) \vdash_{[]} \bot)
\]
and
\[
((\exists z)((x+y)\cdot z=1) \vdash_{x,y} ((\exists z)(x\cdot z=1) \vee (\exists z)(y \cdot z=1))),
\]  
while the theory ${\mathbb T}_{2}$ of integral domains is obtained from $\mathbb T$ by adding the sequents
\[
((0=1) \vdash_{[]} \bot)
\]
\[
((x \cdot y =0) \ \: \vdash_{x,y}\: ((x=0)\vee (y=0))).
\] 
Consider the Grothendieck topologies $J_{1}$ and $J_{2}$ on $\textrm{f.p.} {\mathbb T}\textrm{-mod}(\Set)^{\textrm{op}}$ corresponding respectively to ${\mathbb T}_{1}$ and to ${\mathbb T}_{2}$ as in Remark \ref{transfer}. By Example D3.1.11(a) \cite{El2} and the proof of Proposition 6.4 \cite{OC3}, we have the following descriptions:\\
for any $A\in \textrm{f.p.} {\mathbb T}\textrm{-mod}(\Set)$ and any cosieve $S$ on $A$ in $\textrm{f.p.} {\mathbb T}\textrm{-mod}(\Set)$,\\
(i) $S\in J_{1}(A)$ if and only if $S$ contains a finite family $\{\xi_{i}:A\rightarrow A[{s_{i}}^{-1}] \textrm{ | } 1\leq i\leq n\}$ of canonical inclusions $\xi_{i}:A\rightarrow A[{s_{i}}^{-1}]$ in ${\bf Rng}_{f.g.}$ where $\{s_{1},\ldots, s_{n}\}$ is any set of elements of $A$ which is not contained in any proper ideal of $A$;\\
(ii) $S\in J_{2}(A)$ if and only if either $A$ is the zero ring and $S$ is the empty sieve on it or $S$ contains a non-empty finite family $\{\pi_{a_{i}}:A\rightarrow A/(a_{i}) \textrm{ | } 1\leq i\leq n\}$ of canonical projections $\pi_{a_{i}}:A\rightarrow A/(a_{i})$ in ${\bf Rng}_{f.g.}$ where $\{a_{1},\ldots,a_{n}\}$ is any set of elements of $A$ such that $a_{1}\cdot \ldots \cdot a_{n}=0$.\\
Now, note that we may identify the polynomials with integer coefficients in a finite number of variables with $R$-equivalence classes of terms over $\Sigma$, where $R$ is the equivalence relation on terms given by `$t_{1}$ $R$ $t_{2}$ if and only if $\top \vdash t_{1}=t_{2}$ is provable in $\mathbb T$'; in fact, we shall use this identification below.\\
Then, by Theorem \ref{contruinv} and Remark \ref{transfer}, we have that ${\mathbb T}_{1}\wedge {\mathbb T}_{2}$ is obtained from $\mathbb T$ by adding the sequents
\[
((0=1) \vdash_{[]} \bot)
\]
and
\[
(\mathbin{\mathop{\textrm{\huge $\wedge$}}\limits_{1\leq s\leq m}} P_{s}(\vec{x})=0 \vdash_{\vec{x}}  \mathbin{\mathop{\textrm{\huge $\vee$}}\limits_{1\leq i\leq k}} (\exists y)(G_{i}(\vec{x})\cdot y=1) \vee \mathbin{\mathop{\textrm{\huge $\vee$}}\limits_{1\leq j\leq l}} H_{j}(\vec{x})=0)
\]
where for each $1\leq i\leq k$ and $1\leq j\leq l$, the $G_{i}$ and $H_{j}$ are polynomials in a finite string $\vec{x}$ of variables with the property that if $\vec{x}=(x_{1}, \ldots, x_{n})$ then $\{P_{1}, \ldots, P_{s}, G_{1}\ldots, G_{k}\}$ is any set of elements of ${\mathbb Z}[x_{1}, \ldots, x_{n}]$ which is not contained in any proper ideal of ${\mathbb Z}[x_{1}, \ldots, x_{n}]$ and $(\mathbin{\mathop{\textrm{ $\prod$}}\limits_{1\leq j\leq l}} H_{j})\in (P_{1}, \ldots, P_{s})$ in ${\mathbb Z}[x_{1}, \ldots, x_{n}]$.       

\vspace{10 mm}
{\bf Acknowledgements:} I am very grateful to my Ph.D. supervisor Peter Johnstone for many useful comments and discussions on the subject-matter of this paper.\\

\end{document}